Philipp Pluch

# Some Theory for the Analysis of Random Fields

## With Applications to Geostatistics

# Diplomarbeit

zur Erlangung des akademischen Grades
Diplom-Ingenieur
Studium der Technischen Mathematik

Universität Klagenfurt
Fakultät für Wirtschaftswissenschaften und Informatik
Begutachter: O.Univ.-Prof.Dr. Jürgen Pilz
Institut für Mathematik
September 2004

To my parents, Verena and all my friends

## Ehrenwörtliche Erklärung

Ich erkläre ehrenwörtlich, dass ich die vorliegende Schrift verfasst und die mit ihr unmittelbar verbundenen Arbeiten selbst durchgeführt habe. Die in der Schrift verwendete Literatur sowie das Ausmaß der mir im gesamten Arbeitsvorgang gewährten Unterstützung sind ausnahmslos angegeben. Die Schrift ist noch keiner anderen Prüfungsbehörde vorgelegt worden.

St. Urban, 29 September 2004

# Preface

I remember when I first was at our univeristy - I walked inside this large corridor called 'Aula' and had no idea what I should do, didn't know what I should study, I had interest in Psychology or Media Studies, and now I'm sitting in my office at the university, five years later, writing my final lines for my master degree theses in mathematics. A long and also hard but so beautiful way was gone, I remember at the beginning, the first mathematic courses in discrete mathematics, how difficult that was for me, the abstract thinking, my first exams and now I have finished them all, I mastered them.
I have to thank so many people and I will do so now. First I have to thank my parents, who always believed in me, who gave me financial support and who had to fight with my mood when I was working hard. My mood was not always the best, because I do not like that questions like 'Tell me about your work, what is this?' At this state of my art, if someone has a question what is my mathematical work, just read it. But without my parents, my way in the last five years would not have been possible. I also have to thank my most important person in my live, Verena, who I got to known in my third semester and she gave me so much power, the power that also made this work possible, thank you! Also thanks to my brothers Dieter and Hannes, my sister Kathrin and Claudia.
I also have to thank my dear college Samo, who was studying with me since the beginning and is now one of my best friends - maybe his ambitions to me, to help me with the solutions in the first semesters is also one part of my success. The other colleges who have gone my way with me who gave my the chance to get them known, thank you too, you all are a part of my mathematical development, specially I have to name my friends Martin, Sonja, Richi, Udo and Michi. In my way to statistics I have to tell a short story, namely when I first came in contact with that mathematical theme I was 17 years old in the 'Oberstufe' with my teacher Schicher Ingrid, I always was good in mathematics, but that thing, I didn't like it, I was happy at my last exam on that theme 'Never again' and now that part I have done in an extreme way - no random event, now it is a random field!

One special thank is to Dr. Jürgen Pilz - for all his enthusiasm in his hour long talks in his office, this was the basic thing why I started loving that field of mathematics. I tried to read every book that he told me in the hope that I can discuss with him, but so often after the first twenty pages I didn't understand anything, but we had fruitful discussion. I hope that my academical father, Dr. Jürgen Pilz, will care long long time over my development. When speaking from academical parents I have also to name Dr. Christine Nowak - under students known as the most strict professor, she and her Analysis is one thing I will never forget, the proofs on the blackboard but also supporting talks with her when I was absolutely done, that always pushed my up. I hope that in future there will be a lot of mountains for climbing together. Dr. Pilz showed me how interesting mathematics can be, he is a living reference book and knows to everything anything and Dr. Nowak learned my how I work mathematically correct. Thank you.
I also have to thank Mag. Gunter Spöck who is working with me in a project, he is one person who, how I should say, finds an error in a code by just looking one minute, he is also one part of that wonderful cooking incidents of that work, the long discussions and the learning how to code right were very wonderful, also thanks for all the coffee! Beside Gunter I must name Dr. Werner Müller who is the project leader, Dr. Milan Stehlik the sunny boy at moDa and also a part of the project, thank you for supporting me. Also lots of thanks to Dr. Rose-Gerd Koboltschnig, Dr. Dieter Rasch and Dr. Albrecht Gebhardt the linux ghost. At the end, thanks also to all my friends in the outer mathematical world, that they learned to understand that last years there was less snowboarding, biking and climbing, thanks that you also stayed until now on my side, namely Mike, my brother Hannes und even new Moritz. A lot of thanks also to my friends in the Lidmanskygasse: Stefan, Anja and Johann!
After all that words of thanks, I hope that I have not forgotten any one, but if so, I will make a webpage containing comments and corrections regarding this thesis.

Klagenfurt, September 2004 *Pluch Philipp*

# Contents

# 1

# The Spectrum of Random Fields

## 1.1 Application of Random Fields

Many biological, physical and even social systems have attributes, which viewed on an appropriate scale, exhibit complex patterns of variation in space and/or time. For example substances that constitute the earth's environment possess physical properties which change more or less unpredictably under the influence of meteorological variables. These variables are themselves random processes in space and/or time. For example these are temperature, pressure or fishes in a lake.
When the degree of disorder is sufficiently large, there is usually a merit and economy in probabilistic rather than deterministic models. Random field theory seeks to model complex patterns of variation and interdependence in cases where deterministic treatments are inefficient, and conventional statistics are insufficient. An ideal random field model will capture the essential features of a complex random phenomenon in terms of a minimum number of physically meaningful and experimentally accessible parameters (see VanMarke (1988)).
The methodology of random fields is applicable to phenomena occurring on different temporal or spatial scales. In the case of time, this scale may be the interval between molecular collisions, as in the study of Brownian motion (firstly studied by Norbert Wiener(1938)) , or it may be measured in geological units, for example, to describe the variation of properties and thickness of a layer on a wafer in semi conductor industry. Similarly on a spatial coordinate axis the scale may be subatomic in the study of super heated plasma, while an astronomical scale is needed to describe the property (temperature, density or chemical composition) of matter in interstellar space (see Hawking (1993)).
We can speculate that many phenomena can be examined in the light of random field theory. Each such phenomenon summons the image of a so called distributed disordered system whose attributes display a complex pattern of variation in space. We will handle with dimension two or higher, as well as variation in time. As an example we observed radioactivity at a point with nor-thing, easting, altitude and time. Random field models of complex stochastic phenomena serve multiple purposes:

- particular that in a stochastic description they provide a format for efficient characterization of distributed disordered system
- in a system analysis they provide the basis for predicting response of a distributed disordered system
- in decision analysis they permit assessment of the impact of alternative strategies in decision situations involving distributed disordered systems, as in data acquisition, on line control or design

See VanMarke(1988). Differences between types of random fields stem mainly from the nature of the decisions that the stochastic environment imposes or creates. Uncertainty about the properties of a random medium is essentially of a passive type. An measurement at a given location is deterministic, but its value at a given location is unknown until anybody sees what is measured. Sampling at any location is usually impractical and in field experiments too expensive, and measurement and testing errors tend to dilute the value of the information, see Mueller (2001). The prediction, analysis and decision making must usually proceed on the basis of incomplete information. A spatial process is characterized by active and inherent (or intrinsic) uncertainty, this means that properties at different points in space change randomly to an other point in space. The measurements may be used for forecasting the future values of system wide averages or measures of performance.
If one associates a random variable with each point of a lattice with coordinates $x = (x_1, x_2, ..., x_n)$ in a $n-$dimensional space, the space of these random variables is called a random field. The interdependence among the random variables of such a field will be significant in determining the properties of macroscopic quantities of such a random field.

## 1.2 Random Fields

In the following section we will work out the basic mathematical properties for using a random field in practice. First concepts like set theory (balls, spheres and distance) are introduced. We are going to define the scalar product of two vectors and then introduce the norm of a vector. We will have a closer look at random variables and define them in a pure mathematical way, by developing the $L_2$ we are able to see that random variables in $L_2$ have special properties, they have finite first and second order moments. A main part in this chapter will be the second order properties of random elements, characterised mainly by the scalar product in $L_2$. Through the scalar product we get the norm, with this norm we can characterise what it means that a sequence $X_n$ of random variables converges to a random variable $X$. We will give analytic properties of the Fourier transformation and present the Riemann Lebesgue Lemma, which is also known as the theorem of Mercer. One interesting result in analysis is the Plancherel theorem, which is related to the Fourier transformation. This theorem is the main result in this chapter and it has consequences for Bochner's theorem, for the inverse Fourier transformation.

### 1.2.1 Elements of the Theory of Random Fields

#### Basic Concepts

In this section we will present the basic concepts and statements that are important for developing the theory of random fields. The following statements will have an important influence on the applications in this work.

#### Some Set Theory

Let $\mathbb{R}^n$ be the real Euclidean vector space of $dim(\mathbb{R}^n) = n \geq 1$. On that vector space the scalar product of two vectors $x = (x_1, ...x_n)$ and $y = (y_1, ..., y_n) \in \mathbb{R}^n$ is defined as

$$\langle x, y \rangle = x_1 y_1 + ... + x_n y_n.$$

With the concept of the scalar product we will introduce the length of an vector $x$ as

$$|x| = \langle x, x \rangle^{1/2}$$

and $\rho_{x,y}$ will denote the distance in the above sense between $x$ and $y$. The symbol $dx = dx_1, ..., dx_n$ denotes an element of the Lebesgue measure in $\mathbb{R}^n$. By $\mathbb{R}^n_+$ we will denote the positive octant of $\mathbb{R}^n$ defined as

$$\mathbb{R}^n_+ = \{x \in \mathbb{R}^n : x_i \geq 0, i = 1, ..., n\}$$

and by $\mathbb{Z}^n, n \geq 1$ the integer lattice in $\mathbb{R}^n$.
The sets

$$v_n(r) = \{x \in \mathbb{R}^n : |x| < r\}$$
$$s_n(r) = \{x \in \mathbb{R}^n : |x| = r\}$$

define a ball and a sphere in $\mathbb{R}^n$ of radius $r$ and the centre at the origin. If it is clear which dimension the underlying space has, we will just write $v(r)$ or $s(r)$.
A parallelepiped in $\mathbb{R}^n$ is defined by

$$\Pi[a, b) = \{x \in \mathbb{R}^n : a_i \leq x_i < b_i, i = 1, ..., n\}$$

and in a analog way $\Pi[a,b], \Pi(a,b)$ and $\Pi(a,b]$ are defined and we will write $\Pi(c)$ when $\Pi[0,c]$ is mentioned.
An element of the Lebesgue measure of $s_{n-1}(r)$ will be denoted by

$$dm(x) = r^{n-1}\sin^{n-2}\phi_{n-2}...\sin\phi_1 d\phi_1...d\phi_{n-2}d\phi$$

where the spherical coordinate system $(r,u)$ in $\mathbb{R}^n$ has been used:

$$r \geq 0, u = (\phi, \phi_1, ..., \phi_{n-2}) \in s_{n-1},\ 0 \leq \phi \leq 2\pi$$

and

$$0 \leq \phi_k \leq \pi \text{ for } k = 1, ..., n-2, n \geq 2$$

which is related to the Cartesian coordinate system by the common way.

$$\begin{aligned} x_1 &= r\sin\phi_{n-2}...\sin\phi_1\sin\phi \\ x_2 &= r\sin\phi_{n-2}...\sin\phi_1\cos\phi \\ &... \\ x_{n-1} &= r\sin\phi_{n-2}\cos\phi_{n-3} \\ x_n &= r\cos\phi_{n-2} \end{aligned}$$

More about the Spherical and Cartesian coordinate system can be found in Walter(2002).
We will denote the Lebesgue measure of a measurable set $A \in \mathbb{R}^n$ by $|A|$.
Let $a \in \mathbb{R}^n$ be a non-zero vector and $\Pi_m = \Pi[0,a) + m \otimes a$ the shift of $\Pi[0,a)$ by the vector $m \otimes a = (m_1 a_1, ..., m_n a_n)\, m \in \mathbb{Z}^n$. It is clear that the set of $\{\Pi_m\}, m \in \mathbb{Z}^n$ generates a partition of $\mathbb{R}^n$.
The distance between sets and set and point are defined in the general way by

$$\begin{aligned} \rho(A,B) &= \inf\{|x-y| : x \in A, y \in B\} \\ \rho(x,A) &= \inf\{|x-y| : y \in A\} \end{aligned}$$

## Random Variable

Let $(\Omega, \mathcal{A}, P)$ be a complete probability space, then a random element taking values in a measurable space $(X, \mathcal{B})$ is a mapping $\xi : \Omega \to X$ such that $\{\omega : \xi(\omega) \in A\} \in \mathcal{A}$ for any set $A \in \mathcal{B}$. The measure $P_\xi(A) = P(\xi \in A)$ defined on the $\sigma$-algebra $\mathcal{B}$ is called the distribution of $\xi$. A random element taking values in $(\mathbb{R}^m, \mathcal{B}^m)$ is called a random variable if $m = 1$ and a random vector if $m \geq 2$.
To define $P_\xi(A)$ for a random vector $\xi = (\xi_1, ..., \xi_n) \in \mathbb{R}^n$ it suffices to define the distribution function $F(x) = F_\xi(x_1, ..., x_n) = P(\xi_1 < x_1, ..., \xi_n < x_n)$. A distribution $P_\xi(A)$ is said to be absolutely continuous if

$$P_\xi(A) = \int_A p(x)dx \, , A \in \mathcal{B}^n$$

The function $p(x) = p_\xi(x_1, ..., x_n)\, , x \in \mathrm{I\!R}^n$ is called the density function of the distribution of the random variable $\xi \in \mathrm{I\!R}^n$.
In the following we will use Latin letters instead of Greek letters, if its clear that the notation means a random variable or not. A sequence of random variables $X_m \, , m = 1, 2, ...$ converges in probability to the random variable $X$, if for any $\varepsilon > 0$ $P(|X_m - X| > 0) \to 0$ as $m \to \infty$. A sequence of random variables $X_m \, , m = 1, 2, ...$ converges almost surely (a.s) to a random variable $X$ if $P(X_m \to X) = 1$ for $m \to \infty$. In order for a sequence of random variables $X_m$ to converge to the random variable $X$, it is sufficient to show that the series

$$\sum_{m=1}^{\infty} P(|X_m - X| > \varepsilon)$$

converges for any $\varepsilon > 0$.
Random variables $X_1, ..., X_m$ are said to be independent if $P(X_1 \in B_1, ..., X_m \in B_m) = P(X_1 \in B_1)...P(X_m \in B_m)$ for any $B_i \in \mathcal{B}^1$, $i = 1, ..., m$.

### 1.2.2 Expected Value

Suppose $X$ is a non-negative random variable in the probability space $(\Omega, \mathcal{A}, P)$ and define the integral of $X$ with respect to $P$ by means of

$$\int X \, dP = \sup_{S \in \Sigma} \int S \, dP \tag{1.1}$$

where $\Sigma$ consists of all non-negative simple random variables $S$ such that $S \leq X$ a.s; we know also that $\Sigma \neq 0$. Again $\int X \, dP$ is a non-negative number, possibly $+\infty$.
Following inequality holds

$$\int X \, dP \leq \int Y \, dP \text{ if } X \leq Y \tag{1.2}$$

For a general random variable $X$, $|X|$ is a non-negative random variable, and so are $X^+$ and $X^-$, defined as

$$\begin{aligned} X^+ &= \frac{|X| + X}{2} \\ X^- &= \frac{|X| - X}{2} \\ X &= X^+ - X^- \\ |X| &= X^+ + X^- \end{aligned}$$

Then the integral

$$\int X^+ \, dP \text{ and } \int X^- \, dP$$

is always well defined, although the value $+\infty$ can be assumed.

**Definition 1.1.** *The integral of a real random variable $X$ with respect to the probability $P$ is defined by*

$$\int X \, dP = \int X^+ \, dP - \int X^- \, dP \tag{1.3}$$

*provided at least one of the two integrals on the right hand side is finite.*

**Definition 1.2.** *A random variable X has an expectation EX if $\int X(\omega)\, dP(\omega)$ is finite. The expected value of $X$ is then*

$$EX = \int_{\Omega} X(\omega)\, dP(\omega) = \int_{\mathbb{R}^n} x d\mu_X(x) \tag{1.4}$$

*w.r.t. $P$. $\mu_X$ is a measure, a countable additive nonnegative set function defined on the sigma algebra.*

More generally, if $f : \mathbb{R}^n \to \mathbb{R}$ is Borel measurable and $\int_{\Omega} |f(X(\omega))| dP(\omega) < \infty$ then we have

$$E[f(X)] = \int_{\Omega} f(X(\omega)) dP(\omega) = \int_{\mathbb{R}^n} f(x) d\mu_X(x).$$

Clearly $X$ has expectation if and only if $|X|$ has an expectation. If X has a finite expected value, we will write $X \in L_1(\Omega, A, P)$. Generally we can see, that it is true that

$$E\phi(X) = \int_{-\infty}^{+\infty} \phi(y) dF_x(y). \tag{1.5}$$

**Lemma 1.3.** *$X \in L_1 \Leftrightarrow |X| \in L_1$ and $|EX| \leq E|X|$*

**Theorem 1.4.** *Let $X, Y \in L_1, \alpha, \beta \in \mathbb{R}$ and $\phi$ a convex function, then*

$$\alpha X + \beta Y \in L_1 \text{ and } E(\alpha X + \beta Y) = \alpha EX + \beta EY$$
$$\phi(EX) \leq E(\phi \circ X) \text{ (Jenssen's inequality).}$$

*holds, in particular if $X \in L_1$ then*

$$E(\sigma X + \mu) = \sigma(EX) + \mu.$$

We set by definition $E(X) = (E(X_1), ..., E(X_n))$.

### 1.2.3 Variance

Jenssen's inequality has numerous consequences in the theory, for example taking $\phi(x) = x^2$ we obtain

$$(EX)^2 \leq EX^2$$

if $X^2 \in L_1$. In this case we say that the random variable X is square summable and we will write $X \in L_2(\Omega, A, P)$. That is we can define following set:

$$L_2 = \{X \in L_1 : |X|^2 \in L_1\}$$

**Lemma 1.5.** *For a random variable $X$ with $EX^2 < \infty$ the Chebyshev's inequality hods: for any $\varepsilon > 0$*

$$P(|X - EX| > \varepsilon) \leq \frac{VarX}{\varepsilon^2}$$

If $X$ is square summable then following equation holds

$$EX^2 = \int_{-\infty}^{\infty} y^2 dF_X(y). \tag{1.6}$$

**Definition 1.6.** *If $X \in L_2$, we define its variance by*

$$VarX = E(X - EX)^2.$$

In fact the difference $X - EX$ is a new random variable. Since

$$(X - EX)^2 = X^2 - 2(EX)X + (EX)^2$$

it follows that

$$E(X - EX)^2 = EX^2 - 2(EX)^2 + (EX)^2 = EX^2 - (EX)^2.$$

### 1.2.4 Stochastic Process

**Definition 1.7.** *A stochastic process is a parametrised collection of random variables*

$$\{X_t\}_{t \in T}$$

*defined on a probability space $(\Omega, \mathcal{A}, P)$ and assuming values in $\mathbb{R}^n$.*

The parameter space $T$ is a subspace of $\mathbb{R}^n$ for $n \geq 1$. Note that for each $t \in T$ fixed, we have a random variable

$$\omega \to X_t(\omega); \ \omega \in \Omega$$

On the other hand, fixing $\omega \in \Omega$ we can consider the function

$$t \to X_t(\omega); \ t \in T$$

which is called a path of $X_t$. The (finite-dimensional) distribution of a process $X = \{X_t\}_{t \in T}$ are the measures $\mu_{t_1, t_2, ..., t_k}$ defined on $\mathbb{R}^{nk}, k = 1, 2, ...$ by

$$\mu_{t_1, ..., t_k}(F_1 \times ... \times F_k) = P[X_{t_1} \in F_1, ..., X_{t_k} \in F_k]; \ t_i \in T$$

Here $F_1, ..., F_k$ denote Borel sets in $\mathbb{R}^n$.

### 1.2.5 Second Order Properties

#### Orthogonality

Let $(\Omega, \mathcal{A}, P)$ be a fix probability space, and let $X, Y, Z, ...$ denote complex random variables. All the random variables encountered will be assumed to have second order moments, this means $X, Y, Z, ... \in L_2(\Omega, \mathcal{A}, P)$.

**Lemma 1.8.** $X, Y \in L_2 \Rightarrow X + Y \in L_2$

This follows from the inequality

$$|a+b|^2 \leq 2(|a|^2 + |b|^2) \text{ for all } a, b \in \mathbb{C}$$

Moreover, clearly

$$X \in L_2 \Rightarrow \alpha X \in L_2 \ \forall \alpha \in \mathbb{R}$$

We can see, that $L_2$ is a complex vector space.

We define a new operation on $L_2$, namely an inner product on $L_2$ given by

$$\langle X, Y \rangle = E(X\overline{Y}) \tag{1.7}$$

The inner product is well defined, because of the inequality $2|ab| \leq |a|^2 + |b|^2$ holds for all $a, b \in \mathbb{C}$, i.e. $\langle X, Y \rangle$ is a complex number for any random variables $X, Y \in L_2$. It is easy to show that $\langle ., . \rangle$ possesses the following properties:

- Bilinearity:

$$\langle \alpha X + \beta Y, Z \rangle = \alpha \langle X, Z \rangle + \beta \langle Y, Z \rangle$$
$$\langle X, \alpha Y + \beta Z \rangle = \overline{\alpha} \langle X, Y \rangle + \overline{\beta} \langle X, Z \rangle$$

- Hermitian character

$$\langle X, Y \rangle = \overline{\langle Y, X \rangle}$$

- Positive non degeneracy

$$\langle X, X \rangle \geq 0 \text{ and } \langle X, X \rangle = 0 \Leftrightarrow X = 0$$

**Definition 1.9.** *We call two random variables $X, Y \in L_2$ orthogonal if $\langle X, Y \rangle = 0$*

A consequence of orthogonality for random variables is the following: Any finite collection of pairwise orthogonal random variables $\{X_i\}$ in $L_2$ are linearly independent. In general it is not true that linearly independent elements are necessarily orthogonal.
With the help of the inner product on $L_2$ it is now also possible to introduce the standard norm on $L_2$ by

**Definition 1.10.** *The norm of $X \in L_2$ is defined by*

$$\|X\| = \langle X, \bar{X} \rangle^{1/2}$$

This definition makes sense for any $X$ in $L_2$ by the positive character of the inner product $\langle ., . \rangle$. Clearly $\|X\| \geq 0$, and because of the non degeneracy of $\langle ., . \rangle$ there results $\|X\| = 0 \Leftrightarrow X = 0$. By the bilinearity of the inner product $\langle ., . \rangle$ it follows that $\|\alpha X\| = |\alpha| \|X\|$. The triangle inequality can be obtained directly from the Schwarzs inequality.
The advantage of having a norm $\|.\|$ on a linear space $L_2$ is, that we now are able to define a metric on $L_2$. Given $X, Y \in L_2$ the distance between them is defined by

$$d(X, Y) = \|X - Y\|.$$

The distance $d : L_2 \times L_2 \to \mathbb{R}$ possesses the following properties for $X, Y, Z \in L_2$, namely

- Symmetry

$$d(X, Y) = d(Y, X)$$

- Positive non degeneracy

$$d(X, Y) \geq 0, d(X, Y) = 0 \Leftrightarrow X = Y$$

- Triangle inequality

$$d(X, Z) \leq d(X, Y) + d(Y, Z)$$

A metric in $L_2$ allows us now to introduce convergence, let $\{X_n\}$ be a given sequence in $L_2$ and $X \in L_2$, we will say $X_n$ converges to $X$, and write $X_n \to X$ if $d(X_n, X) \to 0$, when $n \to \infty$. In other words, $X_n \to X$ if and only if given $\epsilon > 0$ there exists $N \geq 1$ such that

$$\|X_n - X\| < \epsilon \text{ if } n \geq N$$

This type of convergence of random variables is in literature known as convergence in quadratic mean. We will write for that $X = \text{l.i.m}_{n \to \infty} X_n$.
The space $L_2$ is a complex vector space with a inner product defined in terms of which both a norm and the convergence can be defined, moreover the resulting space is complete. Thus, $L_2$ is a Hilbert space.

### 1.2.6 Fourier Transformation

In this section we will show some purely mathematical developments which we will find in important applications, like the Bochner theorem for spectral representation of covariance functions and random fields.

### Characteristic Functions

Let $X$ be a real random vector and consider the complex valued random variable $e^{i\langle\xi,X\rangle}$, for fixed $\xi \in \mathbb{R}^q$. This new random variable is bounded and therefore belongs to $L_1$. Define

$$\phi_X(\xi) = E(e^{i\langle\xi,X\rangle}).$$

The complex valued function thus defined is called the characteristic function of $X$, and is defined for any $\xi \in \mathbb{R}^q$. Explicitly

$$\phi_X(\xi) = \int_{\mathbb{R}^q} e^{i\langle\xi,X\rangle} dF_X(x)$$

and if $X$ has a density,

$$\phi_X(\xi) = \int_{\mathbb{R}^q} e^{i\langle\xi,X\rangle} f_X(x)dx$$

By the definition of $\phi_X(\xi) = E(e^{i\langle\xi,X\rangle})$ it follows that

- $\phi_X(0) = 1$
- $|\phi_X(\xi)| \leq 1$ for all $\xi \in \mathbb{R}^q$
- $\phi_{\langle a,X\rangle+b}(\xi) = e^{i\langle b,\xi\rangle}\phi_X(\langle a,\xi\rangle)$

A random vector $X \in \mathbb{R}^q, q > 1$ is called Gaussian if its characteristic function is of the form

$$\phi_X(t) = e^{i\langle t,a\rangle - \frac{1}{2}\langle Bt,t\rangle}$$

$t \in \mathbb{R}^q, a \in \mathbb{R}^q$, the mean and $B$ is a linear self-adjoined nonnegative definite operator called the correlation operator. The Matrix $B$ defined by the operator $B$ is called correlation matrix, given by

$$B = (\mathrm{cov}(X_i, X_j))_{i,j=1,\dots,q}$$

a nonnegative definite symmetric matrix. If $B > 0$, the random vector $X$ has the density function

$$f_{a,B}(u) = (2\pi)^{-q/2}(\det B)^{1/2} e^{-\frac{1}{2}\langle B^{-1}(u-a),u-a\rangle}$$

with $a = (a_1, ..., a_q)$, $a_i = EX_i, i = 1, ..., q$ and $B$ defined by the operator $B$ as above.
Characteristic functions are necessarily continuous, which can easily be seen in the following argument

$$\begin{aligned}|\phi_X(\xi) - \phi_X(\xi + h)| &= |E(e^{i\langle\xi,X\rangle} - e^{i\langle\xi+h,X\rangle})| \\ &= E|e^{i\langle\xi,X\rangle}(1 - e^{i\langle h,X\rangle})| \leq E|e^{i\langle h,X\rangle} - 1|\end{aligned}$$

Clearly $|e^{i\langle h,X\rangle} - 1| < 2$ and $e^{i\langle h,X\rangle} \to 1$ when $h \to 0$ and therefore $E|e^{i\langle h,X\rangle} - 1| \to 0$ by the theorem of dominated convergence. Indeed characteristic functions are uniformly continuous. The regularity of a given random variable is mirrored by the smoothness properties of the corresponding characteristic function. This fact is mentioned in the following theorem:

**Theorem 1.11.** *If $X \in L_n$, then $\phi_X \in C^n$ and*

$$\phi_X^n(\xi) = i^n E(X^n e^{i\langle \xi,X\rangle})$$

An additional interesting result for characteristic functions in the connection with distribution function gives following theorem

**Theorem 1.12.** *A sequence $\{F_n\}$ of distribution functions converges to a distribution function $F$ if and only if the sequence of their characteristic functions $\{\phi_n\}$ converges to a continuous limit $\phi$. In this case, $\phi$ is the characteristic function of $F$, and the convergence $\phi_n \to \phi$ is uniform in every compact interval $I$.*

### Fourier Transformation

In the following we will consider that all functions are complex valued, defined on a Euclidean space $\mathbb{R}^n$, the integration will be understood with respect to the $n$-dimensional Lebesgue measure $dx$. For that we will make the following definitions

**Definition 1.13.** *A function $f : \mathbb{R}^n \to \mathbb{C}$ will be called integrable if and only if $\int |f(x)|dx < \infty$. The notation for integrable functions will be $f \in L_1$*

**Definition 1.14.** *A function $f : \mathbb{R}^n \to \mathbb{C}$ will be called square integrable if and only if $|f|^2 \in L_1$, in symbols $f \in L_2$*

**Lemma 1.15.** *$f \in L_\infty$ if and only if $f$ is bounded on $\mathbb{R}^n$*

Fix $\xi \in \mathbb{R}^n$ and let $f \in L_1$, then $|e^{-2\pi i x} f(x)| = |f(x)|$.

**Definition 1.16.** *The function*

$$\hat{f}(\xi) = \int e^{-2\pi i\langle \xi,x\rangle} f(x)dx \tag{1.8}$$

*will be called the Fourier Transform of f.*

Clearly $|\hat{f}(x)| \le \int |f(x)|dx$ for all $f \in L_1$. Another property will be that $\hat{f} \in L_\infty$ for each $f \in L_1$ and the mapping $f \mapsto \hat{f}$ is a bounded linear operator from $L_1$ into $L_\infty$ with

$$||\hat{f}||_\infty = \sup_{x\in\mathbb{R}^n} |f(x)| \le ||f||_1 = \int |f(x)|dx$$

The inner product in $L_2$ is

$$\langle f, g\rangle = \int f(x)\overline{g(x)}dx$$

and the norm is given in the usual way

$$||f||_2 = \langle f, f\rangle^{1/2}$$

Note that for each $f \in L_1$, the Fourier Transformation $\hat{f}$ is uniformly continuous.
An interesting question will be, which complex continuous function can be the Fourier Transform of an integrable function.

**Theorem 1.17 (Riemann-Lebesgue's Lemma).** *For each* $f \in L_1, \hat{f}(\xi) \rightarrow 0$ *as* $|\xi| \rightarrow \infty$

PROOF: This result can be easily established for special choices of $f$. For instance, let $Q \equiv \Pi[a, b]$ be a parallelepiped in $\mathbb{R}^n$, and let $f(x) = 1$ on $Q$, $f(x) = 0$ if $x \notin Q$. Then

$$\begin{aligned}\hat{f}(\xi) &= \int_a^b ... \int_a^b e^{-2\pi i(\xi_1 x_1 + ... + \xi_n x_n)} dx_1 .... dx_n \\ &= \frac{i}{2\pi\xi_1} e^{-2\pi i(\xi_1 x_1)}|_{a_1}^{b_1} \cdots \frac{i}{2\pi\xi_n} e^{-2\pi i(\xi_n x_n)}|_{a_n}^{b_n}\end{aligned}$$

i.e.

$$|\hat{f}(x)| \leq \frac{2^n}{(2\pi)^n |\xi_1| ... |\xi_n|} \rightarrow 0 \text{ as } |\xi| \rightarrow \infty$$

If $Q_1, ..., Q_k$ are parallelepipeds in $\mathbb{R}^n$ and $f = \sum_{k=1}^n \alpha_k f_k$, with $f_k(x) = 1$ on $Q_k$, $f_k(x) = 0$ outside $Q_k$, then $\hat{f} = \sum_{k=1}^n \alpha_k \hat{f}_k$ and therefore $\hat{f}(\xi) \rightarrow 0$ as $|\xi| \rightarrow \infty$.
A generic $f \in L_1$ can be approximated arbitrarily well in the sense of $||.||_1$ by a sequence of functions of the form considered above. That is, for $f \in L_1$ and each $\varepsilon > 0$ there is a function $g$ of the form, that $g = \sum_{k=1}^n \alpha_k f_k$ such that $||f - g|| \leq \varepsilon/2$, and $|\hat{g}(x)| \rightarrow 0$ as $|\xi| \rightarrow \infty$. Therefore

$$\begin{aligned}|\hat{f}(\xi)| &\leq ||\hat{f} - \hat{g}||_\infty + |\hat{g}(\xi)| \\ &\leq ||f - g||_1 + |\hat{g}(\xi)| \\ &< \frac{\varepsilon}{2} + \frac{\varepsilon}{2} \\ &= \varepsilon\end{aligned}$$

if $|\xi| > M$. Therefore $|\hat{f}(\xi)| \rightarrow 0$ as $|\xi| \rightarrow \infty$.

□

By the Riemann-Lebesgue Lemma follows that not every complex continuous function can be the Fourier Transform of an integrable function $f$: if it does not vanish asymptotically for large values of the argument it is not a Fourier transform. *A constant function is not the Fourier transformation of an integrable function.*
Suppose both $f$ and $x_k f(x)$ are integrable, let $e^{(k)}$ denote the $k$-th unit vector in the usual canonical basis of $\mathbb{R}^n$ then

$$\frac{\hat{f}(\xi + he^{(k)}) - \hat{f}(\xi)}{h} = \int e^{-2\pi i \xi x}(\frac{e^{-2\pi i h x_k} - 1}{h})f(x)dx$$

for all real numbers $h \neq 0$. Notice that

$$\frac{e^{-2\pi i h x_k} - 1}{h} \to -2\pi i x_k \text{ as } h \to 0.$$

By Dominated Convergence

$$\frac{\hat{f}(\xi + he^{(k)}) - \hat{f}(\xi)}{h} \to \int e^{-2\pi i \xi x}(-2\pi i x_k f(x))dx$$

Therefore if $f \in L_1$ and $x_k f(x)$ is also integrable, then $\frac{\partial f}{\partial x_k}$ is the Fourier transform of $-2\pi i x_k f(x)$. Thus, better integrability properties of $f$ give rise to greater regularity of $\hat{f}$.

- If $f \in L_1$ then $\hat{f} \in C_0$ and

$$||\hat{f}||_\infty \leq ||f||_1$$

- If $f \in L_1$, $\frac{\partial f}{\partial \xi_k} \in L_1$ and $f(x) \to 0$ when $|x| \to \infty$, then

$$\frac{\partial \hat{f}}{\partial x_k}(\xi) = 2\pi i \xi_k \hat{f}(\xi)$$

By Fubini's theorem, the function

$$x \mapsto \int f(y)g(x-y)dy$$

is integrable.

**Definition 1.18.** *The convolution of $f$ and $g$ is defined as:*

$$f * g(x) = \int f(y)g(x-y)dy$$

By Fubini's theorem we can infer that:

$$f * g = g * f$$
$$(f * g) * h = f * (g * h)$$

if and only if $f, g, h \in L_1$.
Applying Fubini's theorem

$$\widehat{f * g}(\xi) = \int_{\mathbb{R}^n} \int_{\mathbb{R}^n} e^{-2\pi i \langle \xi, x \rangle} f(y) g(x-y) dx dy$$

and the transformation of $\mathbb{R}^n \times \mathbb{R}^n$

$$u = y, \; v = x - y$$

we get

$$\widehat{f * g}(\xi) = \int_{\mathbb{R}^n} e^{-2\pi i \langle \xi, u \rangle} f(u) du \int_{\mathbb{R}^n} e^{-2\pi i \langle \xi, v \rangle} g(v) dv$$

this means that

- If $f, g \in L_1$ then

$$\widehat{f * g}(\xi) = \hat{f}(\xi) \hat{g}(\xi)$$

**Plancherel Theorem**

Since the Lebesgue measure of $\mathbb{R}$ is infinite, $L_2$ is not a subset of $L_1$, and the definition of the Fourier transformation is therefore not appropriate for every function $f$ from the function space $L_2$. The definition does only apply, if the function is in $L_1 \cap L_2$, then we can follow that $\hat{f} \in L_2$. We can see, that the norms of the Fourier transformation of the function and the norm of the function in the sense of $||.||_2$ coincide. This isometry of $L_1 \cap L_2$ into $L_2$ extends to an isometry of $L_2$ onto $L_2$, and this extension defines the Fourier transformation of every $f \in L_2$, this transformation is also known by the Plancherel transformation. It is a result in harmonic analysis, first proved by Michel Plancherel, a Swiss mathematician.

We now consider square summable functions and we will develop the Fourier inversion theory for functions in $L_2$, an important property of the $L_2$ is that square summable functions may not be necessarily integrable, but they can be approximated in the mean by square integrable functions which are integrable. An other important fact is, that

$$L_1 \cap L_2 \text{ is dense in } L_2$$

An important consequence is the following Lemma:

**Lemma 1.19.** *If $f \in L_1 \cap L_2$ then $\hat{f} \in L_2$ and $||\hat{f}||_2 = ||f||_2$*

A proof of that Lemma can be found in (Stein and Weiss 1971). The foregoing lemma is most important, because it allows us to extend the Fourier transform $f \mapsto \hat{f}$ to the whole of $L_2$. Let $f \in L_2$ be given, because $L_1 \cap L_2$ is dense in $L_2$, there exists a sequence $\{f_k\}$ in $L_1 \cap L_2$ such as, that $||f_k - f||_2 \to 0$ as $k \to \infty$. Consequently $||f_k - f_l||_2 \to 0$ as $k, l \to \infty$, and by the previous lemma it follows that

$$||\hat{f}_k - \hat{f}_l||_2 = ||f_k - f_l||_2 \to 0 \text{ as } k, l \to \infty.$$

By the completeness of $L_2$, there is a unique element $g \in L_2$ such that $\hat{f}_k \to g$ as $k \to \infty$. Define $\hat{f} := g$. Clearly

$$\begin{aligned} ||\hat{f}||_2 &= ||\lim_{k\to\infty} \hat{f}_k||_2 = \lim_{k\to\infty} ||\hat{f}_k||_2 \\ &= \lim_{k\to\infty} ||f_k||_2 = ||f||_2 \end{aligned}$$

for $f \in L_2$.
Thus we have established the important

**Theorem 1.20 (Plancherel's Theorem).** *The Fourier transform $\Phi : L_2 \to L_2$ is a linear isometry, i.e.*

$$||\Phi f||_2 = ||f||_2 \textit{ for every } f \in L_2 \tag{1.9}$$

The complete theorem and the proof can be found in Rudin(1974). Important consequences of Plancherel's Theorem are the following:

- $\Phi$ preserves the inner product in $L_2$.
- The Fourier transform is onto $L_2$.
- The image of the Fourier transform is a closed subspace of $L_2$.
- The Fourier transform is a unitary operator of $L_2$

Important consequences of Plancherel's theorem are:

- *$\Phi$ preserves the inner product in $L_2$*
  For, observe that

$$\langle u, v\rangle = \frac{1}{4}(||u+v||_2^2) - ||u-v||_2^2 + i||u+iv||_2^2 - i||u-iv||_2^2)$$

  Apply it to $u = \Phi f$, $v = \Phi y$, with $f, g \in L_2$, and recall that $||\Phi f||_2 = ||f||_2$ for every $f \in L_2$ and get:

$$\begin{aligned} 4\langle \Phi f, \Phi g\rangle &= ||\Phi(f+g)||^2 - ||\Phi(f-g)||^2 + i||\Phi(f+ig)||^2 - i||\Phi(f-ig)||^2 \\ &= ||f+g||^2 - ||f-g||^2 + i||f+ig||^2 - i||f-ig||^2 \\ &= 4\langle f, g\rangle \end{aligned}$$

  thus proving that $\Phi$ preserves the inner product in $L_2$
  Important properties of the Fourier transformation are given now:

- *The Fourier transform is onto $L_2$*
  Suppose the closed subspace $\mathrm{Im}\Phi$ is not all of $L_2$. Then, its orthogonal complement $\mathrm{Im}\Phi^{\perp}$ contains a nonzero element $g \in L_2$. By Fubini's theorem we can conclude that for $f \in L_2$

$$
\begin{aligned}
\int \hat{g}(x)\bar{f}(x)dx &= \int (\int e^{-2\pi ixy} g(y)dy)\bar{f}(x)dx \\
&= \int g(y)(\int e^{-2\pi ixy}\bar{f}(x)dx)dy \\
&= \int g(y)\hat{\bar{f}}(x)dy \\
&= \langle \hat{\bar{f}}, \bar{g} \rangle \\
&= 0
\end{aligned}
$$

  since $\hat{\bar{f}}(\xi) = \bar{\hat{f}}(-\xi)$ and therefore $\langle \hat{\bar{f}}, \bar{g}\rangle = \int \bar{\hat{f}}(-\xi)g(\xi)d\xi = \overline{\langle g, \hat{f}(-.)}$. Therefore

$$\langle \hat{g}, f \rangle = 0 \text{ for all } f \in L_2$$

  and $\hat{g} = 0$. But by Plancherel's Theorem we deduce that

$$0 \neq ||g||_2 = ||\hat{g}||_2 = 0$$

  which is a contradiction. This shows that $\mathrm{Im}\Phi = L_2$, i.e. $\Phi$ is onto.

- *The image of the Fourier transform is a closed subspace of $L_2$*
  $\mathrm{Im}\Phi$ is clearly a subspace of $L_2$. Let $\{\hat{f}_n\}$ be a sequence in $\mathrm{Im}\Phi$ , with $||\hat{f}_n - g|| \to 0$ for some $g \in L_2$. By Plancherel's Theorem it follows that

$$||f_n - f_m|| = ||\hat{f}_n - \hat{f}_m|| \to 0 \text{ as } m, n \to \infty$$

  hence $\{f_n\}$ converges to some $f \in L_2$. Again by Plancherel's Theorem

$$||\Phi f_n - \Phi f|| = ||f_n - f|| \to 0 \text{ as } n \to \infty$$

  therefore $\Phi f = g$, i.e. $g \in \mathrm{Im}\Phi$, so we see that $\mathrm{Im}\Phi$ is closed.

**Definition 1.21 (Adjoint of a linear Operator).** *The adjoint of a linear operator $T : L_2 \to L_2$ is the unique linear operator $T^* : L_2 \to L_2$ satisfying*

$$\langle Tf, g \rangle = \langle f, T^* g \rangle$$

*for each choice of $f, g \in L_2$.*

Linear isometries of $L_2$ satisfy the condition

$$T^* T = I.$$

Moreover,if they are onto, they also satisfy

$$TT^* = I$$

**Definition 1.22 (Unitary Operator).** *A linear operator $T : L_2 \to L_2$ is said to be unitary if*

$$T^*T = TT^* = I$$

Thus, a unitary operator is necessarily invertible and its invers coincides with its adjoint, i.e. $T^{-1} = T^*$.

- *The Fourier transform is a unitary operator of $L_2$*
  In fact, $\Phi$ preserve the inner product, and therefore

$$\langle \Phi f, \Phi g\rangle = \langle f, \Phi^*\Phi g\rangle = \langle f, g\rangle$$

  for every $f, g \in L_2$. Thus $\Phi^*\Phi = I$. Let $f = \Phi u$, which is possible because $\Phi$ is onto. Therefore

$$||\Phi^* f||_2^2 = ||\Phi^*\Phi u||_2^2 = ||u||_2^2 = ||f||_2^2$$

  hence $\Phi^*$ is also an isometry. Since $\Phi^{**} = \Phi$, it follows that $\Phi^{**}\Phi^* = \Phi\Phi^* = I$. This proves that $\Phi$ is unitary.

With that work we can easyly compute the inverse Fourier Tranform, it suffices to compute the adjoint $\Phi^*$ of $\Phi$. For, by Fubini's theorem

$$\begin{aligned}
\langle \Phi f, g\rangle &= \int (\int e^{-2\pi i\xi x} f(x)dx)\overline{g(\xi)}d\xi \\
&= \int f(x)\overline{(\int e^{2\pi i\xi x} g(\xi)d\xi)}dx \\
&= \langle f, \Phi^* g\rangle
\end{aligned}$$

with

$$\Phi^* g(x) = \int e^{2\pi i\xi x} g(\xi)d\xi.$$

The Fourier inversion formula is then given by

$$f(x) = \int e^{-2\pi i\xi x}\hat{f}(\xi)d\xi$$

This result follows from the fact that there is a symmetric relation between $f$ and $\hat{f}$: If

$$\begin{aligned}
\phi_A(t) &= \int_{-A}^{A} f(x)e^{-2\pi itx}dx \\
\psi_A(x) &= \int_{-A}^{A} \hat{f}(t)e^{-2\pi itx}dt
\end{aligned}$$

then

$$||\phi_A - \hat{f}||_2 \to 0 \text{ and } ||\psi_A - f||_2 \to 0 \text{ as } A \to \infty.$$

A proof for the inversion formula, which is a corollary from the theorem of Plancherel can be found in Rudin(1974). Nonetheless in probability theory the application of characteristic functions are essentially Fourier transformations of the densities, which are $L_1$ functions. We can see, that the $L_2$ theory of the Fourier tranform yields a very elegant solution to the Fourier inversion problem.

# 2

# Random Fields

The 'random fields' that are interesting for our aspects will be random mappings from subsets of Euclidean vector spaces or generally from Riemann manifolds to the real line. Since it is often no more difficult to treat far more general scenarios, they may also be real valued random mappings on any measurable space. We will give a pure mathematical definition of random fields and their realisation and even introduce the concept of stochastically equivalent, which is not often found in the literature. We are going to define the mathematical concepts of expectation and variance for random vectors. By Kolmogorov's extension theorem we can see that a random field is fully determined by the mathematical expectation and mathematical variance function. We will also have a closer look at the properties of a random field, namely the stationarity, we will define strict stationarity and stationarity in the wider sense, which is known in literature as second order stationarity. The properties of a covariance function $C$ are discussed and how they are linking with the variogram functions $\gamma$. A main result will be the construction of random fields, this result may be used for the simulation of a random field by the fast Fourier decomposition. The analytic properties of random fields are shown, for example we are going to introduce what mean square continuous/differentiable means. One main result will be the Karhunen-Loeve representation of random fields, and the spectral representation of covariance functions. The other main result will be the famous theorem of Bochner. At the end we will have a detailed look at the whole distribution of a random field. By introducing Copulas we are able to give the distribution of the random field by some basic assumptions to the underlying process.

**Definition 2.1 (Random Field).** *Let $T \subseteq \mathbb{R}^n$ be a set, a random field is defined as a function*

$$\xi(\omega, x) : \omega \times T \to \mathbb{R}^m$$

*in such way, that $\xi(\omega, x)$ is a random vector for each $x \in T$.*

In the case, when $n = 1$, $\xi(\omega, x)$ is a random process. For $n > 1$ and $m = 1$ the random field $\xi(\omega, x)$ is a scalar random field and for $m > 1$ we wil call it a vector random field. For shorter notation we will write $\xi(x), x \in T$ for $\xi(\omega, x)$. If $\xi(x), x \in T$ is a Gaussian system of random vectors then the field $\xi(x)$ will be called Gaussian field. The finite dimensional distributions for a random field $\xi(x), x \in T$ are defined as the set of distributions

$$P(\xi(x^{(k)}) \in B_k, k = 1, ..., r)$$

where the Borel Algebras $B_k \in \mathcal{B}^m$.

**Definition 2.2 (Stochastically equivalent).** *We call two random fields $\xi(x), x \in T$ and $\eta(x), x \in T$ stochastically equivalent if the finite dimensional distributions of $\xi(x), x \in T$ and $\eta(x), x \in T$ coincide.*

For stochastically equivalent the notation $\xi(x) =^d \eta(x)$, $x \in T$ will be used.

**Definition 2.3 (Realisation).** *For a fixed $\omega$, the function $\xi(x), x \in T$ is called a realization of the random field.*

**Definition 2.4.** *Let $T \in \mathcal{B}$, a random field $\xi(\omega, x) : \Omega \times T \to \mathbb{R}^m$ is measurable if for any $A \in \mathcal{B}^m$*

$$\{(\omega, x) : \xi(\omega, x) \in A\} \in \mathcal{A} \times \mathcal{B}(T)$$

*holds.*

If $\xi(x), x \in T$ is a stochastically continuous on $T$, then there exists a random field $\tilde{\xi}(x), x \in T$ with the property $\tilde{\xi}(x) =^d \xi(x)$ for $x \in T$.

**Definition 2.5 (Separable).** *We will call a random field separable with respect to the set $I \subset T$ if $I$ is countable and dense in $T$ and there exists a set $N \in \mathcal{A}$, $P(N) = 1$ such that for any ball $v(r) \subset \mathbb{R}^n$*

$$\{\omega : \subset_{x \in I \cap v(r)} \xi(x) = \subset_{x \in T \cap v(r)} \xi(x)\} \supset N$$
$$\{\omega : \inf_{x \in I \cap v(r)} \xi(x) = int_{x \in T \cap v(r)} \xi(x)\} \supset N$$

*holds.*

For any random field we can find a separable field so that they are stochastically equivalent.

**Definition 2.6.** *A second order random field over $S \subset \mathbb{R}^d$ is a function $Z : S \to L_2(\Omega, \mathcal{A}, P)$. If $d = 1$ we call $Z$ a second order stochastic process.*

This definition means, that a second oder (complex) random field is specified over $S$, if a random variable $Z(x)$ has been specified for each $x \in S$, with $E|Z(x)|^2 < \infty$. As alternative definition we can say that a second order random field over $S$ is a family $\{Z(x)|x \in S\}$ of square integrable random variables.

**Definition 2.7.** *A real random field $Z$ is Gaussian if for each $n \geq 1$ and each choice of points $x_1, ..., x_n \in S$, the corresponding random vector $(Z(x_1), ..., Z(x_n))^T$ is Gaussian.*

Let $Z$ be a real Gaussian random field and let $x_1, ..., x_n \in S$. Let $m \in \mathbb{C}^n, C \in \mathbb{C}^{n \times n}$ be given by

$$\begin{aligned} m_i &= EZ(x_i),\ i = 1, ..., n \\ C_{ij} &= E(Z(x_i) - m_i)(Z(x_j) - m_j),\ i, j = 1, ..., n \end{aligned}$$

Clearly $C$ is symmetric. Given $w \in \mathbb{R}^n$

$$\begin{aligned} \sum_{i,j=1}^{n} C_{ij} w_i w_j &= E \sum_{i,j=1}^{n} (Z(x_i) - m_i) w_i (Z(x_j) - m_j) w_j \\ &= E|\sum_{i=1}^{n} (Z(x_i) - m_i) w_i|^2 \\ &\geq 0 \end{aligned}$$

hence $C \geq 0$.

Suppose $C > 0$ and let

$$C = LL^T$$

be its Cholesky factorisation. Let us define the random vector $W$ as:

$$W := L^{-1}(Z - m)$$

where $Z$ stands for the vector $(Z(x_1), ..., Z(x_n))^T$. $W$ is a Gaussian random vector with mean 0 and covariance matrix

$$\begin{aligned} EWW^T &= L^{-1}(E(Z - m)(Z - m)^T)L^{-T} \\ &= L^{-1}(LL^T)L^{-T} \\ &= I \end{aligned}$$

In other words, $W \sim N(0, I)$ and $Z = LW + m$.

We know that

$$\phi_W(\xi) = e^{-\frac{1}{2}|\xi|^2}$$

and therefore

$$\phi_Z(\xi) = e^{im\xi - \frac{1}{2}\xi^T C^{-1} \xi}$$

holds, because if $\phi = T \circ S$ is an affine transformation in $\mathbb{R}^n$, namely the compositions of the translation $T : Tx = x + h$ with nonsingular linear transformations $x \mapsto Sx$, with $S$ invertible. Then

$$\widehat{f \circ \phi}(\xi) = \frac{e^{2\pi i h (S^t)^{-1}\xi}}{|\det S|} \hat{f}((S^T)^{-1}\xi)$$

Thus, $m$ and $C$ determine the characteristic function of $(Z(x_1), ..., Z(x_n))^T$, thereby determining the distribution. Thus, the finite dimensional distributions of $Z$ are determined by the pair of functions

$$m(x) = EZ(x)$$
$$C(x, y) = E(Z(x) - m(x))(Z(y) - m(y))$$

called the mean and covariance function.

**Theorem 2.8 (Kolmogorov's Extension Theorem).** *For all* $t_1, ...t_k \in T, k \in \mathbf{N}$ *let* $\nu_{t_1}, ..., \nu_{t_k}$ *be probability measures on* $\mathbb{R}^{nk}$ *then*

$$\nu_{t_{\sigma(1)}}, ..., \nu_{t_{\sigma(k)}}(F_1 \times ... \times F_k) = \nu_{t_1}, ..., \nu_{t_k}(F_{\sigma^{-1}(1)} \times ... \times F_{\sigma^{-1}(k)})$$

*for all permutations* $\sigma$ *on* $\{1, 2, ..., k\}$ *and*

$$\nu_{t_1}, ..., \nu_{t_k}(F_1 \times ... \times F_k) = \nu_{t_1,...,t_k,t_{k+1},...,t_{k+m}}(F_1 \times ... \times F_k \times \mathbb{R}^n \times ... \times \mathbb{R}^n)$$

*for all* $m \in \mathbf{N}$*, where the set on the right hand side has a total of* $k + m$ *factors.*
*Then there exists a probability space* $(\Omega, A, P)$ *and a stochastic process* $\{X_t\}$ *on* $\Omega$ *with* $X_t : \Omega \to \mathbb{R}^n$ *so that*

$$\mu_{t_1,...,t_k}(F_1 \times ... \times F_k) = P[X_{t_1} \in F_1, ..., X_{t_k} \in F_k]$$

*for all* $t_i \in T, k \in \mathbf{N}$ *and all Borel sets* $F_i$*.*

By the foregoing Theorem the whole random field is thus determined by these two functions.

**Definition 2.9.** *Let* $Z$ *be a (complex) second order random field over* $S$*, the mean and covariance functions of* $Z$ *are* $m_Z : S \to \mathbb{C}$ *and* $C_{ZZ} : S \times S \to \mathbb{C}$*, given by*

$$m_Z(x) = EZ(x)$$
$$C_{ZZ}(x, y) = E(Z(x) - m(x))\overline{(Z(y) - m(y))}$$

From above we see, that a Gaussian random field is characterised by its mean and covariance function, indeed $m_z$ and $C_{ZZ}$ determine the finite dimensional distributions of $Z$. In general, this is not the case, nevertheless, a surprisingly vast set of properties of an arbitrary second order random field depends only on $m_z$ and $C_{ZZ}$, but this only holds in the second oder theory.

Without loss of generality, let us now assume that $m_Z = 0$. Rotations of the Euclidean space are defined as linear transformations $g$ on the Euclidean space that do not change their orientations and preserve the distance of the points from the origin.

$$|gx| = |x|$$

As known from algebra, it is clear, that the rotations of $\mathbb{R}^n$ generate a group, commonly denoted by $\mathrm{SO}(n)$. The motions of $\mathbb{R}^n$ are defined as non-homogeneous linear transformations which preserve the distance between points and their orientation. So we can write any motion in $\mathbb{R}^n$ in the form

$$x \to gx + \tau$$

where $g \in \mathrm{SO}(n)$ and $\tau \in \mathbb{R}^n$ may be treated as an element of the group $T = \{\tau\}$ of the shifts $x \to x + \tau$. All motions in the Euclidean space generate a group $M(n)$.

**Definition 2.10 (Strict homogeneous).** *A random field $Z : \Omega \times \mathbb{R}^n \to \mathbb{R}^m$ is called homogeneous in the strict sense, if all of its finite dimensional distributions are invariant with respect to the group $T$ of shifts.*

$$Z(x) =^d Z(x + \tau)$$

**Definition 2.11 (Homogeneous in the wide Sense).** *A random field satisfying $EZ^2(x) < \infty$ is called homogeneous in the wide sense if its mathematical expectation and correlation function are invariant with respect to the group of shifts.*

Equivalent to that definition is the following, where it is assumed that the mathematical expectation is constant over the field:

**Definition 2.12.** *A second order random field is said to be homogeneous, if there is a function $R_Z : \mathbb{R}^d \to \mathbb{C}$, the correlation function, such that*

$$C_{ZZ}(x, y) = R_Z(y - x).$$

*If $d = 1$ we will call the second order stochastic process $Z$ weakly stationary with correlation function $R_Z$.*

If a field $Z$ that is homogeneous in the strict sense possesses a second order moment, it is also homogeneous in the wide sense. In the case of Gaussian one dimensional random fields both concepts of homogeneity coincide. For a homogeneous random field $Z$

$$R_Z(v) = C_Z(v + x, x) = EZ(v + x)\overline{Z(x)}$$

for any $x \in \mathbb{R}^d$. In particular

$$R_Z(0) = E|Z(x)|^2 \text{ for any } x \in \mathbb{R}^d$$

In general, let $C : S \times S \to \mathbb{C}$ be the covariance function of a certain random field $Z$ over $S$, which is a non negative definite kernel on $\mathbb{R}^n \times \mathbb{R}^n$.

**Proposition 2.13.** *A covariance function satisfies the following conditions*

- $C(x,y) = \overline{C(y,x)}$
  *follows by the definition*
- $|C(x,y)| \leq \sqrt{E|Z(x)|^3}\sqrt{E|Z(y)|^2}$
  *as a consequence of the Schwartz inequality*
- *For each choice of $n \geq 1$, of $x_1, ..., x_n \in S$ and of $u \in \mathbb{C}^n$*

$$\sum_{i,j=1}^{n} C(x_i, x_j) u_i \bar{u}_j \geq 0$$

- $C(x,x) \geq 0$

Let $C$ be the covariance function of a homogeneous second order random field over $\mathbb{R}^d$, and let $R$ be the corresponding correlation function. Then

- $R(x) = R(-x)$
- $|R(x)| \leq \sqrt{E|Z(x)|^2}\sqrt{E|Z(0)|^2}$
- $\sum_{i,j=1}^{n} R(x_i - x_j) u_i \bar{u}_j \geq 0$ for any choice of $n \geq 1$, $x_1, ..., x_n \in \mathbb{R}^d$ and $u \in \mathbb{C}^n$.

If the correlation function is continuous at the point $x = 0$, then the field is mean square continuous at each point $x \in \mathbb{R}^n$ and vice versa.

**Definition 2.14.** *Let $Z, W$ be two second order random fields with $m_Z = 0$ and $m_W = 0$, then their cross-covariance function is given by*

$$C_{Z,W} = EZ(x)\overline{W(y)}$$

Let $K$ be the class of all functions which can serve as correlation functions of a homogeneous random field and consider that $B_1(x)$ and $B_2(x) \in K$ and $k_1, k_2$ are constants, then we can infer that

$$k_1 B_1(x) + k_2 B_2(x) \in K$$
$$B_1(x) B_2(x) \in K.$$

**Definition 2.15 (Isotropic in the wide Sense).** *A random field is called isotropic in the wide sense, if $EZ^2(x) < \infty$ and the mathematical expectation $m(x)$ along with the correlation function $R(x,y)$ are invariant with respect to the group $SO(n)$:*

$$m(x) = m(gx)$$
$$R(x,y) = R(gx, gy)$$

*for any $x, y \in \mathbb{R}^n$ and any $g \in SO(n)$.*

There exist isotropic fields which are not homogeneous, in literature we can find an example, namely the Levy multi parametric Brownian motion, that is a Gaussian random field $V(x)$, $x \in \mathbb{R}^n$, with $EV(x) = 0$ and $EV(x)V(y) = \frac{|x|+|y|-|x-y|}{2}$ with $x, y \in \mathbb{R}^n$.
If we have a homogeneous isotropic random field, then the mathematical expectation is constant and the correlation function $R(x, y)$ only depends on the Euclidean distance between the points $x, y \in \mathbb{R}^n$.

## 2.1 Construction of Random Fields

Let $K$ be a compact set in $\mathbb{R}^d$ and let $Z$ be a second order random field over $K$. Let $C$ be the covariance function of $Z$, and it is assumed that the covariance function $C$ is square integrable over $K \times K$ in what follows, i.e.

$$\int_K \int_K |C(x, y)|^2 dxdy < \infty$$

Let $A : L_2(K) \to L_2(K)$ be the integral operator

$$A\phi(x) = \int_K C(x, y)\phi(y)dy.$$

We can infer that $A$ is a self-adjoint, linear compact operator on $L_2(K)$. Let $\lambda_1, \lambda_2, ...$ be the eigenvalues of $A$, with orthonormal eigenfunctions $\psi_1, \psi_2, ...$, i.e.

$$A\psi_n = \lambda_n \psi_n,\ n = 1, 2, ...$$
$$\langle \psi_n, \psi_m \rangle = \delta_{n,m},\ n, m = 1, 2, ...$$

where the inner product is defined as

$$\langle f, g \rangle = \int_K f(x)\overline{g(x)}dx$$

for $f, g \in L_2(K)$.

Suppose moreover that $C$ is continuous. Then

$$\langle Af, f \rangle \geq 0,\ f \in L_2(K)$$

which follows directly from the covariance condition: $\sum_{i,j=1}^n C(x_i, x_j)u_i\overline{u_j} \geq 0$. Thus the eigenvalues are non negative, besides converging to zero as $n \to \infty$. In addition, $Af$ is continuous for each $f \in L_2(K)$, so that $\psi_n$ is continuous if $\lambda_n > 0$. Indeed

$$\psi_n = \lambda_n^{-1} A\psi_n$$

Beside, Mercer's theorem holds, i.e.

$$C(x,y) = \sum_{k=1}^{\infty} \lambda_k \psi_k(x) \overline{\psi_k(y)} \tag{2.1}$$

convergence being absolute and uniform over $K \times K$.

Let $Z_1, Z_2, ...$ be a sequence of (complex) random variables, with

$$\begin{aligned} EZ_n &= 0, n = 1, 2, ... \\ EZ_n \overline{Z_m} &= \lambda_n \delta_{n,m} \text{ for } n, m = 1, 2... \end{aligned}$$

**Theorem 2.16.** *Let $\{\lambda_n\}, \{\psi_n\}$ and $\{Z_n\}$ be a sequence of eigenvalues, orthogonal eigenfunctions and (complex) random variables. In particular, assume $C$ is continuous. Then the series*

$$Z(x) := \sum_{n=1}^{\infty} \psi_n(x) Z_n \tag{2.2}$$

*converges in quadratic mean for any $x \in K$. Thus $Z$ is a second order random field over $K$. Moreover*

$$\begin{aligned} EZ(x)\overline{Z_n} &= \lambda_n \psi_n(x), \text{ for } n = 1, 2, ... \\ \lim_{x \to x_0} E|Z(x) - Z(x_0)|^2 &= 0 \text{ for any } x_0 \in K \end{aligned}$$

PROOF: For each $n \geq 1$, let

$$S_n(x) = \sum_{k=1}^{n} \psi_k(x) Z_k$$

For $n, m \geq 1, n > m$,

$$\begin{aligned} E|S_n(x) - S_m(x)|^2 &= E|\sum_{k=m+1}^{n} \psi_k(x) Z_k|^2 \\ &= \sum_{k=m+1}^{n} \lambda_k |\psi_k(x)|^2 \to 0 \end{aligned}$$

in virtue of (2.1). This proves the convergence in quadratic mean of the series in (2.2), thereby proving that $Z$ is a second order random field. Moreover,

$$ES_n(x)\overline{Z_k} = \sum_{k=1}^{n} \psi_k(x) EZ_k \overline{Z_l} = \lambda_l \psi_l(z).$$

The limit follows upon letting $n \to \infty$. Finally we observe that,

$$|Z(x) - Z(x_0)|^2 = |Z(x)|^2 - Z(x)\overline{Z(x_0)} - Z(x_0)\overline{Z(x)} + |Z(x_0)|^2 \quad (2.3)$$

and therefore

$$\begin{aligned} E|Z(x) - Z(x_0)|^2 &= C(x,x) - C(x,x_0) - C(x_0,x) + C(x_0,x_0) \\ &\to C(x_0,x_0) - 2C(x_0,x_0) + C(x_0,x_0) = 0 \end{aligned}$$

as $x \to x_0$.

□

## 2.2 Analytic Properties

**Definition 2.17.** *A second order random field $Z$ over $S \subset \mathbb{R}^d$ is mean square (m.s.) continuous at $x_0 \in S$ if $|E(Z(x) - Z(x_0))|^2 \to 0$ as $x \to x_0$.*

In order for a field with $EZ^2 < \infty$ to be mean squared continuous, it is necessary and sufficient that the covariance function is continuous along the diagonal $\{(x, y) : x = y\}$, a simplification is given in the following lemma:

**Lemma 2.18.** *A random field $Z$ is m.s. continuous at $x_0$ if and only if its covariance function is continuous at $(x_0, x_0)$*

PROOF: Observe that

$$|Z(x) - Z(x_0)|^2 = |Z(x)|^2 - 2\text{Re}Z(x)\overline{Z(x_0)} + |Z(x_0)|^2,$$

hence

$$E|Z(x) - Z(x_0)|^2 = C(x,x) - 2\text{Re}C(x,x_0) + C(x_0,x_0).$$

Clearly the continuity of $C$ at $(x_0, x_0)$ implies that $E|Z(x) - Z(x_0)| \to 0$ as $x \to x_0$.
Conversely, suppose $Z$ is m.s. continuous at $x_0$, and let $x \to x_0$ and $y \to y_0$. Then

$$C(x,y) = EZ(x)\overline{Z(y)} \to EZ(x_0)\overline{Z(x_0)} = C(x_0,x_0).$$

□

Let $Z$ be the m.s. continuous random field

$$Z(x) = \sum_{n=1}^{\infty} Z_n \psi_n(x)$$

$Z_n$ are the Fourier coefficients of $Z$ with respect to the orthonormal system $\{\psi_n\}$, i.e.

$$Z_n = \int_K Z(x)\overline{\psi_n(x)}dx \tag{2.4}$$

For the evaluation of (2.4) we consider first following integral of the form:

$$\int_G f(x)Z(x)dx \tag{2.5}$$

where $G$ is a bounded region in $\mathbb{R}^d$ and $f : G \to \mathbb{C}$ is continuous. (2.5) will be defined in the sense of Riemann, with convergence in mean square.

Let $\{Q_n\}_{n=1}^{\infty}$ be a partition of $\mathbb{R}^d$ into cubes with edges parallel to the coordinate axes, and let

$$\delta := \sup_{n \geq 1} \operatorname{diam}(Q_n)$$

where

$$\operatorname{diam}(S) = \sup\{|x - y| : x, y \in S\},$$

for any set $S \subset \mathbb{R}$. Number the cubes in such a way that

$$G \subset \sum_{k=1}^{m} Q_k$$

for some $m \geq 1$, this is possible because $G$ is bounded. That is, $G \cap Q_k \neq 0$ for $k = 1, ..., m$ and $G \cap Q_l = 0$ for $l > m$. Let $|Q_k|$ denote the $d$-dimensional volume of the $k$-th cube, and pick $x_k \in G \cup Q_k$, for $k = 1, ..., m$. From the Riemann sum

$$S_m := \sum_{k=1}^{m} f(x_k)Z(x_k)|Q_k| \tag{2.6}$$

whose value is a random variable in $L_2(\Omega\mathcal{A}, P)$, and depends on both the partition of $\{Q_n\}$ and on the choice of $x_k \in G \cap Q_k$, for $k = 1, ..., m$. If (2.6) converges in mean square when $\delta \to 0$.

$$\int_G f(x)Z(x)dx := \text{l.i.m}_{\delta \to 0} \sum_{k=1}^{m} f(x_k)Z(x_k)|Q_k|$$

Observe that

$$\begin{aligned} E|S_m|^2 &= \sum_{k,l=1}^{m} f(x_k)EZ(x_k)\overline{Z(x_l)f(x_l)}|Q_k||Q_l| \\ &= \sum_{k,l=1}^{m} f(x_k)C(x,y)\overline{f(y)}dxdy \end{aligned}$$

If $f(x)C(x,y)f(y)$ is Riemann integrable over $G \times G$, then the last sum converges to

$$\int_G \int_G f(x)C(x,y)\overline{f(y)}dxdy \tag{2.7}$$

The converse is also true:

**Theorem 2.19.** *The integral (2.5) exists in the mean square sense if and only if the integral (2.7) exists in the Riemann sense, and*

$$E\int_G f(x)Z(x)dx = 0$$
$$E|\int_G f(x)Z(x)dx|^2 = \int_G \int_G f(x)C(x,y)\overline{f(y)}dxdy$$

PROOF: Indeed, for $m > n$,

$$S_m - S_n = \sum_{k=n+1}^{m} f(x_k)Z(x_k)|Q_k|$$

and

$$E|S_m - S_n|^2 = \sum_{k=n+1}^{m} \sum_{k=n+1}^{m} f(x_k)C(x_k, x_l)\overline{f(x_l)}|Q_k||Q_l|$$

which spells out the assertion of the theorem, in view of the completeness of both $\mathbb{C}$ and $L_2(\Omega, \mathcal{A}, P)$.

□

Given a compact set $K \subset \mathbb{R}^n$ and a continuous function $C : K \times K \to \mathbb{C}$, the corresponding eigenfunctions $\{\psi_n\}$ are continuous. Therefore the double integral

$$\int_K \int_K \psi_m(x)C(x,y)\overline{\psi_n(y)}dxdy$$

always exists in the sense of Riemann. If $Z$ is a m.s. random field constructed by:

$$Z(x) = \sum_{n=1}^{\infty} Z_n\psi_n(x)$$

then the integrals

$$Z_n = \int_K Z(x)\overline{\psi_n(x)}dx$$

exist for each $n$.

**Theorem 2.20 (Karhunen-Loève Representation).** *A random field $Z$ with continuous covariance function $C$ on a compact set $K \subset \mathbb{R}^n$ is m.s. continuous if and only if*

$$Z(x) = \sum_{n=1}^{\infty} Z_n \psi_n(x)$$

*converges in the mean square sense.*
*Here $\{\psi_n\}$ is the sequence of eigenfunctions associated with the kernel $C(x,y)$, with corresponding eigenvalues $\{\lambda_n\}$. Moreover, $\{Z_n\}$ is a sequence of orthogonal random variables, with zero mean and*

$$EZ_k\overline{Z_l} = \delta_{kl}\lambda_l,\ k,l = 1,2,... \tag{2.8}$$

*Moreover*

$$EZ(x)\overline{Z_k} = \lambda_k \psi_k(x),\ k = 1,2,... \tag{2.9}$$

*In fact the inner product is given by*

$$Z_k(\omega) = \langle Z(.,\omega), \psi_k \rangle,\ k = 1,2,...$$

PROOF: Given a m.s. continuous random field $Z$ over $K$, the sequence of its Fourier coefficients $\{Z_n\}$ is well defined if $C$ is continuous. Then

$$\begin{aligned} EZ_k\overline{Z_l} &= E\int_K \overline{\psi_k(x)}Z(x)dx \overline{\int_K \overline{\psi_l(y)}Z(y)dy} \\ &= \int_K\int_K \overline{\psi_k(x)}C(x,y)\psi_l(y)dxdy \\ &= \int_K (\int_K C(x,y)\psi_l(y)dy)\overline{\psi_k(x)}dx \\ &= \lambda_l \langle \psi_l, \psi_k \rangle \end{aligned}$$

thus proving (2.8). Condition (2.9) is proved in an analogous way. The converse is given in theorem (2.16).

□

**Definition 2.21.** *A (complex) random field $Z$ on a domain $D$, an open connected subset of $\mathbb{R}^d$ is m.s. differentiable at $x_0 \in D$ and its derivative there is the random vector $U \in L_2(\Omega, \mathcal{A}, P)$ if*

$$Z(x_0 + h) = Z(x_0) + Uh + \varepsilon(x_0, h) \tag{2.10}$$

*where*

$$\frac{E|\varepsilon(x_0,h)|^2}{|h|^2} \to 0 \text{ as } h \to 0$$

The derivative of $Z$ at the point $x_0$ will be denoted by $Z'(x_0)$.
By Jenssen's inequality

$$|E\varepsilon(x_0, h)|^2 \leq E|\varepsilon(x_0, h)|^2$$

hence

$$0 \leq \frac{E|\varepsilon(x_0, h)|}{h} \leq \sqrt{\frac{E|\varepsilon(x_0, h)|^2}{|h|^2}} \to 0 \tag{2.11}$$

If $Z$ is differentiable at $x_0$ and by evaluating the expected value of (2.10) we get:

$$EZ(x_0 + h) = EZ(x_0) + EUh + E\varepsilon(x_0, h)$$

which together with (2.11) implies the differentiability of $EZ$ at $x_0$. It holds that:

$$(EZ)'(x_0) = EZ'(x_0)$$

Observe that

$$\begin{aligned}
E| \ Z \ &(x_0 + h) - Z(x_0) - Uh|^2 = \\
&= C(x_0 + h, x_0 + h) - C(x_0 + h, x_0) - C(x_0, x_0 + h) \\
&\quad + C(x_0, x_0)(E(\overline{Z(x_0 + h) - Z(x_0)})U + E(Z(x_0 + h) - Z(x_0))\overline{U})h \\
&\quad + E|Uh|^2 \\
&= C(x_0 + h, x_0 + h) - C(x_0 + h, x_0) - C(x_0, x_0 + h) + C(x_0, x_0) \\
&\quad - E|Uh|^2 + \mathrm{Re}\overline{\varepsilon(x_0, h)}Uh
\end{aligned}$$

upon recalling the differentiability condition (2.10). Then

$$\begin{aligned}
&|\frac{C(x_0 + h, x_0 + h) - C(x_0 + h, x_0) - C(x_0, x_0 + h) + C(x_0, x_0)}{|h|^2} \\
&\qquad - \frac{E|Uh|^2}{|h|^2}| \leq \frac{E|Z(x_0 + h) - Z(x_0) - Uh|^2}{|h|^2} \\
&\qquad\qquad + \frac{E|\varepsilon(x_0, h)|}{|h|}|U|
\end{aligned}$$

Observe that $|Uh|^2$ is a quadratic form of the coordinates of $h$. Let $Q$ be its matrix, necessarily Hermitian.

**Definition 2.22.** *The correlation function $C$ admits a generalised derivative at $x_0 \in D$ if there is a Hermitian matrix $Q \in \mathbb{C}^{d \times d}$ such that:*

$$C(x_0 + h, x_0 + h) - C(x_0 + h, x_0) - C(x_0, x_0 + h) + C(x_0, x_0) = h^T Q h + \eta(x_0, h)$$

*where*

$$\frac{|\eta(x_0,h)|}{|h|^2} \to 0 \ \textit{as} \ h \to 0$$

*Q is denoted by the symbol*

$$\frac{\partial^2 C(x_0,x_0)}{\partial x \partial y}.$$

Thus we have the following characterisation of m.s. differentiability for a second order random field.

**Theorem 2.23.** *Let D be a domain in* $\mathbb{R}^d$*. The random field Z is m.s. differentiable at* $x \in D$ *if and only if the generalised derivative*

$$\frac{\partial^2 C_{ZZ}(x,x)}{\partial x \partial y}$$

*exists. If such is the case, then the following properties hold*

- $EZ'(x) = (EZ)'(x_0)$
- $C_{ZZ'}(x,y) = \frac{\partial C_{ZZ}(x,y)}{\partial x}$
- *The generalised derivative* $\frac{\partial^2 C_{ZZ}(x,y)}{\partial x \partial y}$ *exists and*

$$C_{Z'Z'}(x,y) = \frac{\partial^2 C_{ZZ}(x,y)}{\partial x \partial y}$$

PROOF:

$$C_{ZZ'}(x,y) = EZ(x)Z'(y)$$

therefore we get:

$$\begin{aligned} C_{ZZ}(x+h,y) - C_{ZZ}(x,y) &= E[Z(x+h) - Z(x)]\overline{Z(y)} \\ &= EZ'(x) \cdot h\overline{Z(y)} + E\epsilon(x,h)\overline{Z(y)} \\ &= EZ'(x)\overline{Z(y)} \cdot h + E\epsilon(x,h)\overline{Z(y)} \end{aligned}$$

i.e.

$$C_{ZZ}(x+h,y) - C_{ZZ}(x,y) = C_{Z'Z}(x,y) \cdot h + \delta(x,y;h)$$

where

$$\frac{\delta(x,y;h)}{|h|} \to 0 \text{ as } h \to 0$$

This proves $C_{ZZ'}(x,y) = \frac{\partial C_{ZZ}(x,y)}{\partial x}$.
Observe that:

$$\frac{|C(x+h,y+h) - C(x+h,y) - C(x,y+k) + C(x,y) - (U \cdot h)\overline{U \cdot h}|}{|k||h|}$$

tends to 0 as $h, k \to 0$. Recall that $(U \cdot h)\overline{U \cdot h}$ is a bilinear form in $h$ and $k$, hence there is a matrix $B \in \mathbb{C}^{n \times n}$ such that

$$(U \cdot h)\overline{U \cdot h} = Bh \cdot k$$

thus proving that $\frac{\partial^2 C(x,y)}{\partial x \partial y}$ exists. Indeed, it equals $B$. Moreover

$$\begin{aligned} E[Z(x+h) \cdot Z(x)]\overline{[Z(y+k) \cdot Z(y)]} &= (Z'(x) \cdot h)(Z'(y) \cdot k) + o(|k||h|) \\ &= \frac{\partial^2 C(x,y)}{\partial x \partial y} h \cdot k + o(|h||k|) \end{aligned}$$

which proves $C_{Z'Z'}(x,y) = \frac{\partial^2 C_{ZZ}(x,y)}{\partial x \partial y}$.

□

**Lemma 2.24.** *Let $Z$ be a homogeneous random field over the domain $D \subset \mathbb{R}^d$. Then, $Z$ is m.s. differentiable if and only if the generalised second derivative*

$$\frac{d^2 C}{dx^2}(0)$$

*exists. If this condition holds, then $\frac{d^2 C}{dx^2}$ exists for all $x \in D$, and*

$$\begin{aligned} C_{ZZ'}(x, x+y) &= \frac{dC}{dx}(y) \\ C_{Z'Z'}(x, x+y) &= -\frac{d^2 C}{dx^2}(y) \end{aligned}$$

## 2.3 Spectral Representation of Covariances

The characteristic function of a $d$-dimensional random vector $X$ has been defined as

$$\phi_X(\xi) = Ee^{i\xi X}$$

and therefore

$$\phi_X(\xi) = \int e^{i\xi X} dF_X(x).$$

Also

$$\phi_X(\xi) = \int e^{i\xi X} \mu_X(dx),$$

where $\mu_X$ is a measure defined on the Borel Subsets of $\mathbb{R}$ by

$$\mu_X(B) = P(X \in B)$$

The following theorem says, that the distribution $\mu_X$ of $X$ is determined by the characteristic function of $X$.

**Theorem 2.25.** *Let $\mu, \nu$ be two finite Borel measures such that*

$$\int e^{i\xi X}\mu(dx) = \int e^{i\xi X}\nu(dx) \tag{2.12}$$

*Then $\mu = \nu$*

PROOF: Let $\mathcal{M}$ be the class of all complex valued bounded measurable functions $f$ for which

$$\int f(x)\mu(dx) = \int f(x)\nu(dx). \tag{2.13}$$

Observe that

- $\mathcal{M}$ is not empty and is closed under linear combinations
- $\mathcal{M}$ is closed with respect to the operation of taking limits of uniformly bounded convergent sequences of measurable functions.

Hence it contains all trigonometric polynomials such as

$$\sum_{k=1}^{n} c_k e^{i\xi_k x} \tag{2.14}$$

By the Weierstrass approximation theorem, any bounded continuous function is the uniform limit of a uniformly bounded sequence of trigonometric polynomials like (2.14). Hence $\mathcal{M}$ contains all continuous bounded functions.
If $f$ is any bounded measurable function, there is a uniformly bounded sequence $\{f_n\}_{n=1}^{\infty}$ of continuous functions such that $f_n \to f$ pointwise. Therefore $f \in \mathcal{M}$. Taking $f = \aleph_B$ in (2.13) for an arbitrary $B \in \mathcal{B}$ yields $\mu(B) = \nu(B)$, hence $\mu = \nu$.

□

**Definition 2.26.** *Let $\phi : \mathbb{R}^d \to \mathbb{C}$ be a Borel measurable function. A spectral measure for $\phi$ is a finite measure $\mu : B \to [0, \infty)$ such that*

$$\phi(x) = \int e^{i\xi x}\mu(d\xi)$$

The last theorem says, that if a spectral measure exists, then it is a unique one. If a function $\phi$ has a spectral measure, then it is symmetric, e.g.

$$\phi(y - x) = \overline{\phi(x - y)}.$$

Moreover, for $k \geq 1$ and $v_1, ..., v_k \in \mathbb{C}$ then

$$\begin{aligned}
\sum_{l,m=1}^{k} \phi(x_l - x_m)v_l\overline{v}_k &= \sum_{l,m=1}^{k} \int e^{i(x_l - x_m)\xi}\mu(d\xi)v_l\overline{v}_m \\
&= \int |\sum_{l=1}^{k} v_l e^{ix_l\xi}|^2\mu(d\xi) \\
&\geq 0
\end{aligned}$$

$\phi$ is of positive type.
Thus $\phi$ could very well be the correlation function of a homogeneous random field, indeed

**Proposition 2.27.** *Let $\phi : \mathbb{R}^d \to \mathbb{C}$ have a spectral measure. Then $\phi$ is the correlation function of a homogeneous m.s. continuous random field.*

PROOF: Observe that $\phi$ is continuous, since

$$|\phi(x) - \phi(y)| \leq \int |e^{i\xi x} - e^{i\xi y}|\mu(d\xi)$$

and $|e^{i\xi x} - e^{i\xi y}| \leq 2$. Letting $y \to x$, get $\phi(y) \to \phi(x)$, hence continuity. In addition, $\phi$ is known to be symmetric and of positive type.
Let $W := \mu(\mathbb{R}^d)$. Let $D$ be a random vector in $\mathbb{R}^d$, with the distribution

$$P(D \in B) = \frac{\mu(B)}{W}.$$

Let $V$ be uniformly distributed on $[-\pi, \pi]$, ,independent of $D$. For each $x \in \mathbb{R}^d$, let

$$\begin{aligned} Z(x) &= \sqrt{W} \int_{-\pi}^{\pi} \int_{\mathbb{R}^d} e^{i(\xi x + \eta)} \frac{d\mu(\xi)}{W} \frac{d\eta}{2\pi} \\ &= \frac{\sqrt{W}}{2\pi W} (\int_{-\pi}^{\pi} e^{i\eta} d\eta)(\int_{\mathbb{R}^d} e^{i\xi x} \mu(d\xi)) \\ &= 0 \end{aligned}$$

and

$$\begin{aligned} EZ(x)\overline{Z(y)} &= MEe^{iD(x-y)} \\ &= M \int e^{i\xi(x-y)} \frac{\mu(d\xi)}{M} \\ &= \phi(x - y). \end{aligned}$$

Hence $Z$ is a m.s. continuous random field with correlation function $\phi$.

□

**Theorem 2.28 (Bochner).** *A continuous function of positive type has a unique spectral measure.*

For a proof we refer to Gikhman and Skorokhod (1996).

□

If the spectral measure $\mu$ of $\phi$ has a density $g \in L_1(\mathbb{R}^d)$, i.e.

$$\mu(B) = \int_B g(\xi) d\xi$$

then

$$\phi(x) = \int e^{ix\xi} g(\xi) d\xi.$$

The function $g$ is called the spectral density of $\phi$.
If there is a spectral density $g$, then

$$|\phi(x)| \leq \int |g(\xi) d\xi| = ||g||_1$$

i.e.

$$||\phi||_\infty \leq ||g||_1 \tag{2.15}$$

**Lemma 2.29.** *A sufficient condition for the existence of a spectral density is that* $\phi \in L_1$

## 2.4 Distribution of a Random Field

One interesting idea in the application of stochastic analysis in finance is that it is possible to give the density of a stochastic process. By the theorem of Sklar the multivariate distribution of the process can be represented in terms of its underlying margins by binding them together using a copula function. The copula approach is a modelling strategy where a joint distribution is induced by specifying marginal distributions and a function called copula that binds them together. The copula parameterises the dependence structure of the random variables by capturing all of the joint behaviour. This approach by using copulas is rather new in literature, and the usage in the context of random field is firstly introduced in this work.

### 2.4.1 Sklar's Theorem and Copulas

The copula for a $n$ dimensional multivariate distribution function $F_{12\ldots n}$ with given one dimensional marginal distribution functions $F_1$, ...,$F_n$ is a function that binds together the margins in such a manner as to form precisely the joint distribution function. (The properties of a distribution function and their corresponding density function must hold!) The fact that a copula builds the joint density function from the margins implies, that the dependency of the margins, measured by the correlations, have to be received. Sklar's theorem shows that there exists a copula function which acts to represent the joint distribution function of random variables in terns of the underlying one dimensional margins. Let for an example $F_1(x_1)$ and $F_2(x_2)$ denote the margins, the cumulative distribution function of the random variable $X_1$ and $X_2$ is given by

$$F_i(x_i) = P(X_i \leq x_i)$$

where $i = 1, 2$ and $x_i \in \bar{\mathbb{R}} = \mathbb{R} \cup \{-\infty, \infty\}$ and let

$$F_{1,2}(x_1, x_2) = P(X_1 \leq x_1, X_2 \leq x_2)$$

denote the joint cumulative distribution function. Then the result for the two dimensional case from the theorem of Sklar is, that there exists a unique two place function $C_\theta$, with a parameter vector $\theta$, such that the joint cumulative density function has a representation of the form

$$F_{1,2}(x_1, x_2) = C_\theta(F_1(x_1), F_2(x_2))$$

The function $C_\theta$ is known as the Copula. $\theta$ is not present in either margin, and the copula functions do not depend on the marginal distributions. The theorem of Sklar holds in $n$ dimensions.

**Definition 2.30 (Copula).** *The function $C_\theta$ is a copula if following conditions are fulfilled*

- $C_\theta(u, v) = 0$ *if either or both $u$ and $v$ are zero*
- $C_\theta(1, v) = v$ *and* $C_\theta(u, 1) = u$, *where* $(u, v) \in \Pi$[1]

### 2.4.2 Archimedean Copulas

Of great importance in the context of random fields is the family of Archimedean copulas, because they are captured by an additive generator function $\phi : \Pi \rightarrow [0, \infty]$. This additive generator function is a continuous, convex decreasing function

$$\phi'(t) < 0$$
$$\phi''(t) > 0$$

with the property that $\phi(1) = 0$. In the bivariate case following relation holds:

$$\phi(C_\theta(u, v))) = \phi(u) + \phi(v)$$
$$C_\theta(u, v) = \phi^{-1}(\phi(u) + \phi(v))$$

Copula functions in that class that are interesting for random fields are the Clayton and Joe copula with fat tails and Frank and Gumbel copulas.

**Definition 2.31.** *Clayton Copula The Clayton copula is given in the form*

$$(u^{-\theta} + v^{-\theta} - 1)^{\frac{1}{\theta}}$$

*with* $\theta \in (-\infty, \infty)$ *and the generator* $\phi(t)$ *given by*

$$\frac{1}{\theta}(t^{-\theta} - 1)$$

**Definition 2.32 (Frank Copula).** *The Frank copula is given in the form*

$$-\theta^{-1}\log(1+\frac{(e^{-\theta u}-1)(e^{-\theta v}-1)}{e^{-\theta}-1})$$

*with* $\theta \in (-\infty, \infty)$ *and the generator* $\phi(t)$ *given by*

$$-\log\frac{e^{-\theta t}-1}{e^{-\theta}-1}$$

**Definition 2.33 (Gumbel Copula).** *The Gumbel copula is given in the form*

$$e^{-((-\log(u))^{\theta}+(-\log(v))^{\theta})^{\frac{1}{\theta}}}$$

$\theta \in [1, \infty)$ *and the generator* $\phi(t)$ *given by*

$$(-\log(t))^{\theta}$$

And finally the Joe copula is defined by

**Definition 2.34 (Joe Copula).** *The Joe copula is given in the form*

$$1-((1-u)^{\theta}+(1-v)^{\theta}-(1-u)^{\theta}(1-v)^{\theta})^{\frac{1}{\theta}}$$

$\theta \in [1, \infty)$ *and the generator* $\phi(t)$ *given by*

$$-\log(1-(1-t)^{\theta})$$

### 2.4.3 Density of a Random Field

Under the assumption that a random field is fully defined by the trend function $m(x)$ and by the covariance function $C(x, y)$ we can now build the density of the random function by the margins of the decomposition of the random field in

$$Z(x) = m(x) + \varepsilon(x)$$

If we set the trend function equal to zero, and use a parametric covariance function with covariance parameters $\upsilon_1, ...., \upsilon_n$ we can evaluate by simulation of possible underlying covariance functions the margins of every covariance parameter $\upsilon_i$ and then use a copula to find the corresponding cumulative density function of the random field. For our purpose we have more than $(u, v)$ namely $u_1, ..., u_n$ and our aim is to estimate the dependency parameter $\theta$ for a given copula $C_\theta$. By simulation we get a multivariate observation vector $X_t = (\upsilon_{1t}, ..., \upsilon_{nt})$, $t = 1, ..., T$ with the realisations of possible covariance parameters. The corresponding marginal distribution function $F_k(., \theta_k)$ and density functions $f_k(., \theta_k)$ for $k = 1, ..., n$ can be evaluated by simple smoothing methods. The copula based parametric model of the vector $X$ so that its joint distribution function is expressed as

$$F_k(u1, ..., u_n, \theta, \alpha) = C(F_1(u_1, \theta), ..., F_n(u_n, \theta))$$

with $\alpha$ denoting the dependence parameter. We assume that the mixed partial derivatives of $C$ exists and denote it by

$$c(u_1, ..., u_n) = \frac{\partial^n C(u_1, ..., u_n)}{\partial u_1 ... \partial u_n}$$

then the joint density is given by

$$f(\upsilon_1, ..., \upsilon_n, \theta, \alpha) = c(F_1(\upsilon_1, \theta), ..., F_n(x_n, \theta_n)) \prod_{k=1}^{n} f_k(\upsilon_k, \theta_k)$$

The log likelihood function for our simulated sample of size $T$ is expressed as

$$L(\theta_1, ..., \theta_n) = \sum_{t=1}^{T} \log f(\upsilon_{1t}, ..., \upsilon_{nt}, \theta, \alpha)$$

By maximising the log likelihood function we are able to find the unknown parameters. For that problem we have implemented a R code for the Frank and the Joe copula. We now give the code for random field where the covariance function is the sum of two Matern function

```
densityRF<-function(simulated){
#####
# simulated ... output of varsim with 7 parameters
#####
  nugget <- anpassungKonvex[,1]
  sill1 <- anpassungKonvex[,2]
  range1 <- anpassungKonvex[,3]
  nue1 <- anpassungKonvex[,4]
  sill2 <- anpassungKonvex[,5]
  range2 <- anpassungKonvex[,6]
  nue2 <- anpassungKonvex[,7]

  VerteilungsFunktion<-function(data,findData){
    nug<-data
    fnug<-findData
    s<-summary(nug)
    if(s[1]==s[3]&&s[1]==s[6]&&s[3]==s[6]) {
      nugwahr<-c(rep(1,length(data)))
      return(nugwahr)
    }
    else {
```

```
    unterteilungsPunkte<-length(nug)
    Density<-density(nug,n=unterteilungsPunkte)
    cumDens<-cumsum(Density$y)/sum(Density$y)

    Spline<-spline(Density$x,cumDens)
    plot(Spline,type="l",col="red")

    f<-splinefun(Density$x,cumDens)

    wertF<-c(rep(0,length(fnug)))
    for(i in 1:length(fnug)) {
      wertF[i]<-cumDens[i]
      points(Density$x[i],cumDens[i])
    }
    return(wertF)
  }
}

 Dfrank<-function(th,nugget,sill1,range1,nue1,sill2,range2,nue2){

   hsum<-nue1+nue2+nugget+range1+range2+sill1+sill2+0.0
   hprod<-(th^nue1-1.0)*(th^nue2-1)*(th^nugget-1)*(th^range1-1)*
           (th^range2-1)*(th^sill1-1)*(th^sill2-1)
   hprod2<-1+(1/(th-1.0)^6)*hprod
   l<-log(th)
   erg<-(720*th^hsum*hprod^6*l^6)/((th-1)^42.0*hprod2^7)-(2520*
         th^hsum*hprod^5*l^6)/((th-1)^36*hprod2^6)+(3360*th^hsum
         *hprod^4*l^6)/((th-1)^30*hprod2^5)-(2100*th^hsum*hprod^
         3*l^6)/((th-1)^24*hprod2^4)+(602*th^hsum*hprod^2*l^6)/
         ((th-1)^18*hprod2^3)-(63*th^hsum*hprod*l^6)/((th-1)^12*
         hprod2^2)+(th^hsum*l^6)/((th-1)^6*hprod2)
   return(erg)
   cat(".")
 }

 LogLikli<-function(th,nugget,sill1,sill2,range1,range2,nue1,nue2)
   {
   cat("*")
   nugget<-nugget
   sill1<-sill1
   sill2<-sill2
   range1<-range1
   range2<-range2
   nue1<-nue1
   nue2<-nue2

   dnugget <- varsim[,1]
   dsill1 <- varsim[,2]
   drange1 <- varsim[,3]
```

```
  dnue1 <- varsim[,4]
  dsill2 <- varsim[,5]
  drange2 <- varsim[,6]
  dnue2 <- varsim[,7]
  nt <- length(nugget)
  if(summary(dnugget)[1]==summary(dnugget)[3]&&summary(dnugget)
    [1]==summary(dnugget)[6]&&summary(dnugget)[3]==summary
    (dnugget)[6])
  {

    ddnugget<-c(rep(1,nt))
  }
  else{
    ddnugget <- density(dnugget,n=nt)$y
  }
  ddsill1 <- density(dsill1,n=nt)$y
  ddrange1 <- density(drange1,n=nt)$y
  ddnue1 <- density(dnue1,n=nt)$y
  ddsill2 <- density(dsill2,n=nt)$y
  ddrange2 <- density(drange2,n=nt)$y
  ddnue2 <- density(varsim[,7],n=nt)$y
  cat("*")

  g<-Dfrank(th,nugget,sill1,range1,nue1,sill2,range2,nue2)

  f<-log(g)
  logliki<-f[is.finite(f)]
  logliki<-(-1)*sum(logliki)
  cat(".")
  return(logliki)
}

estCopFrank<-function(th,nugget,sill1,sill2,range1,range2,
             nue1,nue2) {
  cat(".")
  nuggetV<-VerteilungsFunktion(nugget,nugget)
  sill1V<-VerteilungsFunktion(sill1,sill1)
  range1V<-VerteilungsFunktion(range1,range1)
  nue1V<-VerteilungsFunktion(nue1,nue1)
  sill2V<-VerteilungsFunktion(sill2,sill2)
  range2V<-VerteilungsFunktion(range2,range2)
  nue2V<-VerteilungsFunktion(varsim[,7],varsim[,7])
  startwerte<-seq(10,20,0.5)
  thwerte<-c(rep(0,length(startwerte)))
  valuewerte<-c(rep(0,length(startwerte)))

  for(i in 1:length(startwerte)){
    th<-startwerte[i]
    tt1<-try(optim(th,fn=LogLikli,nugget=nuggetV,sill1=sill1V
```

```
              ,sill2=sill2V, range1=range1V,range2=range2V,
              nue1=nue1V,nue2=nue2V,method="BFGS", control=
              list(maxit=200000)))
    cat("\n\tOptimaler Startwert fuer Copula wird ermittelt,
          Iteration",i)
    thwerte[i]<-tt1$par
    valuewerte[i]<-tt1$value
    cat("\n aktueller Copulawert",thwerte[i],
          "und LogFunktionswert",valuewerte[i],"und den
          Startwert",startwerte[i],"\n")
  }
  position <- which.min(valuewerte)

  tt2 <- try(optim(thwerte[position],fn=LogLikli,nugget=nuggetV,
             sill1=sill1,sill2=sill2V,range1=range1V,range2=
             range2V,nue1=nue1V,nue2=nue2V,
             method="BFGS",control=list(maxit=200000)))
  return(tt2)

}

h1<-estCopFrank(th,nugget,sill1,sill2,range1,range2,nue1,nue2)
if(h1$convergence==0) {cat("\n \t Erfolgreiche Berechnung des
                       Copulas \n")}
teh<-h1$par
print(teh)
val<-h1$value
cat("\n\n value")
print(val)


Vnugget<-VerteilungsFunktion(nugget,nugget)
Vsill1<-VerteilungsFunktion(sill1,sill1)
Vrange1<-VerteilungsFunktion(range1,range1)
Vnue1<-VerteilungsFunktion(nue1,nue1)
Vsill2<-VerteilungsFunktion(sill2,sill2)
Vrange2<-VerteilungsFunktion(range2,range2)
Vnue2<-VerteilungsFunktion(varsim[,7],varsim[,7])

hsum<-(Vnue1+Vnue2+Vnugget+Vrange1+Vrange2+Vsill1+Vsill2)+0.0
hprod<-(th^Vnue1-1.0)*(th^Vnue2-1)*(th^Vnugget-1.0)*
       (th^Vrange1-1)*(th^Vrange2-1)*(th^Vsill1-1)*(th^Vsill2-1)
hprod2<-1.0+(1/(th-1)^6)*hprod
l<-log(teh)

erg1<-(720.0*(teh^hsum)*(hprod^6)*(l^6))/(((teh-1)^42)*
      (hprod2^7))
erg2<-(2520.0*(teh^hsum)*(hprod^5)*(l^6))/(((teh-1)^36)*
      (hprod2^6))
```

```
  erg3<-(3360.0*(teh^hsum)*(hprod^4)*(l^6))/(((teh-1)^30)*
        (hprod2^5))
  erg4<-(2100.0*(teh^hsum)*(hprod^3)*(l^6))/(((teh-1)^24)*
        (hprod2^4))
  erg5<-(602.0*(teh^hsum)*(hprod^2)*(l^6))/(((teh-1)^18)*
        (hprod2^3))
  erg6<-(63.0*(teh^hsum)*hprod*(l^6))/(((teh-1)^12)*
        (hprod2^2))
  erg7<-((teh^hsum)*(l^6.0))/(((teh-1)^6)*hprod2)

  ergg<-erg1-erg2+erg3-erg4+erg5-erg6+erg7+0.0

  copula <- ergg

  dichte<-copula*density(nugget,n=length(nugget))$y*
          density(sill1,n=length(sill1))$y*density(range1,n=length
          (range1))$y*density(nue1,n=length(nue1))$y*density(sill2,
          n=length(sill2))$y*density(range2,n=length
          (range2))$y*density(nue2,n=length(nue2))$y
  hilfdichte<-c(rep(0,length(nugget)))
  for(i in 1:length(nugget)) {
    if(dichte[i]>0){
      hilfdichte[i]<-dichte[i]
    }
  }
  dichte<-hilfdichte
  Dichte<-cbind(nugget,sill1,range1,nue1,sill2,range2,nue2,dichte)
  Dichte<-as.data.frame(Dichte)

  return(Dichte)

}
```

A coding for the other copulas can be done in an analogous form, just by interchanging the copula formula in the code. The results for the above derivatives in the codes have been accomplished using the computer algebra system 'Mathematica'.

```
frank[nugget_, sill1_, range1_, nue1_, sill2_, range2_,
      nue2_, \[Theta]_] =
    1/Log[Theta]*
      Log[1 + (1/(([Theta] - 1))^6) ((((Theta^nugget -
                    1)) ((Thetasill1 - 1)) ((Theta^range1 -
                    1)) ((Theta^nue1 - 1)) ((Theta^sill2 -
                    1)) ((Theta^range2 - 1))((Theta^nue2 -
                    1)))))
Dfrank = Derivative[1, 1, 1, 1, 1, 1, 1, 0][frank][nugget, sill1,
    range1, nue1, sill2, range2, nue2, Theta]
```

# 3

# Spectral Representation of Random Fields

The spectral measure associated with every covariance function will become the basic tool for a random field. It is assumed that the random field is m.s. continuous and homogeneous and defined on a finite dimensional Euclidean vector space. In this chapter we will introduce the concept of orthogonal random measure associated with a given positive measure. With that we will be able to show that a m.s. continuous random field is the Fourier Transformation of the random orthogonal measure associated with the spectral measure of its covariance function. This chapter also discuses the concept of random measures and gives mathematical properties for them. Foregoing we will introduce stochastic integrals and define them as linear isometries from a $L_2$ to an other $L_2$, with the help of simple functions. An important theorem shows how the stochastic integral can be evaluated under the use of simple functions and random measures. In analogy to the Karhunen Loeve representation as a main result in this chapter we are going to build a random field $Z$ by a stochastic integral.

## 3.1 Random Measures

In this section we assume $S \subset \mathbb{R}$ and let $\mathcal{B} = \mathcal{B}(S)$ denote the family of its Borel sets. Let a finite positive measure $\lambda : \mathcal{B} \to \mathbb{C}$ be given, as an example $\lambda$ could be the Lebesque's measure:

**Definition 3.1 (Random Measure).** *A random measure in $S$ is a function $\nu : \mathcal{B} \times \Omega \to \mathbb{C}$ such that following conditions hold:*

- $\nu(B, .)$ *is a random variable for each* $B \in \mathcal{B}$
- $\nu(., \omega)$ *is a measure for each* $\omega \in \Omega$

We can deduce that a random measure in $S$ is a random field on $\mathcal{B}$ given by

$$\{\nu(B,.), B \in \mathcal{B}\}$$

and a measure-valued random element

$$\omega \mapsto \nu(.,\omega) \text{ on } (\Omega, \mathcal{A}, P).$$

**Definition 3.2.** *For a random measure for which $\nu(B,.)$ is a square integrable random variable the inner product is given by*

$$E\nu(B,.)\overline{\nu(C,.)} \text{ for } B, C \in \mathcal{B}$$

**Definition 3.3 (Orthogonal random measure).** *An orthogonal random measure in $S$ is a function $\nu : \mathcal{B} \to L_2(\Omega, \mathcal{A}, P)$ such that*

- *$\nu(\sum_{k=1}^{\infty} B_k) = \sum_{k=1}^{\infty} \nu(B_k)$ for every pairwise disjoint sequence in $\mathcal{B}$*
- *there exists a complex measure $\lambda$ on $\mathcal{B}$ such that for every $B, C \in \mathcal{B}$*

$$E\nu(B)\overline{\nu(C)} = \lambda(B \cap C)$$

*for each $B, C \in \mathcal{B}$. In that context $\lambda$ is also called structure measure of $\nu$.*

Note that

$$\lambda(B) = E|\nu(B)|^2 \geq 0$$

is a positive finite measure and particular

$$\lambda(0) = E|\nu(0)|^2 = 0$$

hence $\nu(0) = 0$. The common notation in that theme is $E|\nu(d\xi)|^2 = \lambda(d\xi)$, meaning $E|\nu(B)|^2 = \lambda(B)$ for all $B \in \mathcal{B}$. If a spectral density exists, it is customary to use following notation

$$\lambda(d\xi) = g(\xi)d\xi$$

hence

$$E|\nu(d\xi)|^2 = g(\xi)d\xi$$

In the following $Z$ is a second order random field on $S$, $L_2(Z)$ denotes the closed subspace of $L_2(\Omega, \mathcal{A}, P)$ generated by $\{Z(x), x \in S\}$. $L_2$ is the closure of the set of all linear combination of the form

$$\sum_{k=1}^{n} c_k Z(x_k)$$

for $n \geq 1, x_1, ..., x_n \in S$ and $c_1, ..., c_n \in \mathbb{C}$

**Definition 3.4 (Subordinated).** *An orthogonal random measure $\nu$ is said to be subordinated to the random function $Z$ if $\nu(B) \in L_2(Z)$ for each $B \in \mathcal{B}$.*

## 3.2 Stochastic Integrals

Let $\nu$ be an orthogonal random measure with a structure measure $\lambda$ on $(S, \mathcal{B})$. We are able to construct a linear mapping from $L_2(S, \mathcal{B}, \lambda)$ into $L_2(\Omega, \mathcal{A}, P)$ given by

$$f \mapsto \int f d\nu.$$

This mapping is a linear isometry with the property

$$E|\int f d\nu|^2 = \int |f(x)|^2 \lambda(dx)$$

**Definition 3.5 (Stochastic Integral).** *The linear isometry from $L_2(S, \mathcal{B}, \lambda)$ into $L_2(\Omega, \mathcal{A}, P)$ given by*

$$f \mapsto \int f d\nu.$$

*is called stochastic integral with respect to the orthogonal random measure $\nu$.*

Now we will introduce so called 'simple functions':

**Definition 3.6 (Simple Functions).** *A measurable function $f$ is called simple function if there exists a finite collection of sets, each contained in the domain of $f$ and together covering this domain of definition, such that $f$ assumes a constant finite value on each member of the collection though possible differing from member to member. (see Gikhman and Skorokhod 1969)*

$f$ is a simple function in $L_2(S, \mathcal{B}, \lambda)$ if it is representable in the form

$$f = \sum_{i=1}^{n} a_i \Lambda_{B_k}$$

where $\Lambda_{B_k}$ is the characteristic function of the set $B_k$, $n \geq 1$, $a_1, ..., a_k \in \mathbb{C}$ and $B_1, ..., B_n$ are a finite partition of $S$. An interesting theorem from functional analysis is following

**Theorem 3.7.** *A function $f$ is measurable if and only if the limit of the sequence of simple functions converges everywhere on the given domain.*

Now we can give an alternative definition of the stochastic integral

**Definition 3.8 (Stochastic Integral 2).** *The stochastic integral of $f$ is defined by*

$$\int f d\nu = \sum_{i=1}^{n} a_i \nu(B_i)$$

We can conclude that every simple function is associated with a well-defined random variable in $L_2(\Omega, \mathcal{A}, P)$. So any $f \in L_2(S, \mathcal{B}, \lambda)$ can be approximated by a sequence of simple functions, in the sense that

$$||f - f_n|| = (\int (f(x) - f_n(x))^2 \lambda(dx))^{1/2} \to 0.$$

Because of the construction0 as a linear isometry it follows that

$$E|\int f_m d\nu - \int f_n d\nu|^2 = ||f_m - f_n||^2 \to 0,$$

and this proves that $\{\int f_n d\nu\}_{n=1}^{\infty}$ converges in $L_2(\Omega, \mathcal{A}, P)$ and the limit only depends on $f$.

**Definition 3.9.** *Given $f \in L_2(S, \mathcal{B}, \lambda)$, let $\{f_n\}$ be a sequence of simple functions such that $||f_n - f|| \to 0$. The stochastic integral of $f$ with respect to an orthogonal measure $\nu$ on $(S, \mathcal{B})$ is given by*

$$\int f d\nu = l.i.m_{n \to \infty} \int f_n d\nu$$

The preceding development can be summarised in the following theorem

**Theorem 3.10.** *Let $\nu$ be an orthogonal random measure on $(S, \mathcal{B})$, with structure measure $\lambda$. The stochastic integral with respect to $\nu$ is defined for simple functions $f$ by*

$$\int f d\nu = \sum_i c_i \nu(B_i)$$

*and for arbitrary functions $f \in L_2(S, \mathcal{B}, \lambda)$ by*

$$\int f d\nu = l.i.m_{n \to \infty} \int f_n d\nu.$$

*The correspondence mapping from the function to the integral is a linear isometry from $L_2(S, \mathcal{B}, P)$ into $L_2(\Omega, \mathcal{A}, P)$.*

Now we will construct a new orthogonal random measure from a given one. For an arbitrary function $g \in L_2(S, \mathcal{B}, P)$ and $b \in \mathcal{B}$, the function $g_{\Lambda_B}$ is in $L_2(S, \mathcal{B}, \lambda)$ and therefore the integral

$$\int_B g d\nu = \int \Lambda_B g d\nu$$

is well-defined.

**Definition 3.11.**

$$\nu_g(B) = \int_B g d\nu$$

*for each $B \in \mathcal{B}$*

Above definition defines a function $\nu_g : B \to L_2(\Omega, \mathcal{A}, P)$ and

$$\nu_g(\sum_k B_k) = \sum_k \nu_g(B_k)$$

**Theorem 3.12.** *If $f \in L_2(S, \mathcal{B}, \lambda_g)$, then $fg \in L_2(S, \mathcal{B}, \lambda)$ and*

$$\int f d\nu_g = \int fg d\nu$$

PROOF: For simple $f$

$$\int f d\nu_g = \sum_k c_k \int \Lambda_k g d\nu = \int (\sum_k c_k \Lambda_{B_k}) g d\nu$$

and the statement is true. For a general function $f$, let $\{f_n\}$ be a sequence of simple functions in $L_2(S, \mathcal{B}, \lambda_g)$ converging to $f$. Then

$$\int f_n d\nu_g = \int f_n g d\nu \text{ for } n = 1, 2, ....$$

Since

$$\begin{aligned} E|\int f_n g d\nu - \int f_m g d\nu|^2 &= \int |f_n(x) - f_m(x)|^2 \lambda_g(dx) \\ &= \int |f_n(x)g(x) - f_m(x)g(x)|^2 \lambda(dx) \\ &\to 0 \text{ as } n, m \to \infty \end{aligned}$$

the assertion follows upon taking limits when $n \to \infty$. □
Looking at the Karhunen-Loeve representation for a m.s. continuous random field $Z$ given by

$$Z(x) = \sum_{n=1}^{\infty} \phi_n(x) Z_n$$

we can mention that 'all' the local properties of $Z$ are determined by a single sequence of orthogonal random variables $\{Z_n\}$: it suffices to choose the coefficients appropriately. We will now replace the orthogonal sequence of random variables by an orthogonal random measure $\nu$ on $(S, \mathcal{B})$, where $S$ does not have to be compact, with a structure measure $\nu$. It holds, that

$$E\nu(\mathcal{B}) = 0$$

for each $B \in \mathcal{B}$. In the following we will replace the system of eigenfunctions $\{\phi_n\}$ by a family of Borel functions. But first we will shortly introduce the term Borel.

**Definition 3.13 (Borel Set).** *In a metric space, the $\sigma$-algebra of sets generated by the class $\Sigma$ of open sets is called the $\sigma$-algebra of Borel sets and its elements are called Borel sets.*

Obviously the $\sigma$-algebra generated by the closed sets of a metric space coincides with the $\sigma$-algebra of Borel sets.

**Definition 3.14 (Borel function).** *A function $f(x)$ for $x \in \mathbb{R}$ is called a Borel function if for arbitrary real $a$ the set $\{x : f(x) < a\}$ is a Borel set.*

The family of functions for the replacement is given by

$$\{\phi(.,\xi), \xi \in S\}$$

We assume that $\phi(x,.) \in L_2(S, \mathcal{B}, \lambda)$ for each $x \in S$. It is possible to replace

$$Z(x) = \sum_{n=1}^{\infty} \phi_n(x) Z_n$$

by the stochastic integral

$$\int \phi(x,\xi)\nu(d\xi) = \int \phi(x,.)d\nu$$

By analogy with the Karhunen-Loeve representation we define the random field in the form

$$Z(x) = \int \phi(x,\xi)\nu(d\xi),$$

so $Z$ is a centred second order random field and its covariance functions is given by

$$\begin{aligned} C(x,y) &= E \int \phi(x,\xi)\nu(d\xi) \overline{\int \phi(y,\xi)\nu(d\xi)} \\ &= \int \phi(x,\xi)\overline{\phi(y,\xi)}\lambda(d\xi) \end{aligned}$$

# 4

# Optimal Estimation using Kriging

In the middle of the last century the so called Kriging was developed, named after the African mining engineer Daniel Krige. This development was mainly influenced by the work of mining engineers, geologists, hydrologists and geophysicists and first published by Matheron. This method is now the standard methodology in applied geostatistics science and has been developed in many directions of use. In literature there are many facets of Kriging, namely Simple Kriging, Universal Kriging, Ordinary Kriging, Block Kriging, Co-Kriging and much more, for a detailed look at different Kriging methods we refer to Cressie (1993). In this chapter we will give the basic model assumption for making a prediction at an unobserved location and have a close look at the properties of that prediction. After we introduce the Kriging equations, we look at the properties of this prediction and at the mathematical meaning of Kriging. Kriging will also be identified with Hilbert spaces and proximas. At the end of the chapter we are going to introduce the Bayes Kriging, which was first published in Omre and Halvorsen(1989).

## 4.1 Basic Assumptions

The mean, often called trend function or drift, of the random field $Z(x)$ will be defined as

$$m(x) = EZ(x),$$

where $m(x)$ and $Z(x)$ are defined over all points in the area of interest. The spatial dependence of the random field $Z(x)$ will be modelled by either the covariance function or by the semivariogram. The covariance is defined by

$$\begin{aligned}\mathrm{Cov}(x_1, x_2) &= E[(Z(x_1) - m(x_1))(Z(x_2) - m(x_2))] \\ &= E[Z(x_1)Z(x_2))] - m(x_1)m(x_2)\end{aligned}$$

and the semivariogram by

$$\begin{aligned}\gamma(x_1,x_2) &= \frac{1}{2}\mathrm{Var}[Z(x_1)-Z(x_2)]\\ &= \frac{1}{2}E[(Z(x_1)-Z(x_2))^2]-\frac{1}{2}[(m(x_1)-m(x_2))^2].\end{aligned}$$

The aim of Kriging is to find the best linear unbiased estimation of a linear function of the random field $Z(x)$

- *Linearity* in the sense that the estimator is formed from a linear combination of the observed values.
- *Unbiasedness* This condition requires that the expected value of the estimator have to be equal to the expectation of the true value.
- *Best criterion* is given by the smallest estimation variance.

So our goal is to build a predictor $\hat{Z}(x_0)$ for $Z(x_0)$, where $x_0$ is an unobserved location, where $\hat{Z}(x_0)$ should be linear in the form

$$\hat{Z}(x_0) = \lambda_1 Z(x_1) + ... + \lambda_n Z(x_n)$$

such that

$$E(\hat{Z}(x_0) - Z(x_0))^2 \rightarrow \min_{\lambda \in \mathbb{R}^n} .$$

In order to determine an optimal linear predictor, that means optimal weights $\lambda_i, i = 1, ..., n$ are needed, the assumption of second order stationarity is used.

## 4.2 Kriging of stationary Random Fields

### 4.2.1 Second Order stationarity and the Intrinsic Hypothesis

We will have a stationary random field if following conditions are satisfied

$$\begin{aligned}&EZ(x) = m(x) = \text{ const}\\ &\mathrm{Var}Z(x) = \sigma^2(x) = \text{ const}\\ &\mathrm{Cov}(x_1,x_2) = \mathrm{Cov}(x_1-x_2) \text{ for all } x_1,x_2 \in D\end{aligned}$$

and the Intrinsic Hypothesis holds if

$$\begin{aligned}E[Z(x_1)-Z(x_2)] &= m(x_1-x_2)\\ \mathrm{Var}[Z(x_1)-Z(x_2)] &= 2\gamma(x_1-x_2)\end{aligned}$$

Form the theoretical results we know that second order stationarity implies the intrinsic properties but not vice versa, the relationship is given by

$$\mathrm{Var}[Z(x_1)-Z(x_2)] = \mathrm{Var}[Z(x_1)] + \mathrm{Var}[Z(x_2)] - 2\mathrm{Cov}(x_1,x_2)$$

and in the case of second order stationarity

$$\gamma(x_1 - x_2) = \frac{1}{2}\mathrm{Var}[Z(x_1) - Z(x_2)] = \mathrm{Cov}(0) - \mathrm{Cov}(x_1 - x_2)$$

In the case when the mean $m$ is known and assumed to be constant over the whole domain of the random field the estimator is calculated as the average of the data. This estimator is called Simple Kriging. Simple Kriging is also known in literature as Kriging with known mean. In analogy with multiple regression the weighted average for Simple Kriging is defined as

$$\hat{Z}(x_0) = m + \sum_{i=1}^{n} \lambda_i (Z(x_i) - m)$$

where $\lambda_i$ are weights attached to the residuals $Z(x_i) - m$. By the stationarity assumption the mean is constant at all locations. The estimator is unbiased in the above sense, since

$$\begin{aligned} E[\hat{Z}(x_0) - Z(x_0)] &= m + \sum_{i=1}^{n} \lambda_i (EZ(x_0) - m) - EZ(x_0) \\ &= m + \sum_{i=1}^{n} \lambda_i (m - m) - m \\ &= 0 \end{aligned}$$

The variance of the estimation error is given by

$$\begin{aligned} \sigma_E^2 &= \mathrm{Var}(\hat{Z}(x_0) - Z(x_0)) \\ &= E[(\hat{Z}(x_0) - Z(x_0))^2] \\ &= E[\hat{Z}(x_0)^2 + Z(x_0)^2 - 2\hat{Z}(x_0)Z(x_0)] \\ &= \sum_{i=1}^{n}\sum_{j=1}^{n} \lambda_i \lambda_j C(x_i - x_j) + C(x_0 - x_0) - 2\sum_{i=1}^{n} \lambda_i C(x_i - x_0). \end{aligned}$$

Here we have used the notation $C(x_i - x_j)$ for the term $\mathrm{Cov}(Z(x_i), Z(x_j))$. This expression attains its minimal value by theorem of Fermat at

$$\frac{\partial \sigma_E^2}{\partial \lambda_i} = 0 \text{ for } i = 1, 2, ..., n$$

This leads to the Simple Kriging variance given by

$$\sigma_{SK}^2 = C(0) - \sum_{i=1}^{n} \lambda_i C(x_i - x_0),$$

for details see Wakernagel(1995).

## 4.3 Basic Model Approach

The basic model assumptions in this work are given now, for a detailed look we refer to Chiles and Delfiner (1999), we decompose our random field

$$\{Z(x) : x \in D \subset \mathbb{R}^n\}$$

in

$$Z(x) = m(x) + \varepsilon(x),$$

where $m(x)$ is the trend function as defined above, and typically modelled by a polynomial setup in the form

$$m(x) = \beta_1 f_1(x) + ... + \beta_k f_k(x)$$

The composition is then given by

$$Z(x) = f(x)^T \beta + \varepsilon(x),$$

where $\varepsilon(x)$ is a non-observable random error term with following properties

$$E\varepsilon(x) = 0 \; \forall x \in D \subset \mathbb{R}^n$$
$$\mathrm{Cov}(\varepsilon(x_1), \varepsilon(x_2)) = C(x_1, x_2) \; \forall x_1, x_2 \in D \subset \mathbb{R}^n$$

We can interpret $m(x)$ as the large scale variation of the spatial phenomenon and $\varepsilon(x)$ as the small scale variation.
Let $\mathbf{Z} = (Z(x_1), ..., Z(x_n))$ be a random vector, $\mathbf{F} = (f(x_1), ..., f(x_n))^T$ the design matrix and $\varepsilon = (\varepsilon(x_1), ..., \varepsilon(x_n))^T$ the error term, then the model can be written as

$$\mathbf{Z} = \mathbf{F}\beta + \varepsilon$$

The approach for the prediction of the random variable $Z(x_0)$ is a linear combination of the form

$$\hat{Z}(x_0) = \lambda^T Z$$

where $\lambda = (\lambda_1, ..., \lambda_n)$.

#### Estimation of the mean

The mean value of spatial samples can either be computed by the arithmetic mean or by weighted average, where the weighted average takes into account the knowledge of the spatial correlation. The simplest estimator of the mean is given by

$$\hat{m} = \frac{\sum_{i=1}^{n} Z(x_i)}{n}$$

known by everyone as the arithmetic mean. Since our underlying random field is correlated, a better approach will be an estimator as weighted sum

$$\hat{m} = \sum_{i=1}^{n} w_i Z(x_i)$$

where the $w_i$ are weights with the property that they sum up to 1. To find the best weights, we must solve following problem

$$\min_{\sum_i w_i = 1} \quad \mathrm{Var}(\hat{m} - m)$$

by the Lagrange method.

### Ordinary Kriging

In the case of ordinary Kriging we make the assumption that the trend function $m(x)$ is unknown but constant over the region. Ordinary Kriging is also known in literature as BLUE (best linear unbiased estimator). Our aim is again to estimate an unobserved value of the random field at position $x_0$ by

$$\hat{Z}(x_0) = \sum_{i=1}^{n} w_i Z(x_i) \text{ with } \sum_{i=1}^{n} w_i = 1$$

The estimator of the error variance is given by

$$\begin{aligned} \sigma_E^2 &= \mathrm{Var}(\hat{Z}(x_0) - Z(x_0)) \\ &= C(x_0 - x_0) + \sum_{i=1}^{n}\sum_{j=1}^{n} w_i w_j C(x_i - x_j) - 2\sum_{i=1}^{n} w_i C(x_0 - x_i) \\ &= C(0) + w^T \mathbf{K} W - 2c_0^T w \end{aligned}$$

where

$$\begin{aligned} \mathbf{K} &= (C(x_i - x_j))_{i=1,...,n, j=1,...,n} \\ c_0 &= (C(x_0 - x_1), ..., C(x_0 - x_n))^T \\ w &= (w_1, ..., w_n) \end{aligned}$$

The aim is now to minimise

$$\min_{\sum_i w_i = 1} \quad C(0) + w^T \mathbf{K} w - 2c_0^T w$$

The solution is done by the method of Lagrange multipliers and leads to the ordinary Kriging variance in the case of second order stationarity

$$\sigma_{OK}^2 = C(0) - c_0^T w - \nu,$$

where $\nu$ is the Lagrange parameter. If we assume intrinsic stationarity of the random field $Z(x)$ the Kriging variance is given by

$$\sigma^2_{OK} = w^T \gamma_o - \nu,$$

where $\gamma_0 = (\gamma(x_0 - x_1), ..., \gamma(x_0 - x_n))^T$. One interesting property of ordinary Kriging is, that it is an exact interpolation method when we have no nugget effect (a discontinuity of the variogram at the origin, due to measurement errors and to micro-variability), this means, that

$$\hat{Z}(x_0) = Z(x_k) \Leftrightarrow x_0 = x_k$$

## 4.4 Kriging with non-stationary Random Fields

Under the assumption of universal Kriging, the Krige weights are computed in the way that unbiasedness holds:

$$E\hat{Z}(x_0) = EZ(x_0).$$

This leads to the condition

$$\mathbf{F}^T \lambda = f(x_0),$$

and the mean squared error of prediction should be minimised as above under the setup that there is no knowledge about the trend function. For the estimation of $\beta$ in $\mathbf{Z} = \mathbf{F}\beta + \varepsilon$ generally a BLUE (best linear unbiased estimator) is used in the form

$$\hat{\beta}_{gls} = (\mathbf{F}^T \mathbf{K}^{-1} \mathbf{F})^{-1} \mathbf{F}^T \mathbf{K}^{-1} Z.$$

If we substitute $\beta$ in the simple Kriging predictor by $\beta_{gls}$ we get the universal Kriging predictor in the form

$$\hat{Z}_{UK} = f(x_0)^T \hat{\beta}_{gls} + c_0^T \mathbf{K}^{-1}(Z - \mathbf{F}\hat{\beta}_{gls})$$

with the MSEP

$$C(x_0, x_0) - c_0^T \mathbf{K}^{-1} c_0 + (f(x_0) - \mathbf{F}^T \mathbf{K}^{-1} c_0)^T (\mathbf{F}^T \mathbf{K}^{-1} \mathbf{F})^{-1} (f(x_0) - \mathbf{F}^T \mathbf{K}^{-1} c_0).$$

Summing up there are three classical model assumptions in literature

- Simple Kriging, where the trend function $m(x)$ is known and constant.
- Ordinary Kriging in the case where the trend function $m(x)$ is unknown and constant.
- Universal Kriging when we know nothing about the trend function $m(x)$ and the parameter $\beta$.

We can see that there are two poles in that theory, on the one side we have simple Kriging, where the trend and the parameter $\beta$ is exactly known, and on the other hand we have universal Kriging, where the trend and the parameter $\beta$ is unknown, more on that in later sections.

## 4.5 Thoughts on the best linear prediction

Suppose we have observed a random field $Z(x)$ at $n$ positions $x_1, ..., x_n$ and our aim is to predict at an unobserved location $x_0$. In the case when the law of $Z$ is known, we should use the conditional distribution of $Z(x_0)$ given the observed values $Z(x_1), ..., Z(x_n)$. At the state of art in 'classical' spatial statistics, the calculation of the conditional distribution is very difficult, because of the correlation structure in the field. So it is common to restrict attention to linear predictors, for example Journel and Huijbregts(1978) have written a whole book on linear predictors. Under the assumption that $Z$ has a mean function $m(x)$ and a covariance function $C(x_1, x_2)$ and if they are known it makes sense to apply linear combination of observations of $Z$. We predict $Z(x_0)$ at an unobserved location $x_0$ by a predictor of the form

$$\lambda_0 + \lambda^T \mathbf{Z}$$

The mean squared error of this predictor is the squared mean of the prediction error plus its variance given by

$$[m(x_0) - \lambda_0 - \lambda^T m]^2 + C(0) - 2\lambda^T c_0 + \lambda^T \mathbf{K} \lambda$$

where $m = E\mathbf{Z}$. The squared mean term becomes 0 by taking $\lambda_0 = m(x_0) - \lambda^T m$, so consider choosing $\lambda$ to minimise the variance. For any $\lambda, \nu \in \mathbb{R}^n$ we get:

$$\begin{aligned} \text{Var} \ \ &(Z(x_0) - (\lambda + \nu)^T \mathbf{Z}) = \\ &= c_0 - 2(\lambda + \nu)^T c_0 + (\lambda + \nu)^T \mathbf{K}(\lambda + \nu) \\ &= c_0 - 2\lambda^T c_0 + \lambda^T \mathbf{K}\lambda + \nu^T \mathbf{K} \nu + 2(\mathbf{K}\lambda - c_0)^T \nu \end{aligned}$$

Stein(1999) shows that $c$ lies in the column space of $\mathbf{K}$ and $c$ is orthogonal to the null space of $\mathcal{C}$. So there exists a $\lambda$ such that

$$\mathbf{K}\lambda = c$$

and for such $\lambda$ the inequality

$$\text{Var}(Z(x_0) - (\lambda + \nu)^T \mathbf{Z}) \geq c_0 - 2\lambda^T c + \lambda^T \mathbf{K} \lambda$$

holds for all $\nu$ since $\text{Var}(Z(x_0) - (\lambda + \nu)^T Z) = C(0) - 2\lambda^T c_0 + \lambda^T \mathbf{K}\lambda + \nu^T \mathbf{K} \nu$ and $\nu \mathbf{K} \nu = \text{Var}(\nu^T \mathbf{Z}) \geq 0$. Since $\lambda^T \mathbf{Z}$ achieves this lower bound, we have to minimise the variance of the prediction error. A linear predictor that minimises the mean squared error among all linear predictors is called best linear predictor (BLP).

**Lemma 4.1.** *The BLP is unique in the sense, that if* $\lambda_0 + \lambda^T \mathbf{Z}$ *and* $\nu_0 + \nu^T \mathbf{Z}$ *are both BLPs for* $Z(x_0)$*, then*

$$E(\lambda_0 + \lambda^T \mathbf{Z} - (\nu_0 + \nu^T \mathbf{Z}))^2 = 0$$

If the underlying random field is Gaussian then the conditional distribution is given by

$$N(\lambda_0 + \lambda^T Z, C(0) - c_0^T \mathbf{K} c_0)$$

where $\lambda = \mathbf{K}^{-1} c_0$ and $\lambda_0 = m(x_0) - c_0^T \mathbf{K}^{-1} c_0$.

### 4.5.1 Hilbert Space and Kriging

For studying prediction in the terms of approximation theory in Hilbert spaces, the concept can be seen as projection of an element of a Hilbert space onto a subspace. Let $\mathcal{H}$ denote a Hilbert space and $\mathcal{G}$ a subspace, suppose for a given $h \in \mathcal{H}$, there exists a unique element $g \in \mathcal{G}$ such that

$$||h - g|| = \inf_{g' \in \mathcal{G}} ||h - g'||$$

where $g$ is the projection of $h$ onto $\mathcal{G}$. An important property of the projection $g$ is that it is the unique element in $\mathcal{G}$ satisfying that $h - g$ is orthogonal on $g'$ for all $g' \in \mathcal{G}$. In approximation theory the term 'proximum' is used for $g'$. An interpretation is following: the error of the approximation is orthogonal to all elements in $\mathcal{G}$.

Let $Z$ be a random field defined on $D \subset \mathbb{R}^n$ with mean function $m(x)$ and covariance function $C$. Let $\mathcal{H}_D^0$ be a real linear manifold of $Z(x)$. For $g$ and $h$ in $\mathcal{H}_D^0$ we define the inner product by the expectation

$$\langle g, h \rangle = E(g, h)$$

and the closure of $\mathcal{H}_D^0$ with respect to the inner product is denoted by $\mathcal{H}_D(m, K)$. To characterise BLP, we need to make sure, that constant terms are in this space of possible predictors. Let $Q$ be a set on which the random field $Z$ is observed, all linear predictors of $h \in \mathcal{H}_D(m, K)$ are of he form

$$c + g,$$

where $c$ is a scalar and $g \in \mathcal{H}_Q(m, K)$.

Let $g(h)$ be a unique element in $\mathcal{H}_Q(m, K)$ satisfying

$$\text{Cov}(h - g(h), g') = 0 \; \forall g' \in \mathcal{H}_Q(m, K)$$

and set $c(h) = Eh - Eg(h)$. Then $c(h) + g(h)$ is the BLP of $h$, which follows from

$$E((h - c(h) - g(h))(c' + g')) = 0 \; \forall c', g' \in \mathcal{H}_Q(m, K)$$

For more details we refer to Deutsch(2001).

### 4.5.2 Best Linear unbiased prediction

Best linear unbiased prediction is called Kriging in geostatistical literature as mentioned above. If the trend function $f(x) \equiv 1$, so the mean of the process is assumed to be an unknown constant, then best linear unbiased prediction is called ordinary Kriging as introduced in earlier sections. For more general $m(x)$ it is known as universal Kriging, and best linear prediction is called simple Kriging when $m(x)$ is known or 0. A useful characterisation of BLP is that its error is orthogonal to all possible linear predictors (uncorrelated). The BLUP has a similar characterisation. Suppose for a given random field

$$Z(x) = f(x)^T \beta + \varepsilon(x)$$

The random variable

$$\sum_{i=1}^{n} w_i Z(x_i)$$

is called a contrast if it has mean 0 for all $\beta$, or equivalently if $\sum_{i=1}^{n} w_i m(x_i) = 0$. A BLUP has the property that its error is a contrast and its error is orthogonal to any contrast of the observations. To construct a BLUP it is necessary to know the covariance structure of all contrasts of the random field, but this concept is of value, when using intrinsic random fields.

**Lemma 4.2.** *If a LUP exists, then the BLUP exists and is unique.*

The Bayesian interpretation of BLUP can be found in Omre(1987), in this context suppose that $Z$ is given in the form

$$Z(x) = f(x)^T \beta + \varepsilon(x)$$

where $\varepsilon(x)$ is Gaussian and independent of the random vector $\beta$, which has a prior distribution of the form $N(\mu, \sigma^2 V)$, for simplifications in notation we define

$$W(\sigma^2) = (\mathbf{F}^T \mathbf{K}^{-1} \mathbf{F} + \sigma^{-2} \mathbf{V}^{-1})^{-1}$$

where $\mathbf{F} = (f(x_1), ..., f(x_n))^T$, then the posterior distribution of $\beta$ given $Z$ is

$$\beta | Z \sim N(W(\sigma^2)(\mathbf{F}^T \mathbf{K}^{-1} \mathbf{Z} + \sigma^{-2} V^{-1}), W(\sigma^2)).$$

The predictive distribution of $Z(x_0)$ given $\mathbf{Z}$ is given by

$$\begin{aligned} Z(x_0)|\mathbf{Z} \sim N( & c_0^T \mathbf{K}^{-1} \mathbf{Z} + \gamma^T W(\sigma^2)(\mathbf{F} \mathbf{K}^{-1} \mathbf{Z} + \sigma^{-2} V^{-1} \mu), \\ & C(0) - c_0^T \mathbf{K}^{-1} c_0 + \gamma^T W(\sigma^2) \gamma), \end{aligned}$$

where $\gamma$ is given by $f(x_0) - \mathbf{F}^T \mathbf{K}^{-1} k$. $k$ is given by $\mathrm{Cov}(\mathbf{Z}, Z(x_0))$
Letting $\sigma^2 \to \infty$ means letting the prior on $\beta$ gets increasingly uninformative, and the limiting predictive distribution of $Z(x_0)$ given $\mathbf{Z}$ is Gaussian with the BLUP as its conditional expectation and the conditional variance is given by

$$C(0) + c_0^T \mathbf{K}^{-1} c_0 + \gamma^T (\mathbf{F}^T \mathbf{K}^{-1} \mathbf{F})^{-1} \gamma$$

which is the mean squared error of the BLUP.

## 4.6 Bayes Kriging

As we have mentioned above there is a bridge between simple Kriging and universal Kriging. Between these two poles there is the so called Bayes Kriging, first introduced by Omre and Halvorsen (1989), the approach for that model is a decomposition of the random field as mentioned in earlier chapters in

$$Z(x) = m(x) + \varepsilon(x)$$

with

$$\begin{aligned} &\mathrm{E}(\varepsilon(x)) = 0 \ \forall x \in D \subset \mathbb{R}^n \\ &\mathrm{Cov}(\varepsilon(x), \varepsilon(y)) = C(x, y) \ \forall x, y \in D \subset \mathbb{R}^n \end{aligned}$$

and the trend function $m(x)$ is modelled by

$$\begin{aligned} m(x) &= f_1(x)\beta_1 + ... + f_n(x)\beta_n \\ &= f(x)^T \beta \end{aligned}$$

where $\beta \sim P$. Since in practice we always know something of the spatial phenomena we have some a priori distribution over $\beta$, namely $P_\beta$. Under the assumption that the first two moments of this distribution are known, we can build the Bayes predictor.

$$\begin{aligned} \mathrm{E}(\beta) &= \mu \\ \mathrm{Var}(\beta) &= \Phi \end{aligned}$$

The Bayes Kriging predictor for an unobserved location $x_0$ is then given by

$$\hat{Z}_B(x_0) = f(x_0)^T \mu + (c_0 + \mathbf{F}\Phi f(x_0))^T (\mathbf{K} + \mathbf{F}\Phi\mathbf{F}^T)^{-1} (\mathbf{Z} - \mathbf{F}\mu)$$

and the total mean squared error of prediction is given by

$$\begin{aligned} \mathrm{TMSEP}(\hat{Z}_B) = \ & C(x_0, x_0) + f(x_0)^T \Phi f(x_0) \\ & -(c_0 + \mathbf{F}\Phi f(x_0))^T (\mathbf{K} + \mathbf{F}\Phi\mathbf{F}^T)^{-1} (c_0 + \mathbf{F}\Phi f(x_0)) \end{aligned}$$

In the case of $\Phi = 0$, i.e. 'perfect' knowledge of the trend, we have simple Kriging and corresponding to $\Phi^{-1} = 0$, i.e. 'nothing' is known a priori about the trend, we have universal Kriging. Omre and Halvsen also have shown that $\hat{Z}_B(x_0)$ minimises the total MSEP, i.e. the MSE averaged with respect to the trend parameter $\beta$.

# 5

# Covariance and Variogram Structure

From the chapter on analytical properties we know that a random field is fully described by the trend function and the covariance function, which is in connection to the variogram function as seen in earlier chapters. In following chapter we will shortly explain different covariance and variogram models and also look at their estimation. We are going to discuss the behaviour of the random field near the origin with the link to the covariance/variogram function and introduce spherical and derived models. After giving a detailed overview about commonly used covariance/variogram functions we simulate random fields with that function. For the simulation we are going to construct a dense grid on a domain and simulate at each point a realisation of the random field. We are give give the code for the simulation using the R package 'RandomFields'.

## 5.1 Covariance and Variogram Functions

As we already know the spatial dependency of the observations is modelled by the covariance function

$$C(x_i, x_j) = \mathrm{Cov}(Z(x_i), Z(x_j))$$

where $x_i$ and $x_j$ are in the given domain $D \subset \mathbb{R}^n$. For estimation of the covariance function the so called semivariogram is used,

$$\gamma(x_i, x_j) = \frac{1}{2}\mathrm{Var}(Z(x_i) - Z(x_j)),$$

and under the assumption of stationarity it follows that

$$\begin{aligned}\gamma(h) &= \gamma(x, x+h) \\ &= \frac{1}{2}E(Z(x+h) - Z(x))^2\end{aligned}$$

$$
\begin{aligned}
&= \frac{1}{2}(E(Z(x+h)-m)^2 + E(Z(x)-m)^2) - \\
&\quad E((Z(x+h)-m)(Z(x)-m)) \\
&= C(0) - C(h)
\end{aligned}
$$

as seen in the chapters before and it holds that

$$C(0) = \lim_{h \to \infty} \gamma(h)$$

where $h$ is the distance between the observations. For the estimation of the variogram there are following estimators in the literature, first the classical empirical semivariogram introduced by Matheron,

$$\hat{\gamma}(h) = \frac{1}{2|N(h)|} \sum_{N(h)} (Z(x_i) - Z(x_j))^2,$$

where $N(h)$ is defined by the set

$$N(h) = \{(i,j) : x_i - x_j = h\}$$

and $|N(h)|$ gives the number of pairs in $N(h)$. Cressie and Hawkins(1980) gave a robust version for the estimation in the form

$$\hat{\gamma}(h) = \frac{(\frac{1}{|N(h)|}\sum_{N(h)} |Z(x_i) - Z(x_j)|^{1/2})^4}{2(0.457 + \frac{0.494}{|N(h)|})}.$$

In literature there is not a large number of other robust approaches for the estimation of the variogram structure, so we hope that the interesting reader will invest work on that theme, or we will do this in coming papers.

## 5.2 Behaviour near the Origin

An interesting thing in research is the behaviour of the covariance of the variogram near the origin, this behaviour is linked to the continuity and to the spatial regularity of the regionalized variable. We can characterise following typical behaviours:

- A parabolic behaviour, this characterises a highly regular random field that is usually differentiable at least piecewise.
- A linear behaviour, the random field is continuous, at least piecewise, but not differentiable. It is less regular than in the parabolic case.
- A discontinuity at 0, nugget effect. The random field is generally not continuous and is thus very irregular.

The theoretical variogram function is then fitted to the empirical estimated variogram function. We can use ordinary least squares for the fitting in the form

$$Q(b) = \sum_{i=1}^{m} [\hat{\gamma}(h_i) - \gamma(h_i, b)]^2$$

where $b$ is chosen such, that $Q(b)$ is minimised. Or we can apply generalised least squares, we take into account the correlations between the different values of the sample variogram, we choose $b$ in the way that

$$Q(b) = [\hat{\gamma} - \gamma(b)]^T \mathbf{V}^{-1} [\hat{\gamma} - \gamma(b)]$$

is minimised. Here $\hat{\gamma}$ is the vector of elements $\hat{\gamma}(h_i)$,$\gamma(b)$ the vector of elements $\gamma(h_i, b)$ and $\mathbf{V}$ is the variance-covariance matrix of $\gamma(\hat{h}_i)$. A compromise between the efficiency of generalised least squares and the simplicity of ordinary least squares is weighted least squares, namely minimisation of

$$Q(b) = \sum_{i=1}^{m} w_i^2 [\hat{\gamma}(h_i) - \gamma(h_i, b)]^2$$

The weights $w_i$ are equal to the reciprocals of $\mathrm{Var}[\gamma(\hat{h}_i)]$. In literature we can also find the maximum likelihood approach, this method is implemented in geoR, Bayesian and Fuzzy fitting methods. In following we will present some of the typical covariance and variogram models that are used in practise. The covariances and variograms are presented as functions of $r = |h|$ and the corresponding spectra as functions of $\rho = |u|$.

## 5.3 Spherical Models and Derived Models

By the autoconvolution of the indicator function of the sphere of $\mathbb{R}^n$ with diameter $a$, namely in terms of the modulus $\zeta = |x|$, of the function

$$w_n(\zeta) = \begin{cases} 1 & \text{if } \zeta \leq \frac{a}{2} \\ 0 & \text{if } \zeta > \frac{a}{2} \end{cases}$$

we obtain the spherical covariogram of $\mathbb{R}^n$ as a function of $r = |h|$ as

$$g_n(r) = \begin{cases} a^n v_{n-1} \int_{r/a}^{1} (1 - u^2)^{\frac{n-1}{2}} du & \text{if } r \leq a \\ 0 & \text{if } r \geq a \end{cases}$$

where $v_n$ is the volume of the unit ball of$\mathbb{R}^n$ given by

$$v_n = \frac{\pi^{n/2}}{2^{n-1} n \Gamma(\frac{n}{2})}$$

This model is used in practice corresponding to $n = 1, 2, 3$.

- Triangle model
  This model is also known under the name tent covariance and is only valid in $\mathbb{R}^1$

$$C_1(r) = \begin{cases} 1 - \frac{r}{a} & \text{if } r \le a \\ 0 & \text{if } r \ge a \end{cases}$$

- Circular model
  This model is only valid in $\mathbb{R}^2$

$$C_2(r) = \begin{cases} \frac{2}{\pi}(\text{arc } \cos(\frac{r}{a}) - \frac{r}{a}(1 - \frac{r^2}{a^2})^{1/2}) & \text{if } r \le a \\ 0 & \text{if } r \ge a \end{cases}$$

- Spherical model
  This model is only valid in $\mathbb{R}^3$

$$C_3(r) = \begin{cases} 1 - \frac{3r}{2a} + \frac{r^3}{2a^3} & \text{if } r \le a \\ 0 & \text{if } r \ge a \end{cases}$$

The corresponding variograms have a linear behaviour near the origin and reach their sill at $r = a$, this means the scale parameter $a$ coincides with the range $r$. More regular models corresponding to $q$ times m.s. differentiable second order stationary random field are obtained by the Radon transformation (see Chiles and Delfiner(1999) for more details) of every even order. The Radon transformation of order $2q$ of the indicator function of $\mathbb{R}^n$ is given by

$$w_{n,2q}(\zeta) = \begin{cases} v_{2q}(a^2 - 4\zeta^2)^q & \text{if } \zeta \le \frac{a}{2} \\ 0 & \text{if } \zeta \ge \frac{a}{2} \end{cases}$$

as a function of $\zeta = |x|\, x \in \mathbb{R}^{n-2q}$ and it represents the volume of the sphere of $\mathbb{R}^{2q}$ with radius $(\frac{1}{4}a^2 - \zeta^2)$ and only depends on $q$. Looking at the covariogram shows that $g_{n,2}$ can be expressed in terms of $g_n$ and $g_{n+1}$ as

$$g_{n,2}(r) = \pi((a^2 - r^2)g_n(r) - \frac{v_{n-1}}{v_{n+1}} g_{n+2}(r))$$

where we have used the fact that $g_{n,2}(r) = 2\pi \int_r^\infty u g_n(u) du$. By using the concept of recursion $g_{n,2q}$ can be expressed in terms of $g_n, g_{n+2}, ..., g_{n+2q}$ in the form

$$g_{n,2q} = \frac{\pi^q}{q!} \sum_{i=1}^{q} (-1)^i \binom{q}{i} \frac{v_{n-1}}{v_{n+2l-1}} (a^2 - r^2)^{q-i} g_{n+2i}(r) \text{ for } r \le a$$

and $g_{n,2q}(r) = 0$ for $r \ge a$. These models are valid in $\mathbb{R}^{n-2q}$. In $\mathbb{R}^3$ we obtain following model as an example

- Cubic model, for the choice of $q = 1$ and $n = 5$ with the lowest irregular term $r^3$.

$$C(r) = \begin{cases} 1 - \frac{7r^2}{a^2} + \frac{35r^3}{4a^3} - \frac{7r^5}{2a^5} + \frac{3r^7}{4a^7} & \text{if } r \le a \\ 0 \text{ if } r \ge a \end{cases}$$

## 5.4 Behaviour of $\gamma(\boldsymbol{h})/|\boldsymbol{h}|$

The variogram of a stationary random field is finite, the variogram of a non stationary random field can increase to infinity, Yaglom (1987) gives following property

- If $\gamma(h)$ is the variogram of a m.s. continuous intrinsic random field,

$$\frac{\gamma(h)}{|h|^2} \to 0 \text{ as } |h| \to \infty$$

- If $\gamma(h)$ is the variogram of a m.s. differentiable intrinsic random field then there exists a constant $A$ with following property

$$\gamma(h) \leq A|h|^2 \quad \forall h \in \mathbb{R}^n$$

## 5.5 Exponential Model and Derived Models

The exponential model with scale parameter $a > 0$ is defined by

$$C(r) = e^{-\frac{r}{a}},$$

and it is a covariance in $\mathbb{R}^n$ $n \geq 1$, since it corresponds to the positive spectral density

$$f_n(\rho) = \frac{2^n \pi^{(n-1)/2} \Gamma(\frac{n+1}{2}) a^n}{(1 + 4\pi^2 a^2 \rho^2)^{(n+1)/2}}.$$

For more details see Yaglom(1987). The variogram reaches its sill only asymptotically when $r$ goes to infinity and the practical range (95 percent of the sill or equivalently a correlation of 5 percent) is approximately at $3a$. If $n = 1$ the exponential model is the covariance of a time continuous Markov process. (see Gikhman and Skorokhod(1996)).

## 5.6 Gaussian Model

The Gaussian model with scale parameter $a > 0$ is defined by following formula

$$C(r) = e^{-\frac{r^2}{a^2}}$$

it is a covariance in $\mathbb{R}^n$ for any $n$. Its practical range is about $1.73a$. This model is associated with an infinitely differentiable second order random field and thus is extremely regular. Also after a Radon transformation of oder $m$ the covariance remains Gaussian. It is possible to evaluate all partial derivatives of all orders.

## 5.7 Covariance Models Overview

Now we want to present a short overview over commonly used covariance models. We illustrate the models by a graphical plot applied to a given simulated random field on a grid. On that grid we simulated a random field under the usage of the R package 'RandomField' by Martin Schlather. We used a Gaussian random field with mean equal 1, variance 2, nugget 1 and scale value equal 2.

```
step <- 0.1
x <- seq(0,20,step)
y <- seq(0,20,step)
```

We used following function for evaluation of a grid

```
gridMaker <- function(Xmin,Xmax,Ymin,Ymax,nX,nY) {
  laengeX <- Xmax-Xmin
  laengeY <- Ymax-Ymin
  coord <- NULL
  for(i in 1:nX) {
    for(j in 1:nY) {
      coord1 <-c(Xmin+laengeX/nX*i,Ymin+laengeY/nY*j)
      coord <- rbind(coord,coord1)
    }
  }
}
```

and get following dense grid with 40000 points, we used such an amount for a more 'beautiful' presentation of the pictures. The random field is then simulated by

```
f <- GaussRF(x=x, y=y, model=model, grid=TRUE, param=param )
```

where in the variable model the underlying covariance model is saved and in the variable parameter the covariance parameters are saved. The empirical variogram and the true variogram and covariance function for that random field have been evaluated by

```
bins <- seq(0,10,0.5)
binned <- EmpiricalVariogram(x=x, y=y, data=f, grid=TRUE, gridTriple=FALSE,
          bin=bins)

trueVariogram <- Variogram(binned$c$, model=model, param=param)
trueCov <- CovarianceFct(binned$c$, model=model, param=param)
```

### Gaussian Variogram and Covariance Function

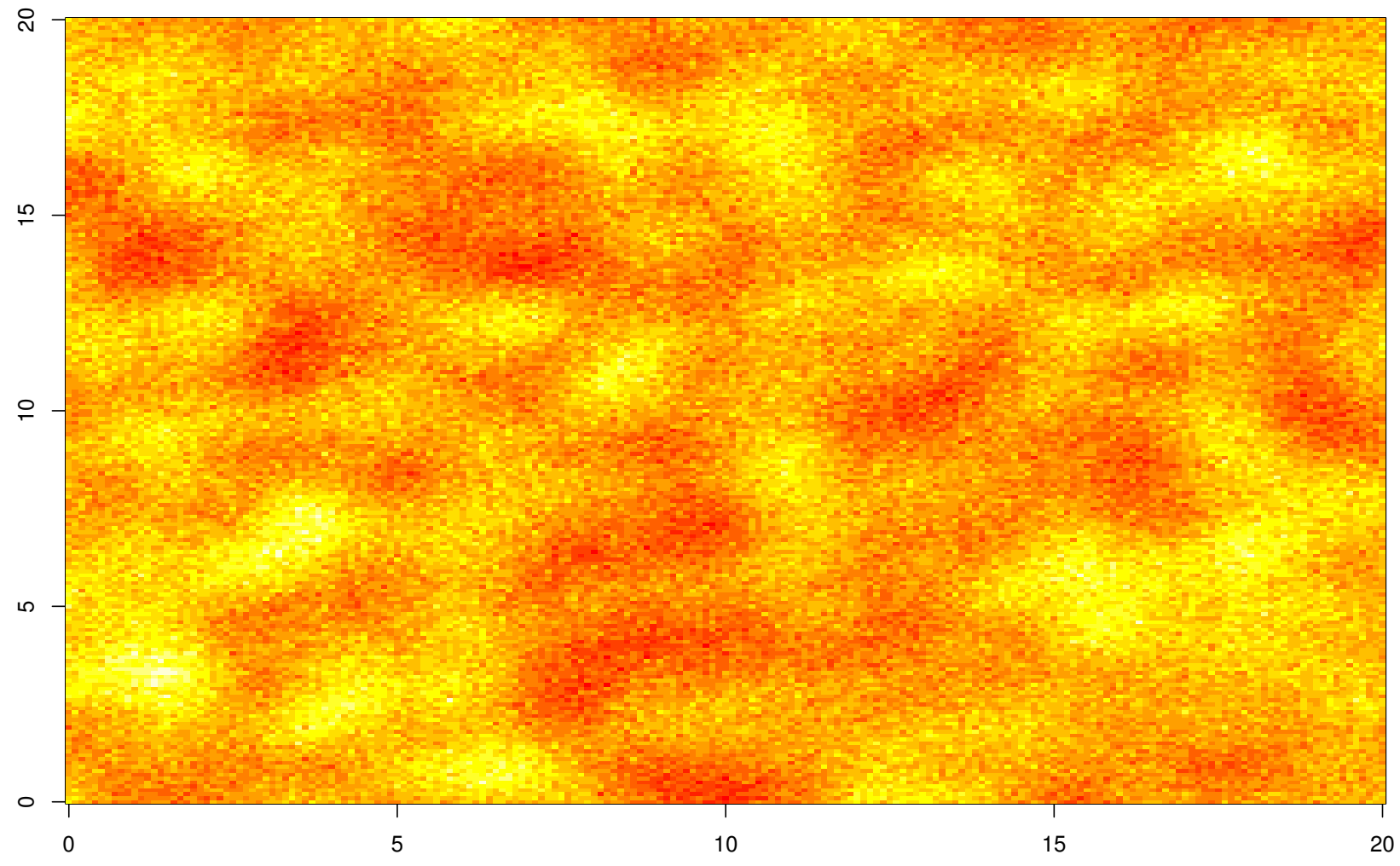

**Fig. 5.1.** Simulated Gaussian random field with Gaussian covariance function

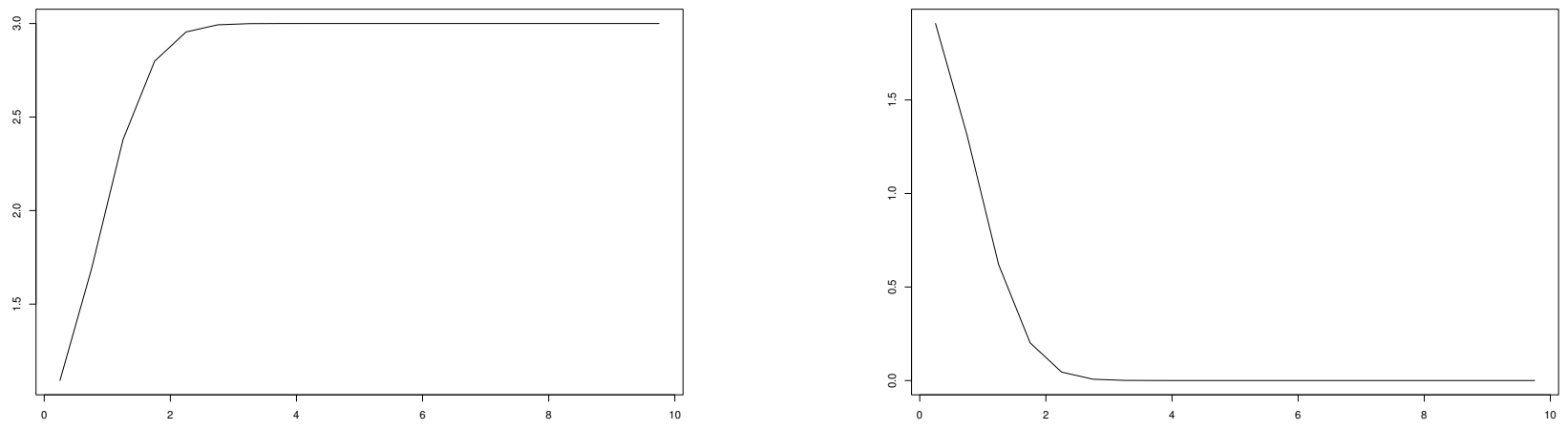

**Fig. 5.2.** Variogram and covariance function of a Gaussian random field

### Cauchy Variogram and Covariance Function

$$C(r) = (1 + \frac{r^2}{a^2})^{-\alpha}$$

where $a > 0$ and $\alpha > 0$. This model is valid in $\mathbb{R}^n$ for $n \geq 1$ since it can be expressed in the form

$$C(r) = \int_0^{\infty} e^{-\frac{r^2}{t^2}} \mu(dt)$$

where $\mu$ is an arbitrary bounded positive measure (see Yaglom(1987)). This model is very regular near the origin and reaches its sill slowly.

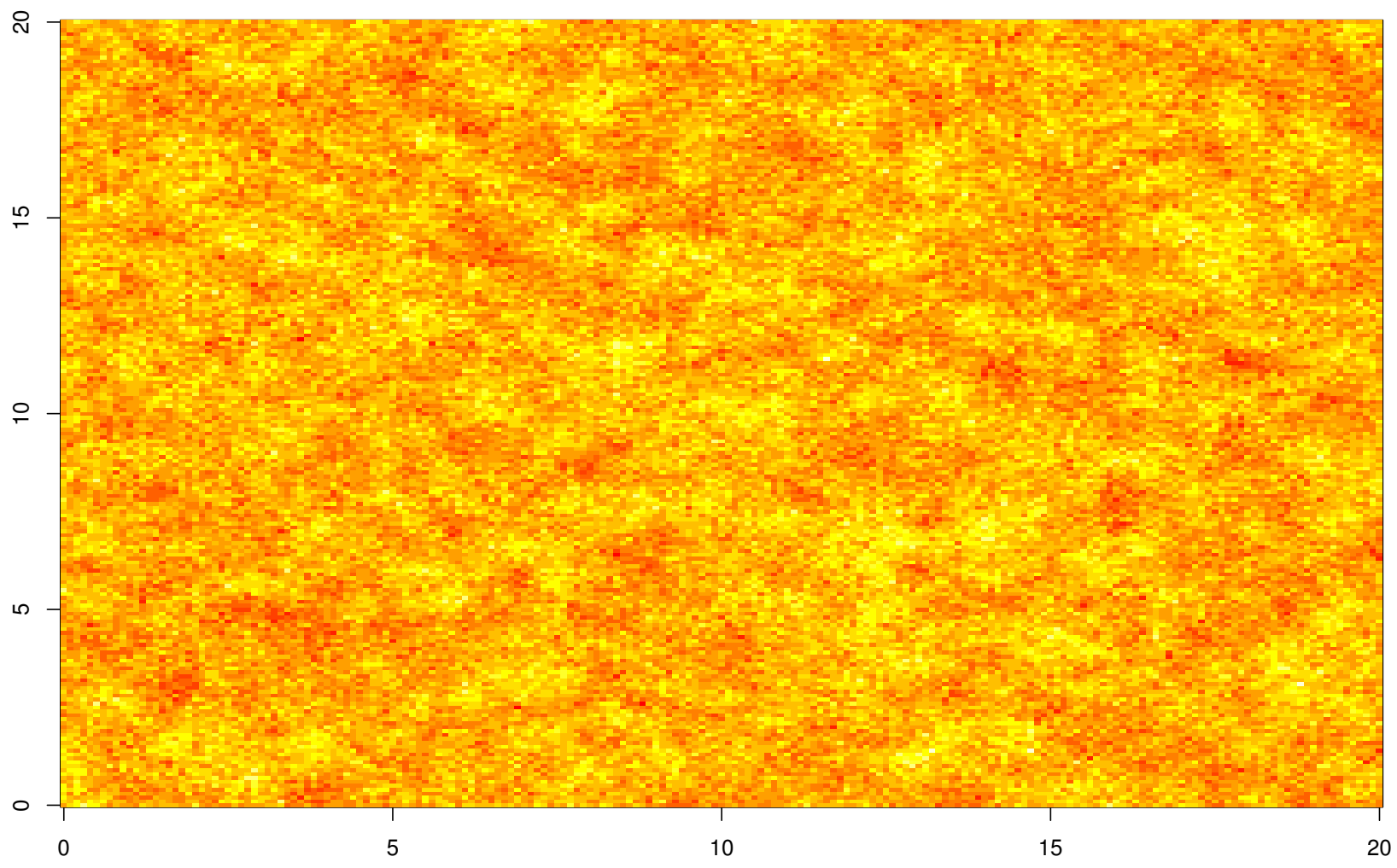

**Fig. 5.3.** Simulated Gaussian random field with Cauchy covariance function

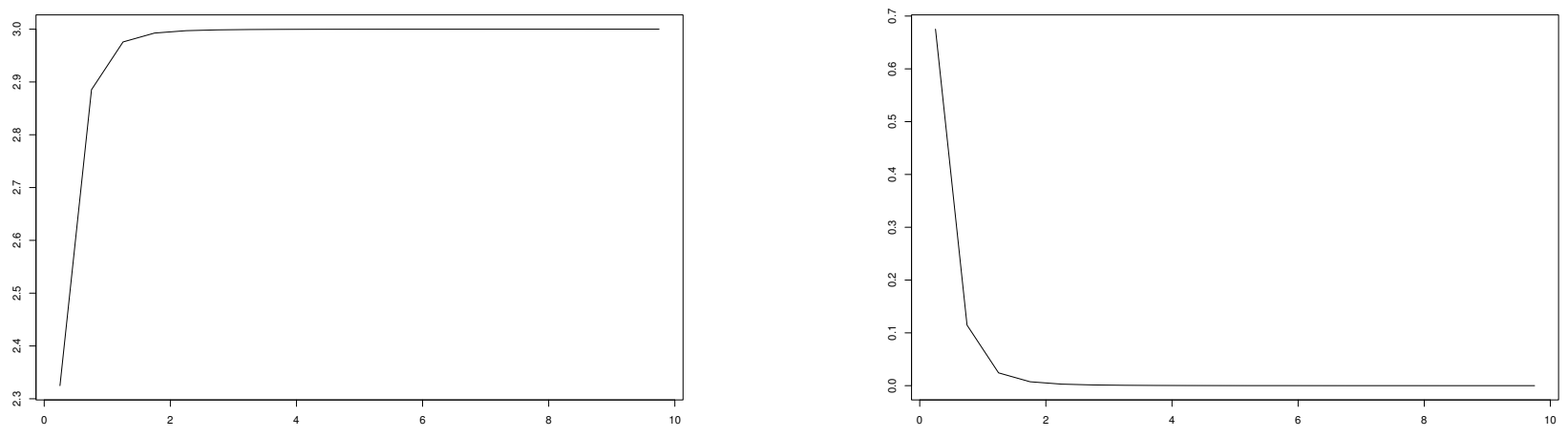

**Fig. 5.4.** Variogram and covariance function of a Gaussian random field with Cauchy Covariance function

### Circular Variogram and Covariance Function

$$C(x) = \begin{cases} 1 - \frac{2}{\pi}(x(1-x^2)^{\frac{1}{2}} + \text{acrsin}(x)) & \text{if } 0 \leq x \leq 1 \\ 0 & \text{else} \end{cases}$$

This is an anisotropic covariance function and is only valid for $n \leq 2$.

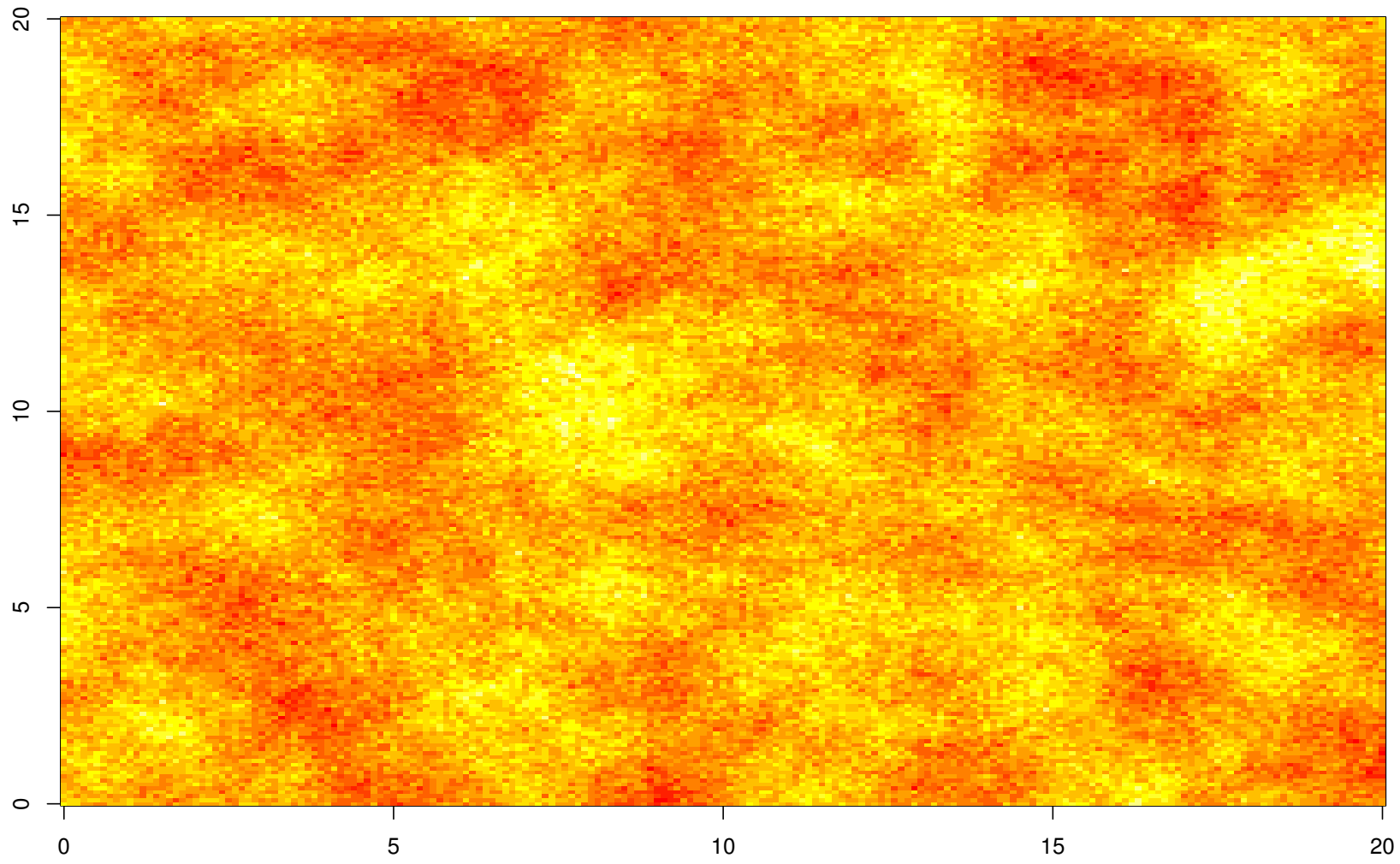

**Fig. 5.5.** Simulated Gaussian random field with circular covariance function

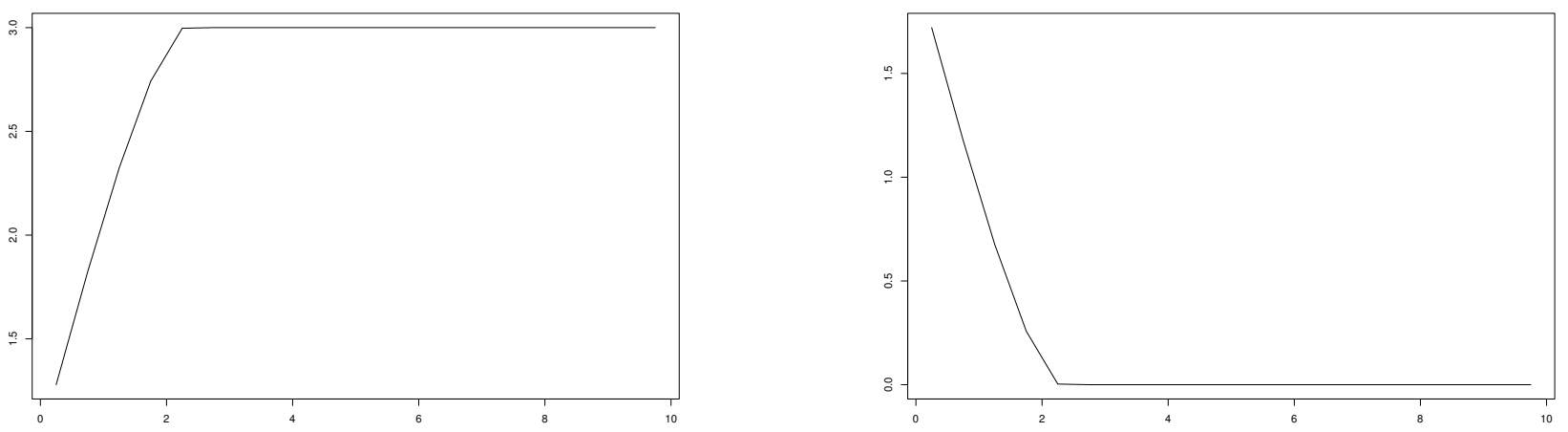

**Fig. 5.6.** Variogram and covariance function of a Gaussian random field with circular covariance function

### Cubic Variogram and Covariance Function

$$C(x) = \begin{cases} 1 - 7x^2 + 8.75x^3 + 3.5x^5 + 0.75x^7 & \text{for } 0 \leq x \leq 1 \\ 0 & \text{else} \end{cases}$$

This is a two times differentiable covariance function with a compact support, but the cubic covariance function is only valid for $n \leq 3$.

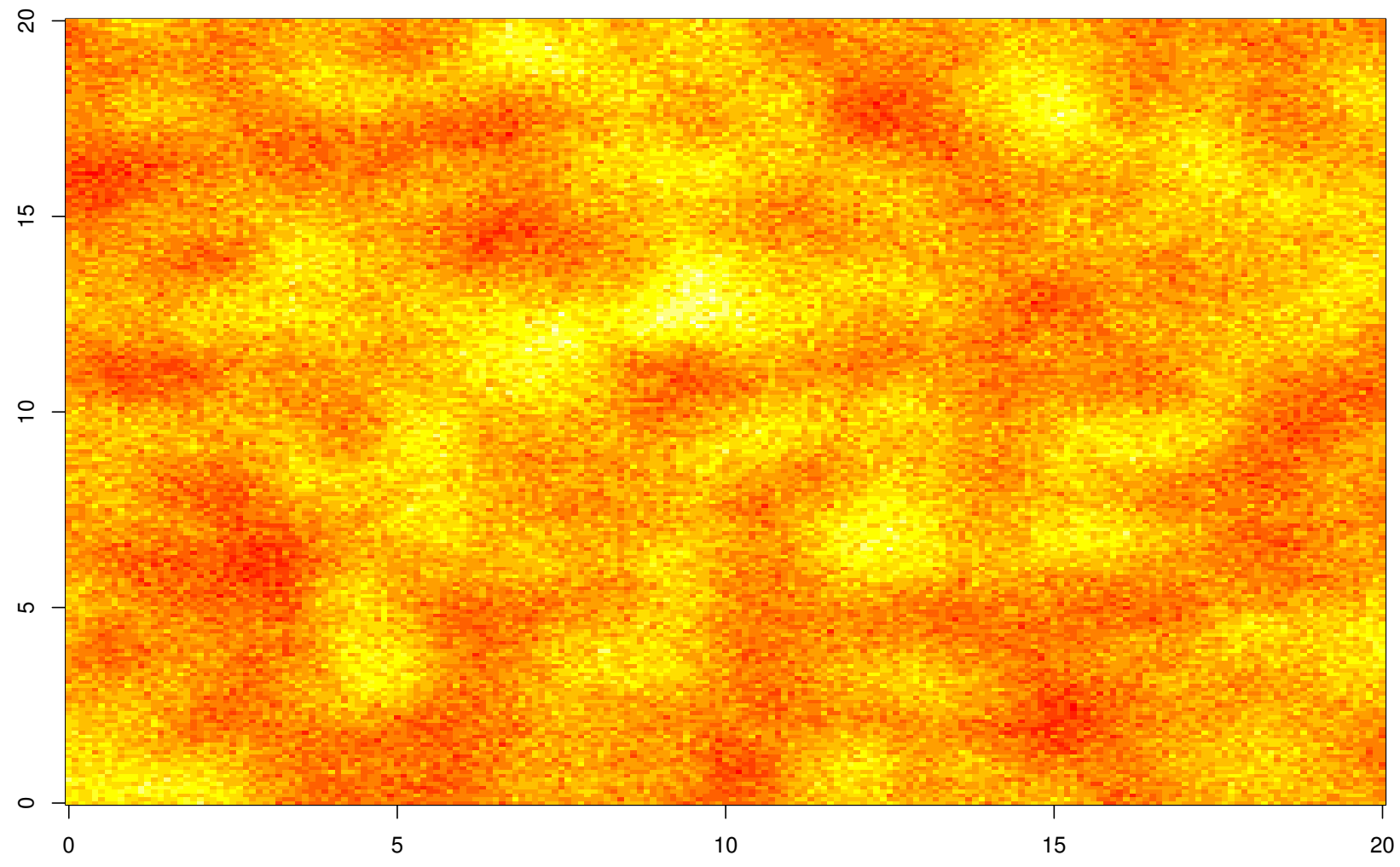

**Fig. 5.7.** Simulated Gaussian random field with cubic covariance function

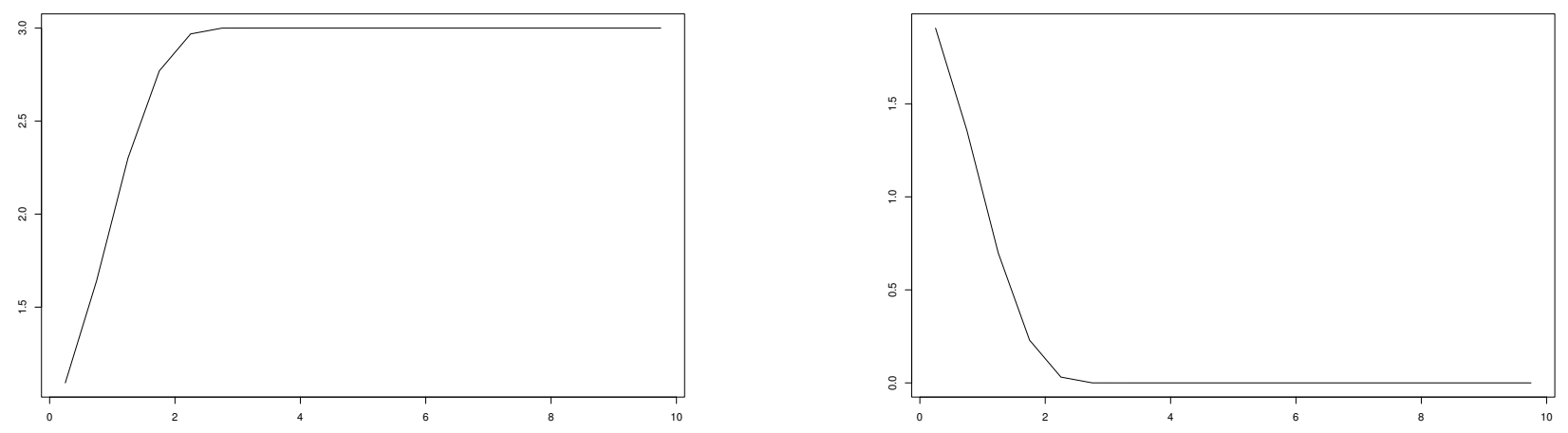

**Fig. 5.8.** Variogram and covariance function of a Gaussian random field with cubic covariance function

## Exponential Variogram and Covariance Function

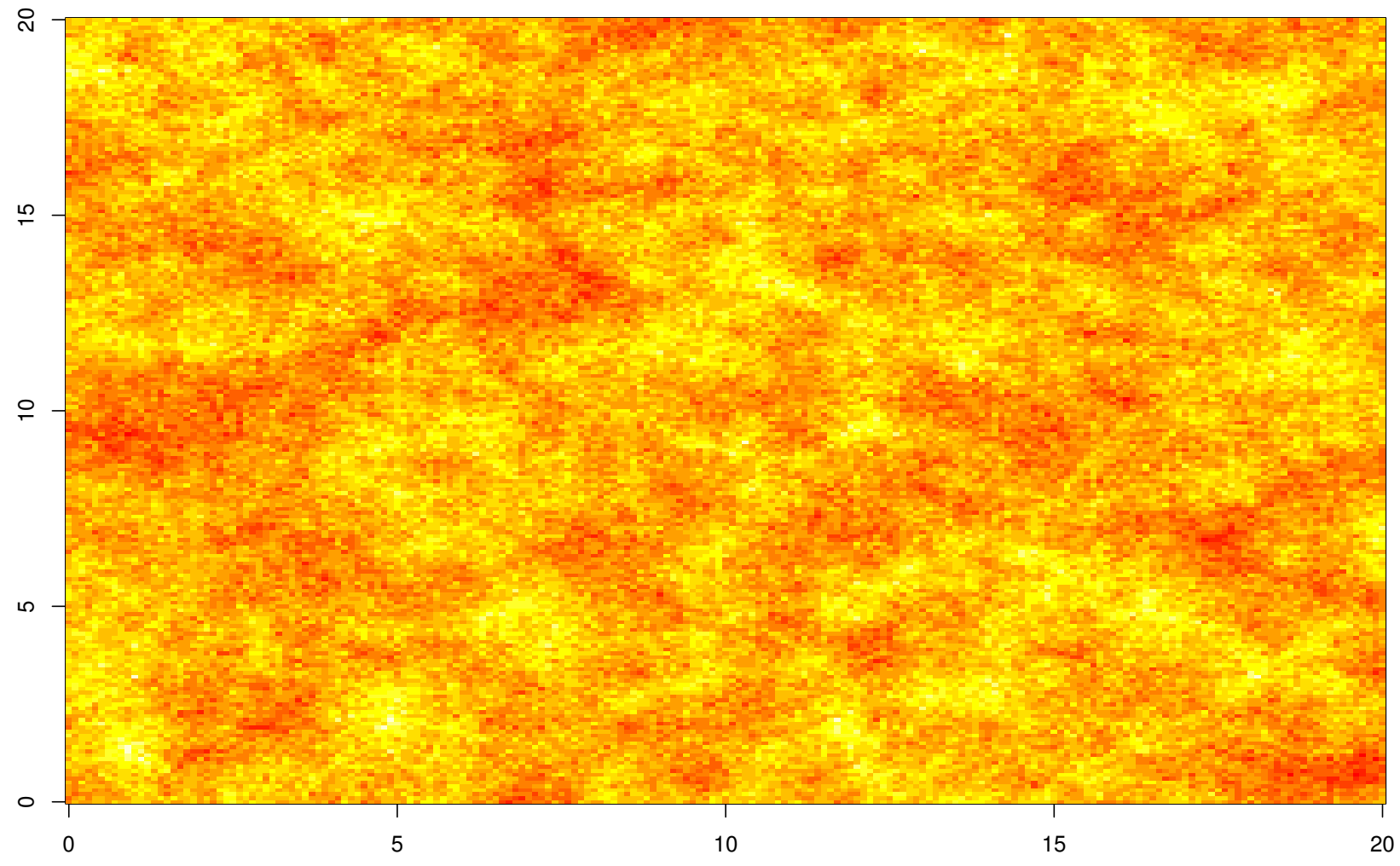

**Fig. 5.9.** Simulated Gaussian random field with exponential covariance function

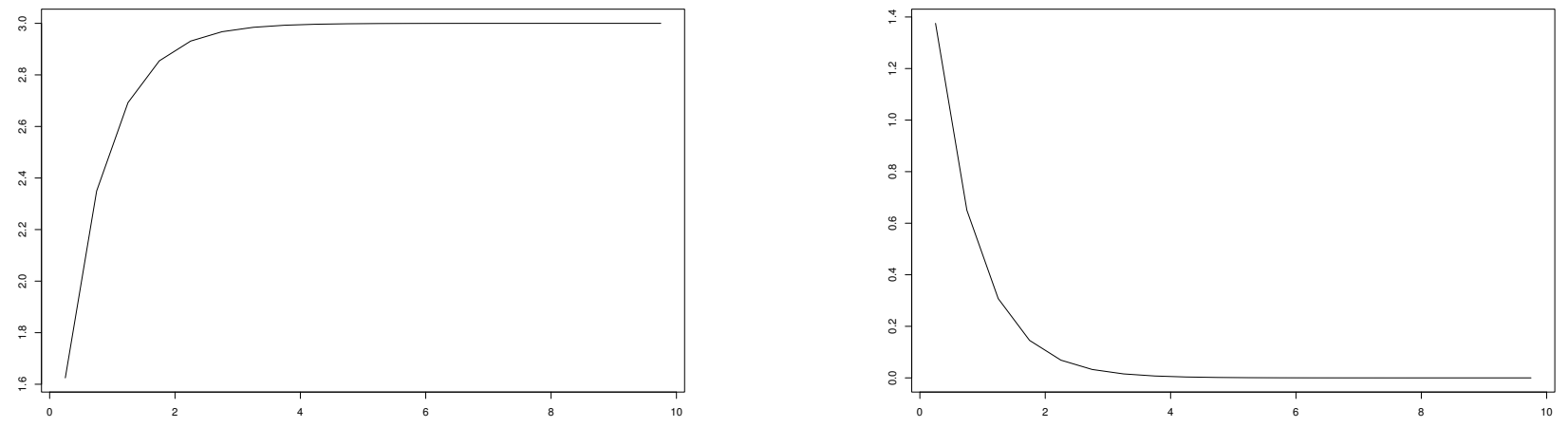

**Fig. 5.10.** Variogram and covariance function of a Gaussian random field with exponential covariance function

### Penta Variogram and Covariance Function

$$C(x) = \begin{cases} 1 - \frac{22}{3}x^2 + 33x^4 - \frac{77}{2}x^5 + \frac{33}{2}x^7 - \frac{11}{2}x^9 + \frac{5}{6}x^{11} & \text{for } 0 \leq x \leq 1 \\ 0 & \text{else} \end{cases}$$

This is a 4 times differentiable covariance function with a compact support, only valid for $n \leq 3$.

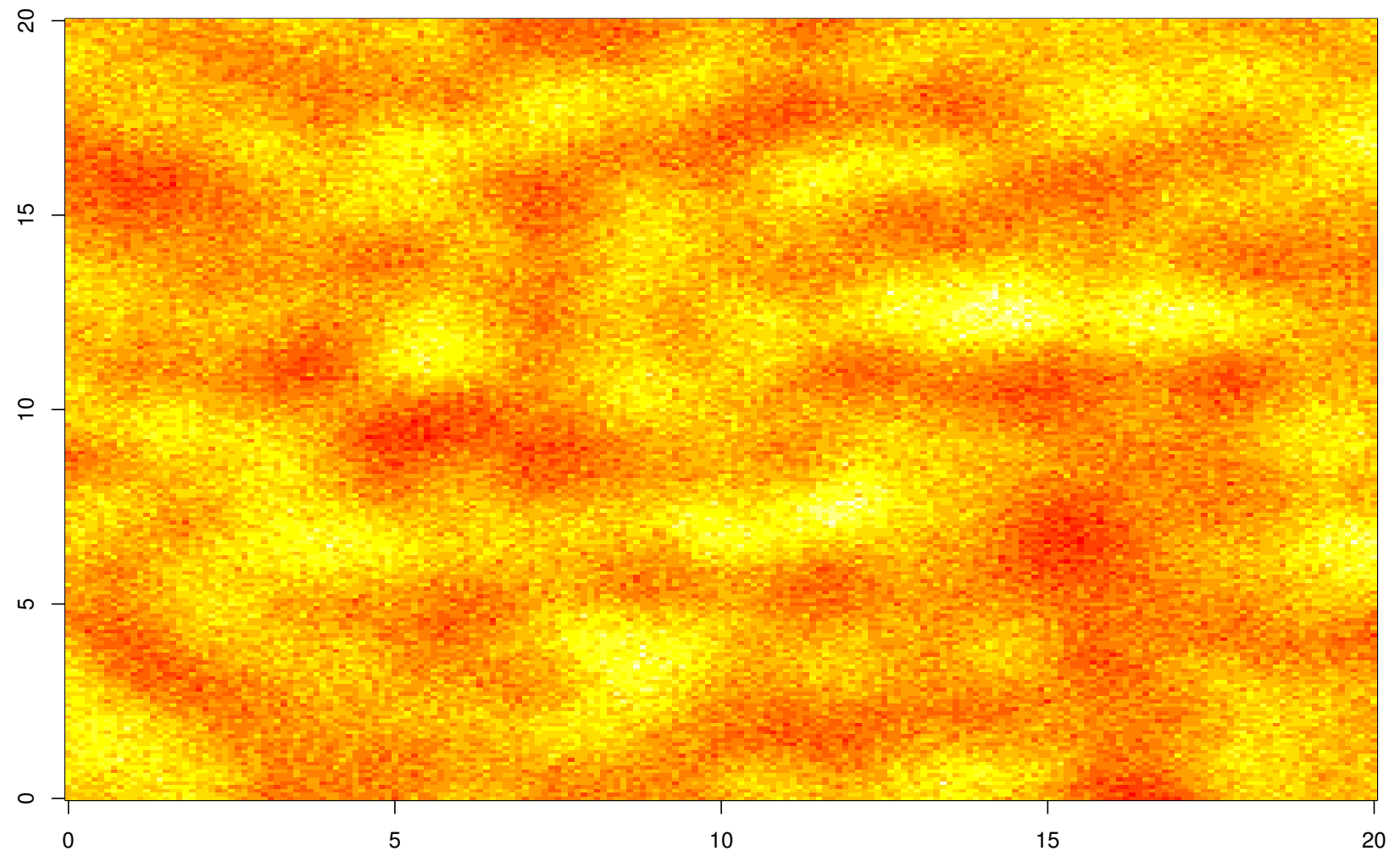

**Fig. 5.11.** Simulated Gaussian random field with penta covariance function

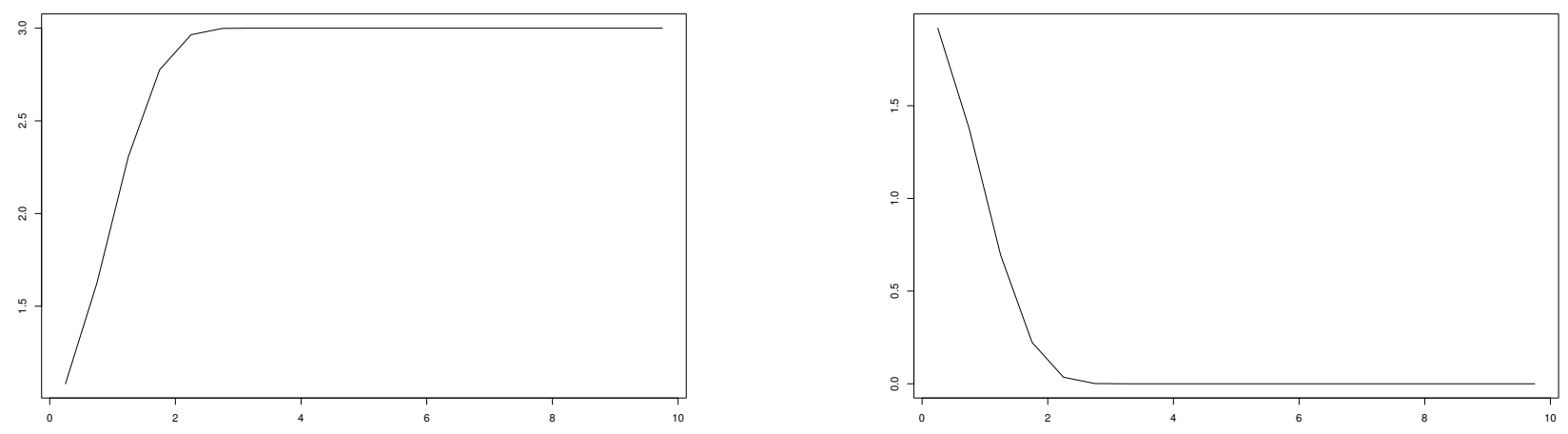

**Fig. 5.12.** Variogram and covariance function of a Gaussian random field with penta covariance function

### Power Variogram and Covariance Function

$$C(x) = \begin{cases} (1-x)^a & \text{for } 0 \le x \le 1 \\ 0 & \text{else} \end{cases}$$

This covariance function is valid for dimension $n$ if $a \ge (n+1)/2$. For $a = 1$ we get the well-known triangle (or tent) model, which is only valid on the real line.

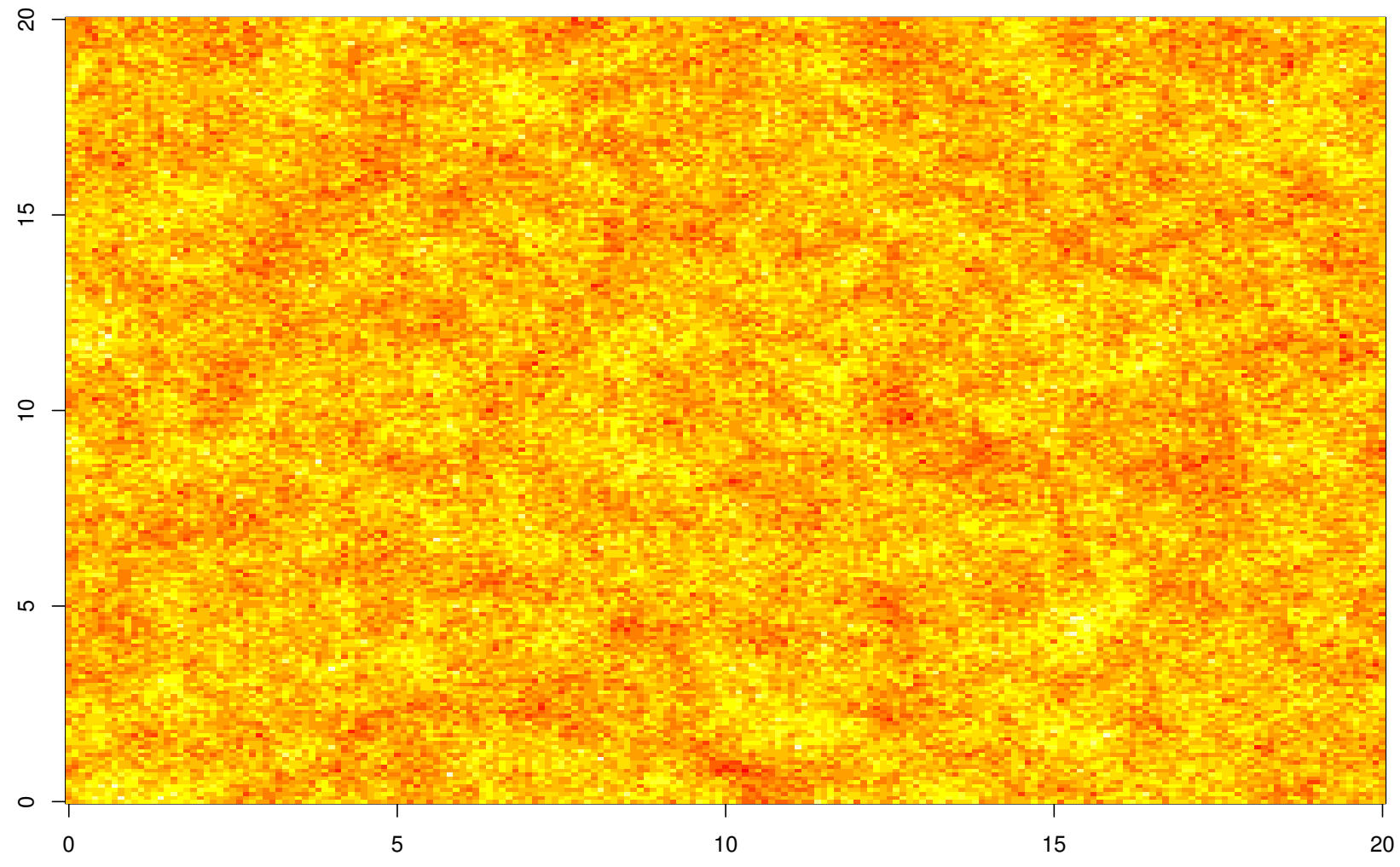

**Fig. 5.13.** Simulated Gaussian random field with power covariance function

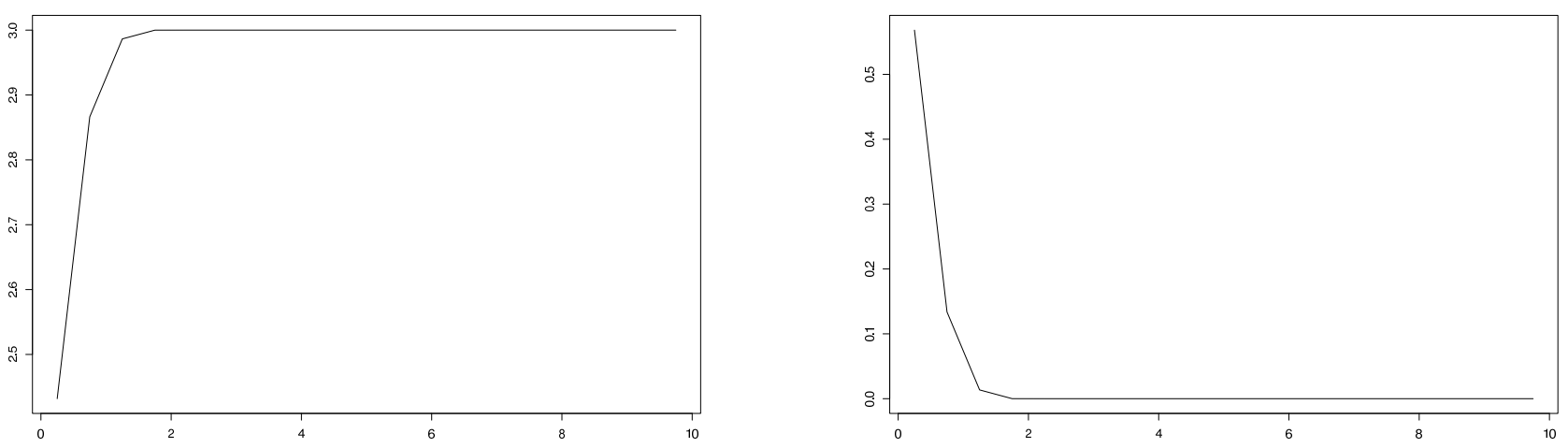

**Fig. 5.14.** Variogram and covariance function of a Gaussian random field with power covariance function

## Spherical Variogram and Covariance Function

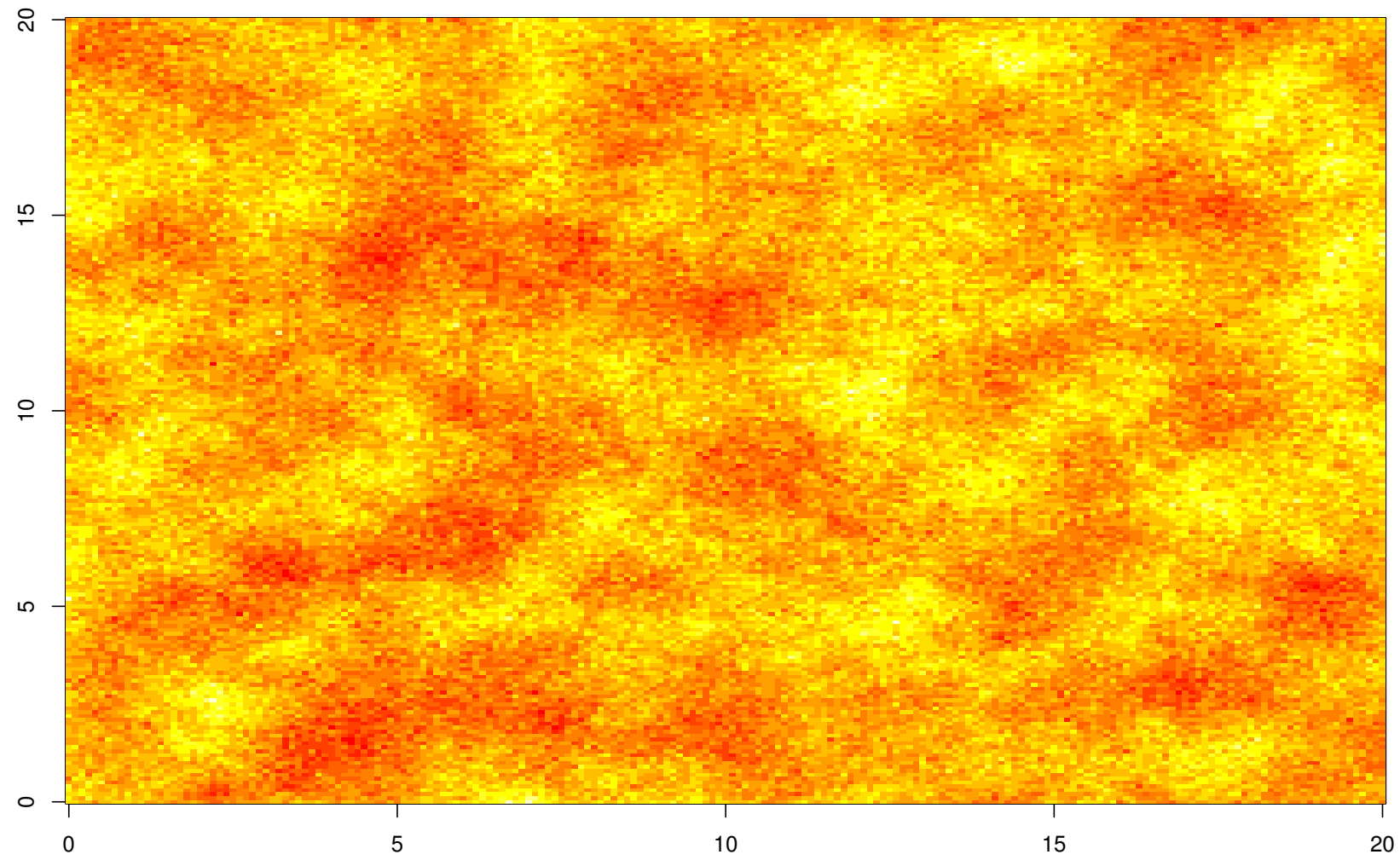

**Fig. 5.15.** Simulated Gaussian random field with spherical covariance function

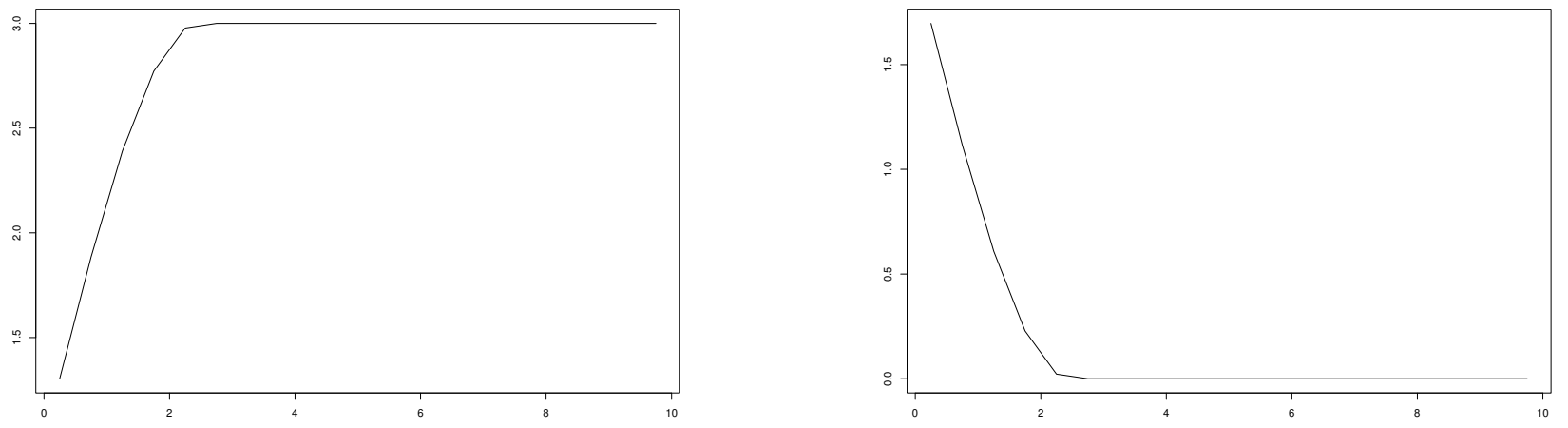

**Fig. 5.16.** Variogram and covariance function of a Gaussian random field with spherical covariance function

## Stable Variogram and Covariance Function

$$C(x) = e^{-x^a} \; \forall a \in [0, 2]$$

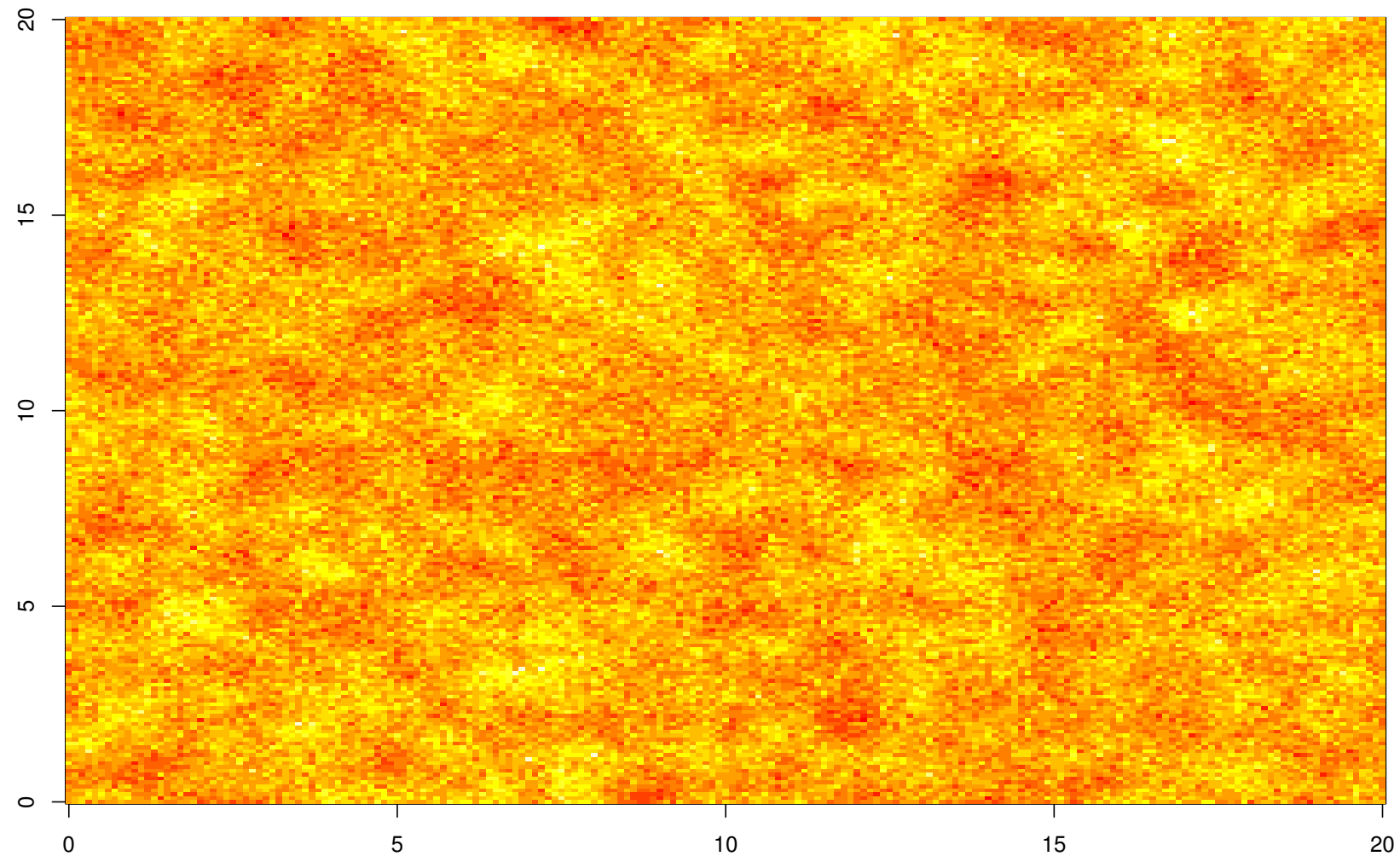

**Fig. 5.17.** Simulated Gaussian random field with stable covariance function

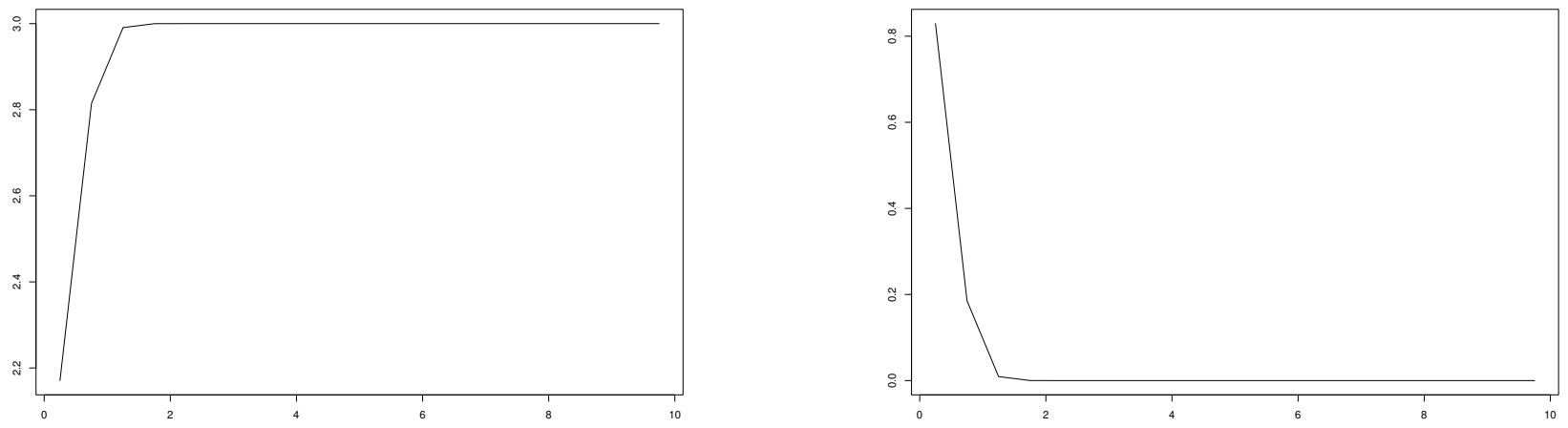

**Fig. 5.18.** Variogram and covariance function of a Gaussian random field with stable covariance function

### Wave Variogram and Covariance Function

$$C(x) = \frac{\sin(x)}{x} \quad \text{if} x > 0$$

In practise we can find many more covariance/variogram functions, like cardi-

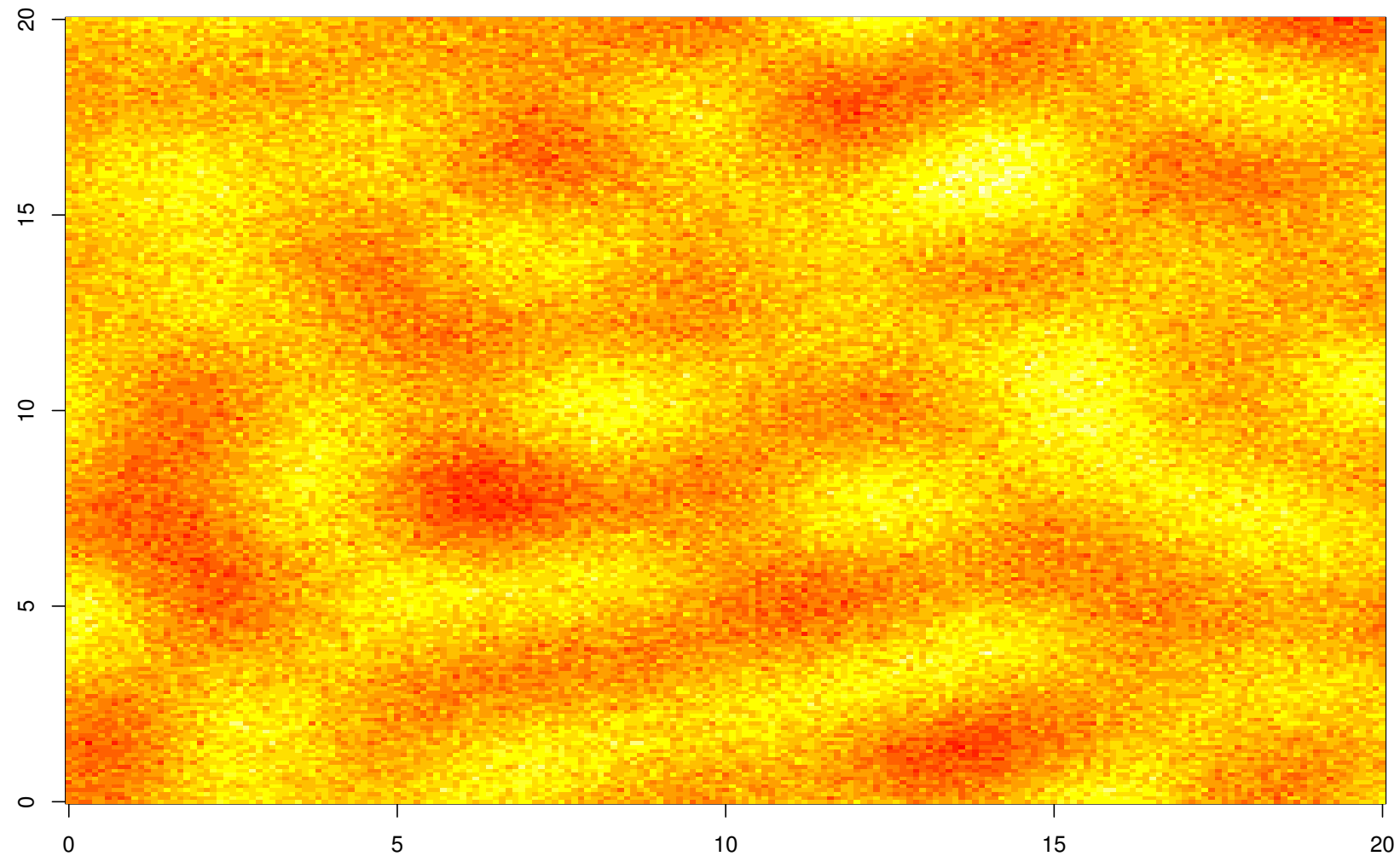

**Fig. 5.19.** Simulated Gaussian random field with wave covariance function

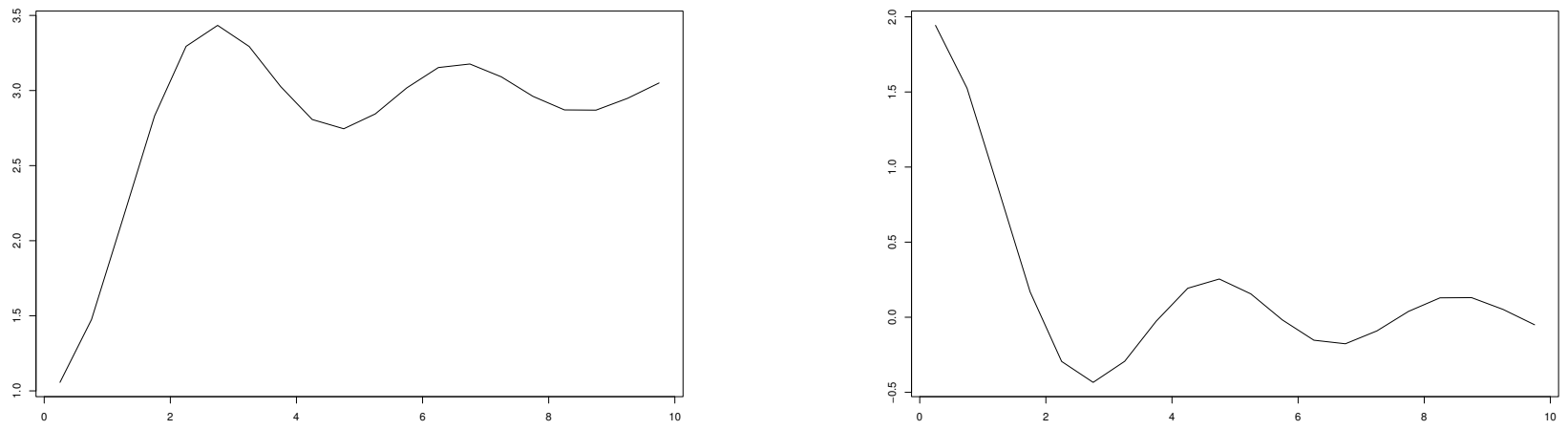

**Fig. 5.20.** Variogram and covariance function of a Gaussian random field with wave covariance function

nal sine, gamma, hyperbolic, J Bessel, K Bessel, linear, logarithmic, quadratic,

and so on and so on. For an overview we refer to Cressie(1993). For a covariance function, only the theoretical properties in foregoing chapters must be fulfilled, so everybody can construct his/her own covariance function.

## 5.8 Robust Variogram

Armstrong and Definer (1980) introduced a robust variogram estimator in the sense of Huber. Let $(a, b)$ be a pair such that $x_a - x_b = h$, if

$$\frac{(z(x_a) - z(x_b))^2}{2}$$

exceeds a threshold $c^2\gamma(h)$, then this threshold value is used for calculation of the variogram. This result is a downward bias of the variogram, a correction is introduced that reduces the bias in the case of a bi-Gaussian random field $Z$. The term $\frac{(z(x_a)-z(x_b))^2}{2\gamma(h)}$ is then the square of a standard Gaussian random field, and therefore

$$E[\min(\frac{(z(x_a) - z(x_b))^2}{2}, c^2\gamma(h))] = f(c)\gamma(h)$$

with $f(c)$ defined by $f(c) = E[\min(U^2, c^2)]$, where $U \sim N(0, 1)$. The robust variogram $\hat{\gamma}_h$ is defined as the value that achieves the equality

$$\text{Mean}_{x_a - x_b = h}[\min(\frac{(z(x_a) - z(x_b))^2}{2}, c^2\hat{\gamma}_H(h))] = f(c)\hat{\gamma}_H(h)$$

An implementation of that estimator is not yet done in R.
A very important covariance function is the Matern Covariance function, which is discussed in following section.

## 5.9 Matern Covariance Function

In this section we will introduce a general class of covariance functions that we feel is very flexible in modelling Gaussian random fields. This class is motivated by the smooth nature of its spectral density. This model includes a parameter that allows for any degree of differentiability for the random field and includes the exponential model as a special case and the Gaussian model as a limiting case. This class is name after Bertil Matern. The spectral density is given by

$$f(\omega) = \phi(\alpha^2 + \omega^2)^{-\nu - \frac{1}{2}}$$

for $\nu > 0, \phi > 0$ and $\alpha > 0$. The corresponding covariance function is

$$C(t) = \frac{\pi^{\frac{1}{2}}\phi}{2^{\nu-1}\Gamma(\nu+\frac{1}{2})\alpha^{2\nu}}(\alpha|t|)^{\nu}K_{\nu}(\alpha|t|)$$

where $K_\nu$ is a modified Bessel function. Here $\nu$ is the smoothness parameter of the random field $Z$, the larger $\nu$ is, the smoother $Z$ is. $Z$ will be $m$ times m.s. differentiable if and only if $\nu > m$, since

$$\int_{\mathbb{R}} \omega^{2m} f(\omega) d\omega < \infty$$

if and only if $\nu > m$. In the case that $\nu$ is of the form $m+\frac{1}{2}$ with $m$ a positive integer, the spectral density is rational and the covariance function is of the form $e^{-\alpha|t|}$ times a polynomial in $|t|$ of degree $m$. If we choose $\nu = \frac{1}{2}$ we get

$$C(t) = \pi\phi\alpha^{-1}e^{-\alpha|t|}$$

and for $\nu = \frac{3}{2}$ we get

$$C(t) = \frac{\pi}{2}\phi\alpha^{-3}e^{-\alpha|t|}(\alpha|t|+1).$$

Now we will discuss the m.s. smoothness of the process through the behaviour at the origin. Using the series representation of the Bessel and the gamma functions, we can rewrite the covariance function $C$ with a $\nu \in \mathbb{R}$ in the interval $(m, m+1)$ as

$$C(t) = \sum_{i=1}^{m} b_i t^{2i} - \frac{\pi\phi}{\Gamma(2\nu+1)\sin(\nu\pi)}|t|^{2\nu} + O(|t|^{2m+2}) \text{ as } t \to 0$$

for real constants $b_i$ depending on $\phi\nu$ and $\alpha$. For the choice of $\nu = m+1$, a positive integer, we get

$$C(t) = \sum_{i=0}^{m} b_i t^{2i} + \frac{2(-1)^m\phi}{(2m+2)!}t^{2m+2}\log|t| + O(t^{2m+2}) \text{ as } t \to 0$$

for constants $b_i$ depending on $\phi, m$ and $\alpha$. Stein(1999) mentioned the coefficients multiplying the principal irregular term does not depend on $\alpha$ and that

$$f(\omega) \sim \phi|\omega|^{-2\nu-1}$$

as $|\omega| \to 0$, so that high frequency behaviour of $f$ does not depend on $\alpha$. Under the use of following theorem (Stein(1999)):

**Theorem 5.1.** *If $C$ is positive definite on $\mathbb{R}$ and*

$$C(t) = \sum_{i=1}^{n} c_i t^{2i} + o(t^{2n})$$

*as $t \to 0$, then $C$ has $2n$ derivatives.*

We can conclude from the above theorem and the series expansion that

$$(\alpha|t|)^{\nu} K_{\nu}(\alpha|t|)$$

is $2m$ times differentiable if and only if $\nu > m$. Another interesting result in literature is the following one:

$$\mathrm{Var}(Z^m(h) - Z^m(0)) \sim \frac{2\phi}{\Gamma(2\nu - m + 1)\sin(\nu\pi)} h^{2(\nu-m)}$$

as $h \to 0$ and for $m < \nu < m + 1$. This result follows from term wise differentiation of above series expansion. For $\nu = m + 1$ we get

$$\mathrm{Var}(Z^m(h) - Z^m(0)) \sim 2\phi h^2 \log(h)$$

as $h \to 0$. So the continuous parameter $\nu$ has a direct interpretation as a measure of the smoothness of the process. One special thing in using the Matern class is, that it is very flexible, we are able to model nearly every random field with it. We have done simulation studies and simulated 50000 empirical variograms and fitted over 98 percent of them by the sum of theoretical Matern variograms.

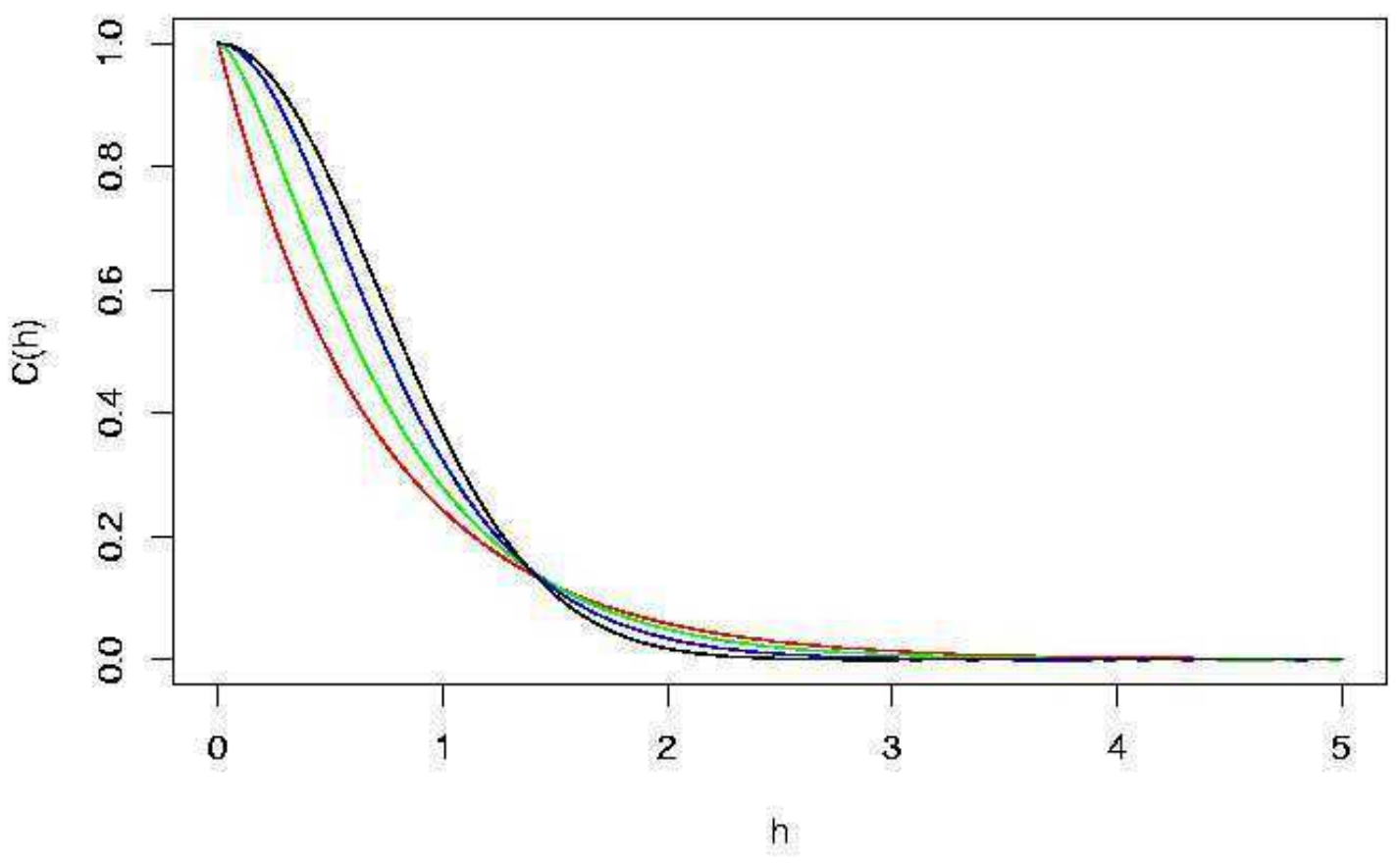

**Fig. 5.21.** Realisations of the Matern covariance function

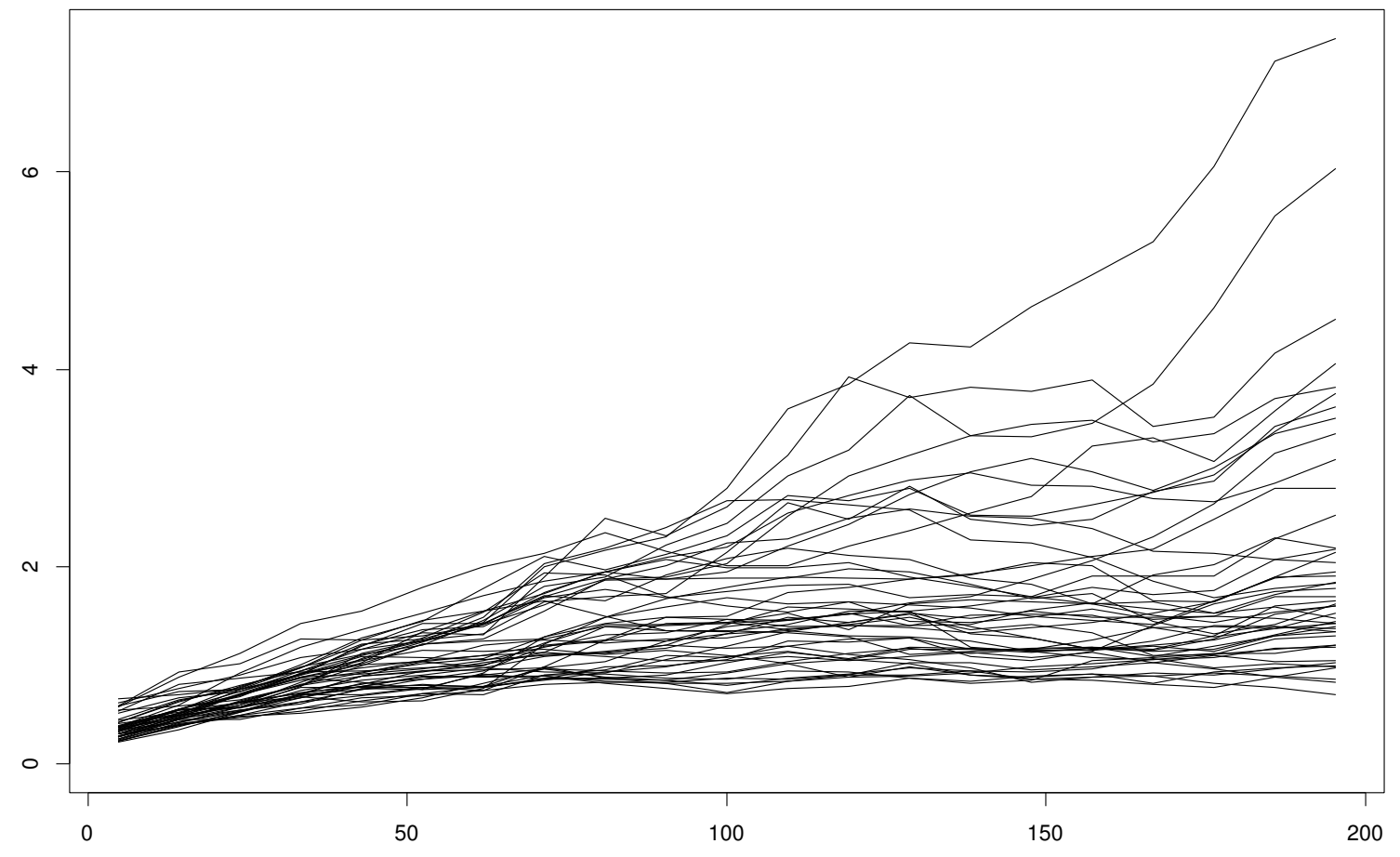

**Fig. 5.22.** Simulated empirical variograms

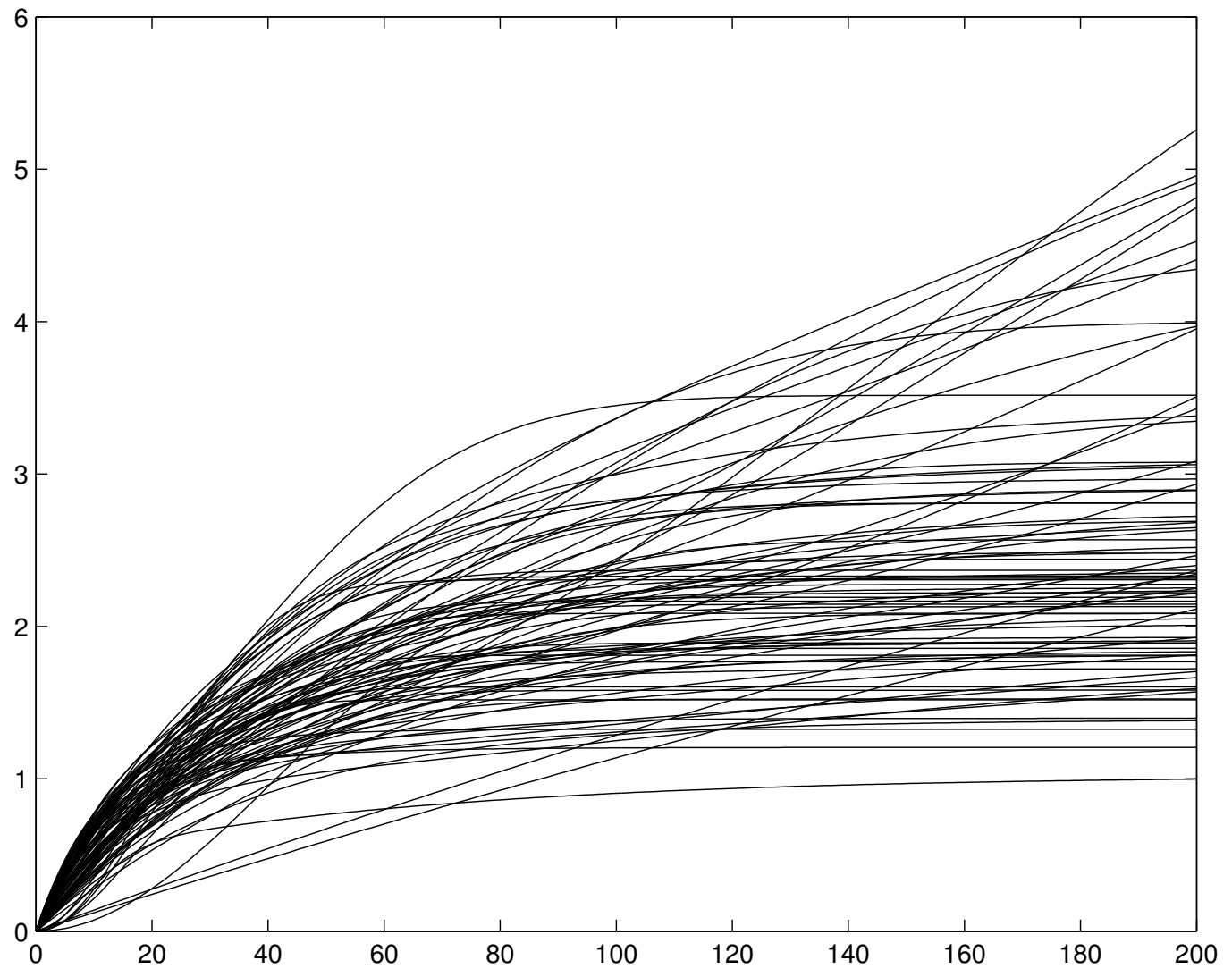

**Fig. 5.23.** Fitted theoretical variograms

# 6

# Link between Kriging and Splining

The main part of spatial statistics is engaged with models for prediction and smoothing of data. Since Krige(1951) introduced his prediction method for spatial data, a lot of other methods were developed, including stochastic methods for spatial prediction (for an overview see Cressie, (1990)), Bayesian non-parametric smoothing (Weerahandi and Zidek (1988)) and non-stochastic methods for spatial prediction, such as moving average, Median-polish plating (Cressie (1990)), Splines (Wahba (1978)) etc. Based on the discussion between Cressie and Wahba (kriging vs. splining) in The American Statistician, this chapter explores theoretical connections between these methods for spatial prediction. Horiwitz (1996) discussed in which parts of mathematics splines have a place. Here it will be shown how splines link with different branches of mathematics and stochastics, even how smoothing splines and thin plate splines can be seen as minimum variance estimators. Starting from special variation theory (Smirnow (1990)), radial functions, partial differences equations and Markov Random Fields the particular role of splines will be discussed.

## 6.1 Introduction

Spatial Statistics refers to a class of models and methods for spatial data which aims at providing quantitative descriptions of natural variables distributed in space or space and time. A family of techniques, stochastic and non stochastic ones were developed in geostatistics for that interpolation problem. The general approach is to consider a class of unbiased estimators, usually linear in the observations and to find the one with minimum uncertainty, as measured by the error variance. A group of techniques, known loosely as kriging, is a popular method among different interpolation techniques developed in geostatistics by Krige (1951), Matheron (1963) and Journel and Huijbregts (1978). An interesting comparison of ten classes of interpolation techniques with characteristics can be found in Burrough and McDonnell (1998) and in a lot of papers published recently a comparison of several interpolation techniques was made.
The goals of Kriging are like those of nonparametric regression, the understanding of spatial estimation is enriched by the interpretation as smoothing estimates. On the other hand random processes models are also valuable in setting uncertainty estimates for function estimates, specially in low noise situations (see Nychka (1998)). There are close connections between different mathematical subjects such as Kriging, radial basis functions (RBF) interpolations, spline interpolations, reproducing hilbert space kernels (rhsk), PDE, Markov Random Fields (MRF) etc. (see more in Matheron, 1980; Powell, 1992 and Royer 1981). A short discussion of these links is given in Horiwitz et al. and see 6.1. Splines link different fields of mathematics and statistics - and shortly they are used in statistics for spatial modelling (see more in Wahba(1990)).

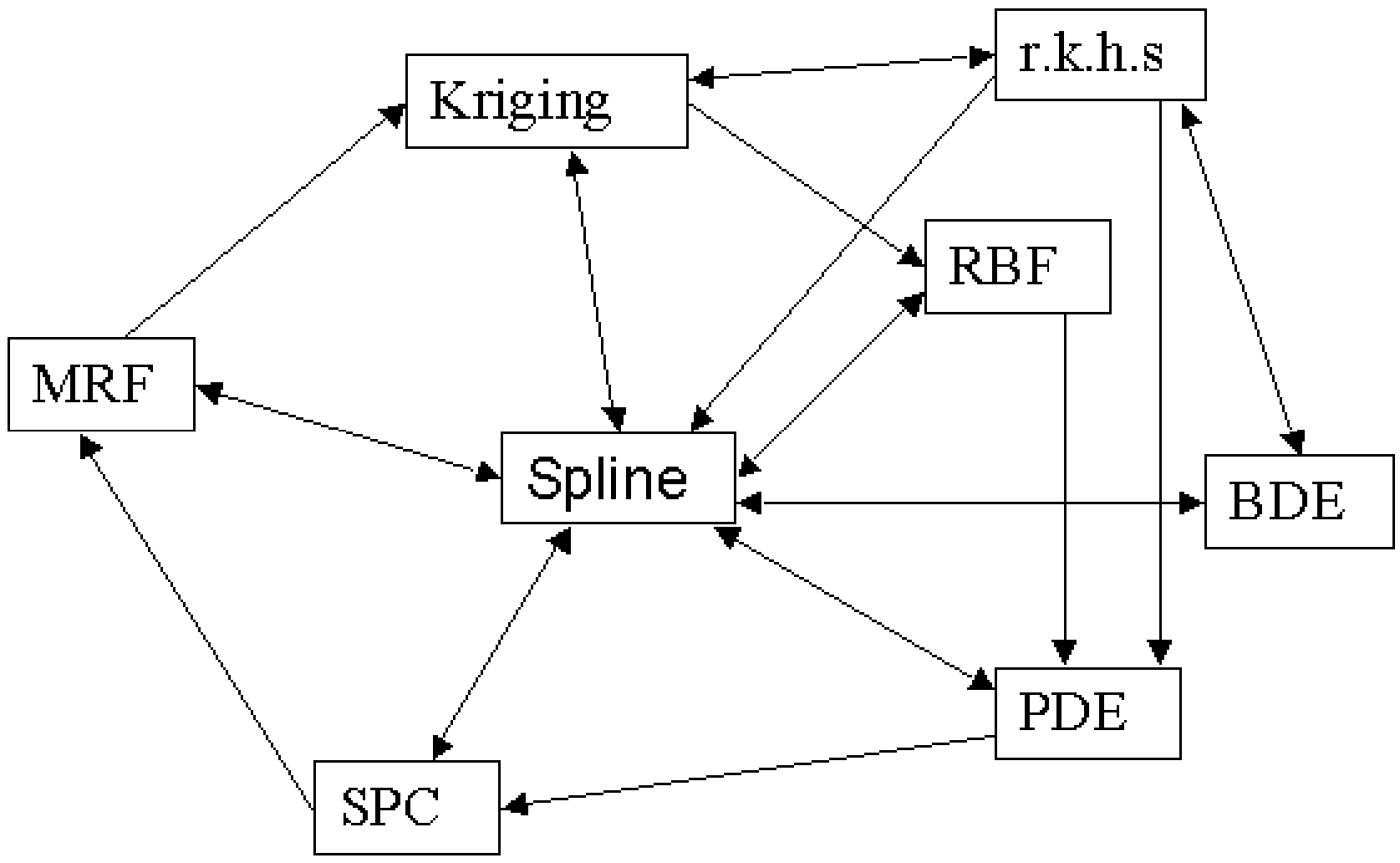

**Fig. 6.1.** Where splines link

## 6.2 Interpolation techniques

Interpolation techniques can be divided into techniques based on deterministic and stochastic models. The stochastic approach regards the data $\{y_i\}_{i=1}^n$ as a realisation of a random field on $R^d$ at $t_i = \{x_{i1}, ..., x_{id}\}$, $i = 1, ..., n$ and set $g(t)$ to be the best unbiased linear predictor of the random field at site $t$ given the measurements.
We assume here an intrinsic random field (second-order stationary). Splines are smooth real valued functions $g(t)$. Define a roughness penalty based on the sum of integrated squared partial derivatives of a given order $n$. The choice of $g(t)$ which interpolates the data and minimises the roughness penalty is known as the smoothing thin plate spline introduced by Reinsch (1967).

### 6.2.1 Univariate Spline Function

A theoretical definition and historical motivation can be found in Haemmerlin, Hoffmann(1992). The natural spline $S(x) = S_n$ is a real valued function $S : [a, b] \to R$ with $n$ knots $-\infty \leq a < x_1 < x_2 < ... < x_n \leq \infty$with following properties

$$\begin{aligned}
&S \in \Pi^{m-1} \text{ for } x \in [a, x_1] \text{ and } x \in [x_n, b] \\
&S \in \Pi^{2m-1} \text{ for } x \in [x_i, x_{i+1}], i = 1, ..., n-1 \\
&S \in C^{2m-2} \text{ for } x \in (-\infty, \infty) \\
&f(x_i) = f_i \text{ for } i = 1, ..., n
\end{aligned}$$

$\Pi^q$ is the class of polynomials with degree $q$ and $C^p$ the class of continuous functions of order $p$. The historical problem (Schoenberg, 1964) is to find a function $f$ in a function space with continuos derivatives of order $(m-1)$ with the minimum of $\int_a^b (f^{(m)}(x))^2 dx$ under all functions with the property $f(x_i) = f_i$ for $i = 1, ..., n$. Schoenberg shows that the solution is our natural spline.

The statistical approach turns attention to smooth the data and not to interpolate it. The first access to this problem was made by Reinsch(1967) with finding a function $g$ that minimises

$$\int_{x_0}^{x_n} g''(x)^2 dx + \rho\{\sum_{i=0}^{n} (\frac{g(x_i) - y_i}{\delta y_i})^2 + z^2 - \mathcal{S}\} \tag{6.1}$$

where $z$ and $\mathcal{S}$ are auxiliary variables and $\rho$ is a Lagrange parameter. The solution for (6.1) is the cubic spline:

$$g(x) = a_i + b_i(x - x_i) + c_i(x - x_i)^2 + d_i(x - x_i)^3 \text{ for } x_i \leq x \leq x_{i+1} \tag{6.2}$$

### 6.2.2 Multivariate Spline

In analogy with the interpolation problem in one dimension we now handle a data recording three dimensions with a grid of points $(x_i, y_i)$, $i = 1, ..., n$. For polynomial interpolation we can use a polynomial with low degree $r$ of the form

$$P_r(x, y) = \sum_{p+q=0}^{r} a_{pq} x^p y^q \tag{6.3}$$

On an arbitrary grid it is in general not possible to construct a unique solution. For our views we assume a rectangular region on which we build a rectangular $(n+1)(k+1)$ grid of the form

$$\begin{aligned} a = x_0 < x_1 < ... < x_n = b \\ c = y_0 < y_1 < ... < y_k = d \end{aligned}$$

where in $x$ direction and in $y$ direction a spline will be constructed. As an analogue to the univariate B-Spline we can give the following spline base (see Haemmerlin, Hoffmann (1992))

$$B_{1\nu\kappa} = \begin{cases} \frac{(x_{\nu+2}-x)(y_{\kappa+2}-y)}{(x_{\nu+2}-x_{\nu+1})(y_{\kappa+2}-y_{\kappa+1})} & \text{for } (x,y) \in I_1 \\ \frac{(x-x_\nu)(y_{\kappa+2}-y)}{(x_{\nu+1}-x_\nu)(y_{\kappa+2}-y_{\kappa+1})} & \text{for } (x,y) \in I_2 \\ \frac{(x-x_\nu)(y-y_\kappa)}{(x_{\nu+1}-x_\nu)(y_{\kappa+1}-y_\kappa)} & \text{for } (x,y) \in I_3 \\ \frac{(x_{\nu+2}-x)(y-y_\kappa)}{(x_{\nu+2}-x_{\nu+1})(y_{\kappa+1}-y_\kappa)} & \text{for } (x,y) \in I_4 \end{cases} \tag{6.4}$$

### 6.2.3 Additive Model

Now we will handle a $(d+1)$-dimensional data record $(x_i, Y_i)$, with $x_i = (x_{i1}, x_{i2}, ..., x_{id})^T$; $x_1, x_2, ..., x_n$ are independent realisations of the random vector $X = (X_1, X_2, ..., X_d)$, where $x_i$ and $Y_i$ fulfil

$$Y_i = g(x_i) + \varepsilon_i \tag{6.5}$$

with $1 \le i \le n$, where $g$ is an unknown smoothing function mapping from $R^d$ to $R$, and $\varepsilon_1, \varepsilon_2, ..., \varepsilon_n$ are independent error terms. Our aim is to find a valid additive approximation of $g$ by

$$g(X) \approx g_0 + \sum_{i=1}^{d} g_i(X_i) \tag{6.6}$$

The additive model is defined by

$$Y = \alpha + \sum_{i=1}^{d} g_i(X_i) + \varepsilon \tag{6.7}$$

Like in the multilinear regression case the error term is independent of $X_i$, $E(\varepsilon) = 0$ and $var(\varepsilon) = \sigma^2$. From (6.7) we can conclude that for any predictor and for all dimensions there exists at least one function $g$. More can be found in Hastie (1990).

## 6.3 Smoothing Spline

In this paragraph we consider the following problem: We want to find a function $g(x)$ under all two times continuously differentiable functions that minimise (6.1).

### 6.3.1 Univariate Approach

First we have to define:

$$h_i = x_{i+1} - x_i \text{ for } i = 1, ..., (n-1)$$

$$\Delta = \begin{cases} \Delta_{i,i} = \frac{1}{h_i} \\ \Delta_{i,i+1} = -(\frac{1}{h_i} + \frac{1}{h_{i+1}}) \quad \text{a } (n-2) \times n \text{ matrix} \\ \Delta_{i,i+2} = \frac{1}{h_{i+1}} \end{cases}$$

$$C = \begin{cases} C_{i-1,i} = \frac{h_i}{6} \\ C_{i,i} = \frac{1}{6}(h_i + h_{i+1}) \quad \text{a } (n-2) \times (n-2) \text{ matrix} \\ C_{i,i-1} = \frac{h_i}{6} \end{cases}$$

With this representation the minimisation problem in (6.1) is equivalent to

$$\| y - g \| + \rho g^T K g \tag{6.8}$$

where K is a quadratic penalty matrix,

$$K = \Delta^T C^{-1} \Delta.$$

Calculation of the inverse of $C$ is possible because it is strictly diagonal dominant (see Stoer (1983)). The solution to (6.8) is given by

$$\hat{g} = Sy \tag{6.9}$$

where $S$ is the smoothing Matrix of the form

$$S = (I + \rho K)^{-1}. \tag{6.10}$$

Equation (6.10) we will later find in the solution of the kriging prediction problem and equation (6.10)can be also found in ridge regression with natural basis and Demmler-Reinsch basis (see more in Nychka(2000))

### 6.3.2 Multivariate Approach

When we use this approach with the additive model we are able to find (with help of equation (6.8)) the following result

$$g(X) \approx \sum_{i=1}^{d} \rho_i g_i^T K_i g_i + \sum_{i=1}^{d} (Y_i - g_0 - \sum_{j=1}^{n} g_j(x_{ij}))^2 \tag{6.11}$$

Minimisation leads to

$$\hat{g}_l = S_l(y - g_0 - \sum_{j=1 j \neq l}^{d} \hat{g}_j) \tag{6.12}$$

$$\text{with } S_l = (I + \rho_l K_l)^{-1} \text{ for } l = 1, ..., d \tag{6.13}$$

$S_l$ can be seen as smoothing matrix and $K_l$ as penalty matrix. Then the following matrixequation system must be solved:

$$\begin{pmatrix} I & S_1 & S_1 & \dots & S_1 \\ S_2 & I & S_2 & \dots & S_2 \\ \vdots & \vdots & \vdots & \ddots & \vdots \\ S_d & S_d & S_d & \dots & I \end{pmatrix} \begin{pmatrix} g_1 \\ g_2 \\ \vdots \\ g_d \end{pmatrix} = \begin{pmatrix} S_1 y \\ S_2 y \\ \vdots \\ S_d y \end{pmatrix} \tag{6.14}$$

shortly we can write $Pg = Qy$ where $P$ is a $(nd) \times (nd)$ and $Q$ is a $(nd) \times (nd)$block matrix. A solution of that problem can be found by Backfitting Algorithm (see Schimek and Stettner(1993) or Green and Yandell(1985)). A solution to equation (6.14) is given by

$$\hat{g}(x) = \hat{Y} = g_0 + \sum_{j=1}^{d} O_j^{-1} R_j^T g_j \tag{6.15}$$

where $R_j$ is a reduction matrix and $O$ an ordermatrix (the solution can be found in Stingel (1994)).

### 6.3.3 Abstract Minimisation

In Wahba(1990) a more general approach to the problem of spline smoothing is given. She discusses the problem of finding a function $f$ in the Sobolev function space of the form

$$W_m : W_m[0,1] = \{f, f', f'', ..., f^{(m-1)} \text{ absolutely continuous, } f^{(m)} \in L^2\}$$

where the solution can be found by abstract analysis and abstract optimisation. Wahba gives the following equivalent smoothing problem via r.h.k.s.

$$\frac{1}{n} \sum_{i=1}^{n} (y_i - \langle \eta_i, f \rangle)^2 + \rho \parallel P_1 f \parallel_R^2 \rightarrow \min_{f \in W_m} \tag{6.16}$$

with the solution given by

$$\begin{aligned} f_\rho &= \sum_{\nu=1}^{M} d_\nu \phi_\nu + \sum_{i=1}^{n} c_i \xi_i \\ \xi_i &= P_1 \eta_i, i = 1, ..., n \\ d &= (d_1, d_2, ...d_M)^T = ((T^T M^{-1} T^T)^{-1} T^T M^{-1}) y \\ c &= (c_1, ..., c_n)^T = M^{-1}(I - T(T^T M^{-1} T)^{-1} T^T M^{-1}) y \\ M &= \Sigma + n\rho I \\ \Sigma &= (\langle \xi_i, \xi_j \rangle), i = 1, ...n \text{ and } j = 1, ...M \end{aligned}$$

$P_1$ is a projection matrix. A more general way to describe the roughness penalty function is in the form

$$J_{r+1}^{d}(g) = \sum_{|m|=r+1} \frac{(r+1)!}{m!(r+1-m)!} \int_{R^d} \left( \frac{\partial^{r+1} g(t)}{\partial t_1^{m_1} \cdots \partial t_d^{m_d}} \right)^2 dt \tag{6.17}$$

When a particular penalty is chosen then the result is invariant under rotations and translations of $t$.

## 6.4 Kriging

Let $\{Y(t), t \in R^d\}$ be an intrinsic or stationary random field, we will additionally demand a polynomial drift of order $n > 0$. In that case the drift is a linear combinationof $t^m$ for $\mid m \mid < n$ with unknown coefficients. The subspace of the polynomial drift has an order less than $n$ with dimension

$$\nu = \binom{d+n}{d}$$

Let

$$V_r = \begin{pmatrix} f(t_0)^T \\ F \end{pmatrix}$$
$$F = (f(t_1), ..., f(t_n))^T$$
$$EY(t) = \beta^T f(t),$$

$V_r$ a $(n+1) \times \nu$ matrix of the drift at locations $t_0, ..., t_n$, $F$ a $(n \times \nu)$ matrix that describes the drift at locations $t_1, ..., t_n$ and $u_{r,0}$ is a vector of lentgh $\nu$ with the elements $t_0^m$, $\mid m \mid < n$. The covariance function $\sigma_{i,j}$ is defined as $\sigma(t_i - t_j)$ for $i, j = 1, ..., n$ and $\varphi_{i,j} = \sigma(t_i - t_j)$ for $i, j = 0, ..., n$.

$$\Phi = \begin{pmatrix} \sigma^2 & \sigma_0^T \\ \sigma_0 & \Sigma \end{pmatrix}$$

where $\sigma^2 = \sigma(0)$ and $\sigma_0$ is a $(n \times 1)$ vector with elements $\sigma(t_0 - t_i), i = 1, ..., n$. $\Sigma$ has the elements $\sigma_{i,j}$.

If $\{Y(t)\}$ is a stationary random field with a polynomial drift, then kriging involves predicting $Y(t_0)$ by a linear combination $\hat{Y}(t_0) = \alpha^T y$. The goal is to minimise the prediction mean squared error subject to an unbiasedness constraint.

$$var(Y(t_0) - \hat{Y}(t_0)) = var(Y(t_0) - \alpha^T y) = var(\beta^T z) = (c^T \Phi c) \tag{6.18}$$

under the unbiasedness constraint $Ec^T z = c^T V_r \beta = 0$. It's straightforward to minimise (6.18) by Lagrange multipliers to give

$$\alpha = AF_{r,0} + B\sigma_0 \tag{6.19}$$

where

$$A = \Sigma^{-1}F(F^T\Sigma^{-1}F)^{-1} \tag{6.20}$$

$$B = \Sigma^{-1} - \Sigma^{-1}F(F^T\Sigma^{-1}F)^{-1}F^T\Sigma^{-1} \tag{6.21}$$

for more details see Kent and Mardia(1994) who also give a solution for intrinsic random fields where $\sigma(h)$ is a polynomial in $h$ with degree $2p$. $A$ and $B$ can be found with the use of Moore-Penrose generalised inverse.

$$B = [(I - U_r(U_r^T U_r)^{-1}U_r^T)\Sigma(I - U_r(U_r^T U_r)^{-1}U_r^T)]^- \tag{6.22}$$

$$A = (I - B\Sigma)U_r(U_r^T U_R)^{-1} \tag{6.23}$$

## 6.5 Link between Kriging and Thin Plate Splining

The following theorem, which identifies kriging solution and thin-plate spline is one of the main results in literature (see Kent and Mardia (1994)) for the discussion 'kriging vs. splines'.

**Theorem 6.1.** *Let $r + 1 > \frac{1}{2}d$ (d is the dimension, r order of the polynomial drift)and set $\alpha = r + 1 - \frac{d}{2} > 0$. Then the problem of interpolating data $(y_i, t_i)$, $i = 1, ..., n$ subject to minimising the roughness penalty (6.17) has a solution $g^*(t_0)$ given by*

$$g^*(t_0) = y^T Au_{r,0} + y^T B\sigma_0 \tag{6.24}$$

where $A, B, u_{r,0}$ and $\sigma_0$ are determined as before. The link between kriging and thin-plate splines also holds for the smoothing problem as well as the interpolation problem shown by the last theorem. In the thin-plate spline approach a smoothing function $g(t)$ with square-integrable $(r + 1)^{th}$order derivates has to be found such that

$$F(g, \rho) = \sum \mid y_i - g(t_i) \mid^2 + \rho J^d_{r+1}(g) \tag{6.25}$$

is minimised. The solution to this problem can be found in Kent and Mardia (1994) or in an equivalent formulation in Nychka (2000).
The optimal choice of $g$ is given by

$$g(t_0) = y^T(I + \kappa B)^{-1}Au_0 + y^T(I + \kappa B)^{-1}B\sigma_0 \tag{6.26}$$

which is the same as the kriging predictor with a nugget effect.

## 6.6 Example for Kriging and Splining

Now we want to compare these two methods with the help of a modification of Wendelberger's test function, this function will be disturbed by some noise ($\sim N(0, 1)$). This data set $x, y$ and $z$ will be smoothed by thin plate splines and with the help of GCV introduced in Wahba (1990).

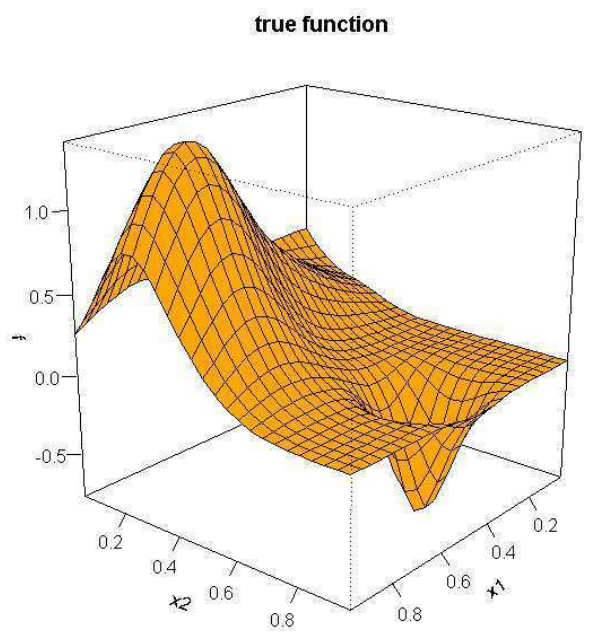

**Fig. 6.2.** Modification of Wendelberger's test function

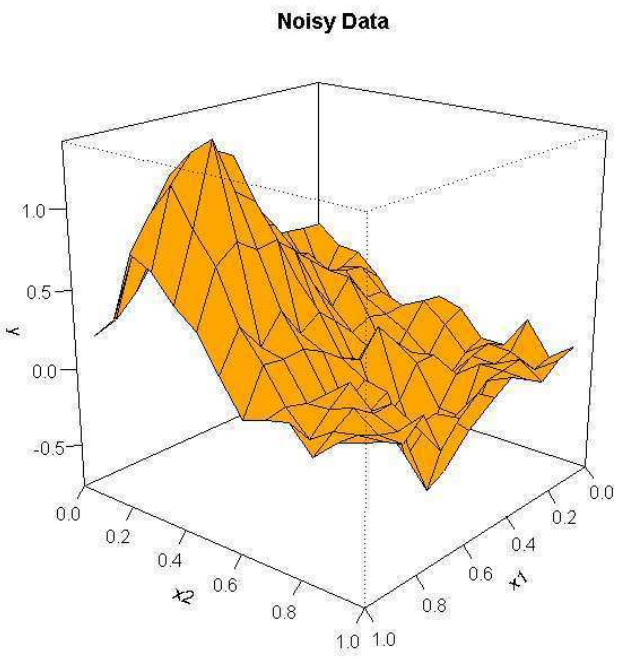

**Fig. 6.3.** Wendelberger's test function with noise

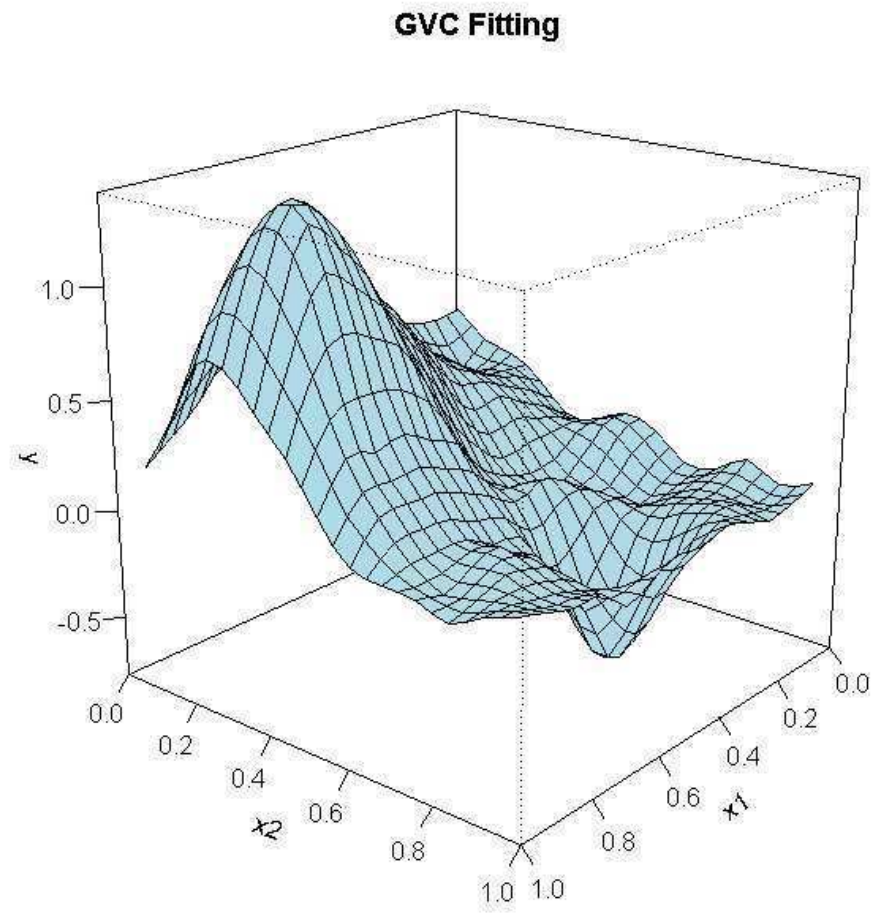

**Fig. 6.4.** GCV fitted Thin Plate Spline

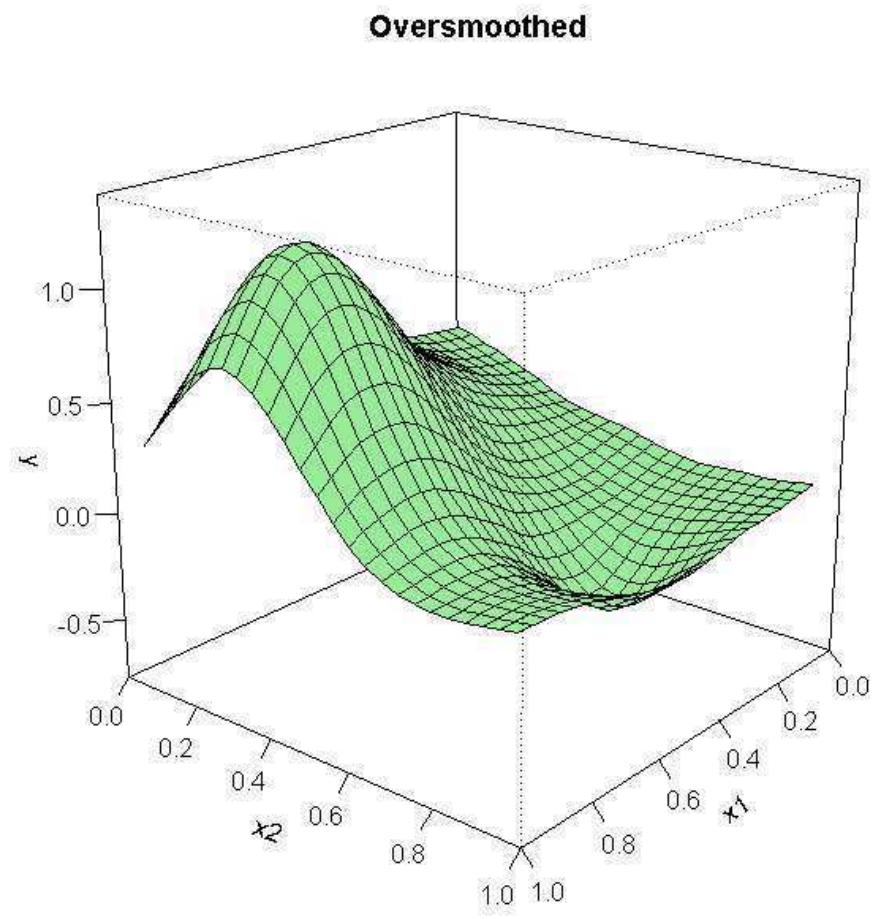

**Fig. 6.5.** Oversmoothed Thin Plate Spline

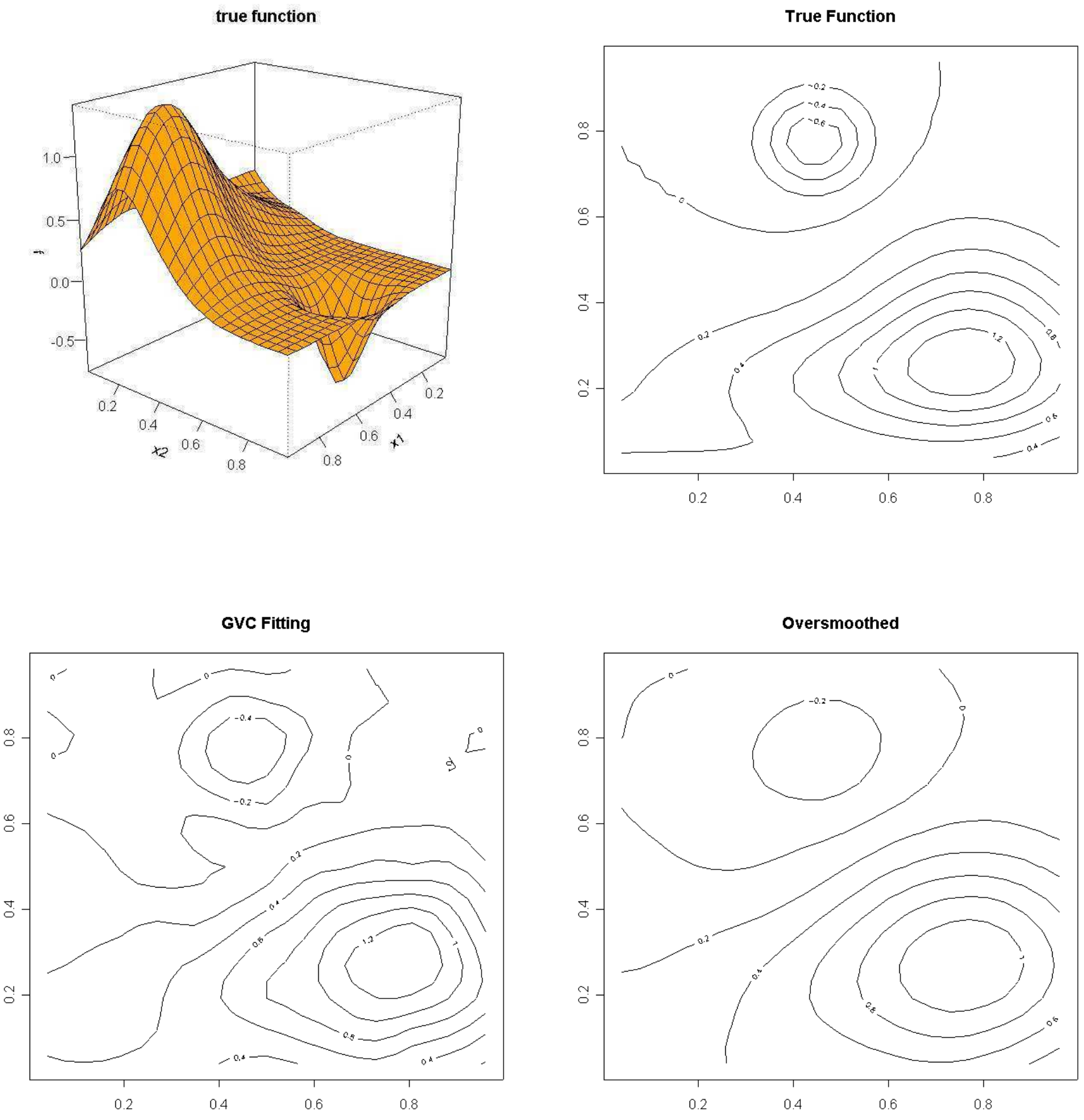

**Fig. 6.6.** Contourplots

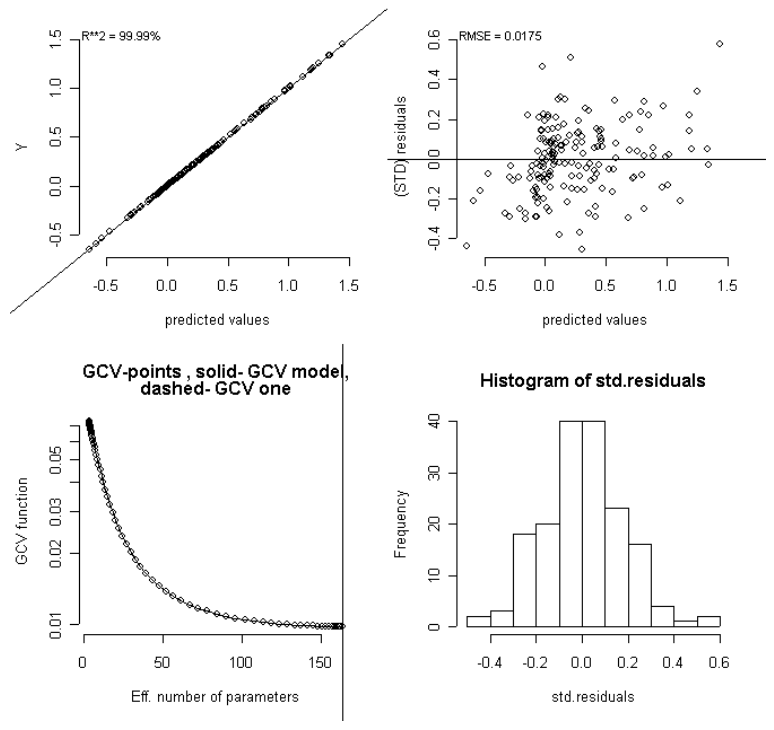

**Fig. 6.7.** Summary of Kriging with Matern Covariance Function

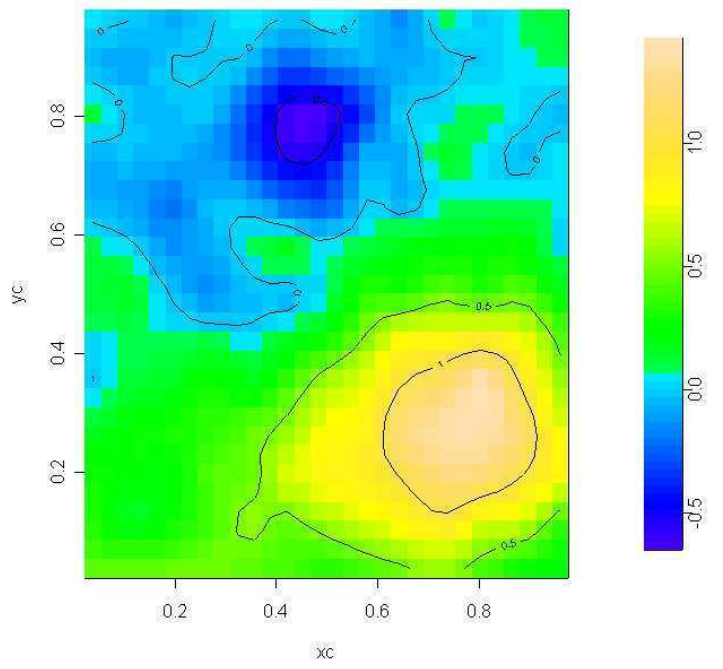

**Fig. 6.8.** Contourplot

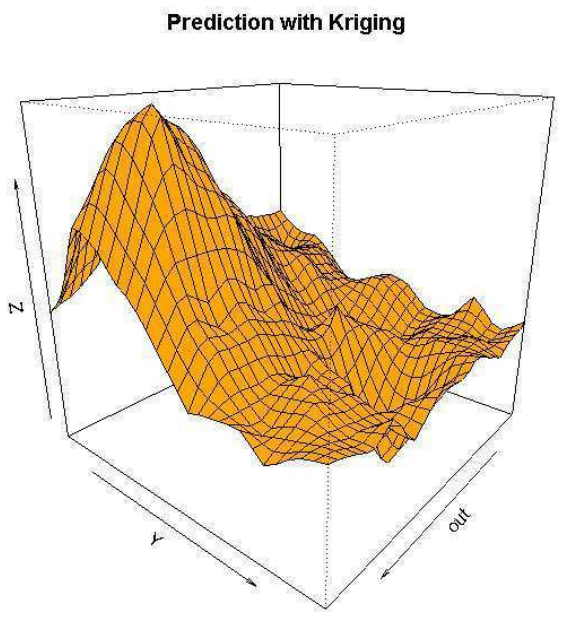

**Fig. 6.9.** Prediction of the Surface

# 7

# Kriging and Partial Stochastic Differential Equations

In physics laws are often expressed in terms of partial differential equations, these laws are completely determined by the observable process only, when we have sufficiently many boundary conditions and we can use numerical methods to solve these equations. Kriging is a method for interpolating random processes based on the estimated expectation function and covariance function. When applying Kriging, we do not care about the underlying physical law which is governing the process. It can be shown that linear partial differential equations impose restrictions on the class of admissible covariance and expectation functions, and the known physical law can help to select a physical reasonable covariance function. Under these admissible classes Kriging estimations solve differential equations in the mean square sense.

## 7.1 Stochastic Differential Equations

A differential system is a set of equations relating to a variable $t$, often interpreted as time, a number of functions $f(t)$ and their derivatives. Normally this system is given in the form

$$\begin{cases} \frac{dZ(t)}{dt} = F(Z(t), Y(t)) \\ \quad Z(t_0) = Z_0 \end{cases}$$

where $Z(t)$ is the vector of unknown functions with $\dim Z(t) = m$, and the known functions are given in the vector $Y(t)$ with $\dim Y(T) = p$. $F$ is a vector transformation operator. The usual objective is to find the expression of $Z(t)$ as a function of $Y(t)$ and initial conditions. This problem becomes a stochastic one if the known functions, the equation parameters or the initial conditions are made random.

## 7.2 Differential Equations in Mean Square Sense

To consider differential equations of a random field we need the concept of differentiation in mean square sense, a random function $f : \mathbb{R}^d \to \mathbb{R}$ or some other Frechet space is called mean square differentiable in $x_0$ in direction $x$ if the differential quotient

$$\lim_{h\to 0} \frac{f(x_0 + hx) - f(x_0)}{h}$$

exists in the first two moments. That is

$$| \lim_{h\to 0} \mathrm{E}(\frac{f(x_0 + hx) - f(x_0)}{h})| < \infty$$
$$| \lim_{h\to 0} \mathrm{Var}(\frac{f(x_0 + hx) - f(x_0)}{h})| < \infty$$

are finite. Higher order derivatives, partial derivatives and gradients are defined in a common way. If the gradient

$$\nabla_x f(x) = (\frac{\partial}{\partial x_1} f(x), ..., \frac{\partial}{\partial x_d} f(x))^T$$

exists and has finite variation it can be shown that

$$\mathrm{E}(\nabla_x f(x)) = \nabla_x \mathrm{E}(f(x))$$
$$\mathrm{Cov}(\nabla_x f(x), f(y)) = \nabla_x C(x, y)$$
$$\mathrm{Cov}(\nabla_x f(x), \nabla_y f(y)) = \nabla_x \nabla_y C(x, y)$$

where $C$ is the covariance function, details can be found in Chiles, Delfiner (1999). Note that a function $f$ is $n$ times differentiable in m.s. sense, if and only if the derivative

$$\frac{d^n}{dx^n} \frac{d^n}{dy^n} C(x, y)|_{x=y} < \infty$$

exists, details can be found in Christakos (1992).

## 7.3 Partial linear differential Equations

We are going to consider linear partial differential equations in the form

$$L_x f(x) = \sum_{i=0}^{n} \alpha_{ij}^{j_1 ... j_i}(x) \frac{\partial^i}{\partial x_{j_1} ... \partial x_{j_i}} f(x) = k_j(x)$$

with tensor valued functions $\alpha_i : \mathbb{R}^d \to \mathbb{R}^{p\times d\times ...\times d}$, where $d\times ...\times d$ is $i$ times, $k : \mathbb{R}^d \to \mathbb{R}^p$ and $f : \mathbb{R}^d \to \mathbb{R}$ and $x \in \mathbb{R}^d$ is a given location.

**Definition 7.1 (Homogeneous linear partial differential Equation).** *We call a linear partial differential equation homogeneous if it is in the form*

$$\sum_{i=0}^{n} \alpha_{ij}^{j_1 \ldots j_i}(x) \frac{\partial^i}{\partial x_{j_1} \ldots \partial x_{j_i}} f(x) = 0$$

**Definition 7.2 (Stationary linear partial differential Equation).** *We call a linear partial differential equation stationary if $\alpha_i$ and $k$ do not depend on the location $x$.*

**Definition 7.3 (Isotrope linear partial differential Equation).** *We call a linear partial differential equation isotrope if the equation is invariant under orthogonal transformations of the coordinate system.*

**Definition 7.4 (Even linear partial differential Equation).** *We call a linear partial differential equation even, if*

$$\sum_{i=0}^{n} \alpha_{ij}^{j_1 \ldots j_i}(x) \frac{\partial^i}{\partial x_{j_1} \ldots \partial x_{j_i}} f(x) = \sum_{i=0}^{n} \alpha_{ij}^{j_1 \ldots j_i}(-x) \frac{\partial^i}{\partial x_{j_1} \ldots \partial x_{j_i}} f(-x)$$

*This is equivalent to $\alpha_i = 0$ for odd $i$.*

As an example the Laplace equation given in the form

$$\Delta_x f(x) = \left(\frac{d^2}{dx_1^2} + \frac{d^2}{dx_2^2}\right) = 0$$

is homogeneous, stationary, even and isotrope.

**Theorem 7.5 (Moments of linear transformed Processes).** *For two random functions $f(x)$ and $g(y)$ for which*

$$E(f(x)) < \infty$$
$$E(g(x)) < \infty$$
$$Cov(f(x), g(y)) < \infty$$

*exist and the differential operator $L_x$ exists such that*

$$L_x E(f(x)) < \infty$$
$$Cov(L_x f(x), g(y)) = L_x Cov(f(x), g(y))$$

*it holds in the m.s. sense that*

$$E(L_x f(x)) = L_x E(f(x))$$
$$Cov(L_x f(x), g(y)) = L_x Cov(f(x), g(y))$$

A proof of this theorem can be found in Christakos (1992).

## 7.4 Covariance Functions

From the theory of stochastic differential equations it is known that linear equations impose restrictions on the distribution. This restrictions are summarised in the following. We are going to have a look at instationary covariance functions and stationary covariance functions.

### 7.4.1 Instationary Covariance Functions

Suppose that our process solves

$$L_x f(x) = k(x)$$

and since $k(x)$ does not depend on the realization of the process $f$ we can conclude that also $L_x f(x)$ does not depend on the realization of the process and that it has variance equal to zero.

$$0 = \mathrm{Var}(L_x f(x)) = L_x L_y C(x,y)|_{x=y}$$
$$0 = \mathrm{Cov}(L_x f(x), f(y)) = L_x C(x,y) \ \forall y$$

Since the partial differential equation holds we can conclude that

$$L_x \mathrm{E}(f(x)) = E(L_x f(x)) = k(x)$$

if we set

$$f_v(x) = f(x) - E(f(x))$$

we get

$$L_x f_v(x) = 0$$

since the expectation of $\mathrm{E} f_v(x) = 0$ for all $x$ and all its derivatives have zero expectation the following formula holds

$$L_x \mathrm{E} f_v(x) = 0$$
$$\mathrm{Var}(L_x f_v(x)) = L_x L_y C(x.y)|_{x=y} = 0$$
$$L_x f_v(x) = \mathrm{E} L_x f_v(x) = L_x \mathrm{E}(f_v(x)) = 0$$

This leads to the following theorem

**Theorem 7.6.** *Following three conditions are equivalent*

- *$f(x)$ solves $L_x f(x) = k(x)$ in the m.s. sense*
- *following formulas holds*

$$L_x E(f(x)) = k(x)$$
$$L_x L_y C(x,y)|_{x=y} = 0$$

- *following formulas holds*

$$L_x E(f(x)) = k(x)$$
$$L_x C(x,y) = 0 \ \forall y$$

### 7.4.2 Stationarity

Let us consider the case of stationary variogram and stationary differential equations. We get the following equation

$$L_x L_y \gamma(x,y) = L_x L_y \gamma(y-x) = L_x L_h \gamma(y-x)$$

If we assume that $L_x$ is even we obtain that

$$L_x(L_h\gamma)(y-x) = L_x(L_h\gamma)(x-y) = (L_h L_h \gamma)(0)$$

and in the same way we can conclude

$$L_x L_y C(x,y) = (L_h L_h C)(0)$$

and we get an equivalent theorem to the theorem above

**Theorem 7.7.** *If $L_x$ is stationary and even then the following conditions are equivalent*

- *$f(x)$ solves $L_x f(x) = k(x)$ in the m.s. sense*
- *The following conditions hold simultaneously*

$$L_x E(f(x)) = k(x)$$
$$(L_h L_h C)(0) = 0$$

- *the following conditions hold simultaneously*

$$L_x E(f(x)) = k(x)$$
$$(L_h C)(0) = 0 \; \forall h$$

where $h$ denotes the distance between $x, y$, for more details on the proofs see Boogaart (2001). Thus the differential equation has a finite number of conditions on the derivatives of $C$ at the origin.

## 7.5 Kriging for Solving differential Equations

The implications on the covariance function $C$ are interesting for kriging in two ways, first they help us to find correct kriging weights, since they help to find appropriate covariance models, and second they imply that the kriging results solve the differential equations.

**Theorem 7.8.** *The universal kriging interpolation solves the linear partial differential equations*

$$L_x f(x) = k(x)$$

*if the covariance function is admissible for the differential equation in the form*

$$L_x L_y C(x,y)|_{x=y} = 0$$

*and the trend functions* $g_1, ..., g_p$ *solve the homogeneous equation*

$$L_x g_i(x) = 0 \; \forall i = 1, ..., p$$

*and the fixed trend part* $g_0$*solves the heterogeneous differential equation*

$$L_x g_0(x) = k(x)$$

A proof can be found in Boogaart (2001), in general $C(x,y)$ can be replaced by any generalised covariance function with respect to the trend function $g(x)$ without changing the results.
An other important application of Kriging is that conditional simulations solve partial differential equations in the m.s. sense, for details on that we refer to Lantuejoul (2002) on the simulation and to Boogaart (2001) on the connection to Kriging and partial differential equations.

## 7.6 The Link to Green Functions

Let us now come back to the regular system in the form

$$\begin{cases} \frac{dZ(t)}{dt} = F(Z(t), Y(t)) \\ Z(t_0) = Z_0 \end{cases}$$

if $Z(t)$ is a solution of that stochastic system, then $Z(t)$ is a multidimensional stochastic process. We are interested now in its spatial distribution knowing that of $Y$, as well as the joint distribution of $Y$ and $Z$ processes. In the Gaussian case, the distribution is fully determined by the first and second moment, the moments of $Z$ and those of $(Y, Z)$ can be calculated either from the expression of $Z$ as a function of the other variables, or as averages over realizations (Chiles Delfiner 1999). We consider following linear system

$$\begin{cases} \frac{dZ(t)}{dt} + A(t)Z(t) = Y(t) \\ Z(t_0) = Z_0 \end{cases}$$

where $A(t)$ is a matrix of continuous functions, $Y(t)$ a m.s. continuous stochastic process with finite second oder moments and $Z_0$ a deterministic or random vector with $\mathrm{Var} Z_0 < \infty$. Since $Y$ is involved linearly in the differential equation, the solution in the m.s. sense is of the form

$$Z(t) = G(t, t_0)Z_0 + \int_{t_0}^{t} G(s,t)Y(s)ds$$

where $G$ is a matrix and $t \geq t_0$. The first term is the product of a deterministic matrix and a random vector independent of $t$, the second term is a linear functional of $Y$. To simplify this, let us suppose a system given in canonical form

$$\begin{cases} \frac{d^m Z(t)}{dt^m} + a_1(t)\frac{d^{m-1}Z(t)}{dt^{m-1}} + ... + a_n(t)Z(t) = Y(t) \\ Z(t_0) = \frac{dZ}{dt}(t_0) = ... = \frac{d^{m-1}Z}{dt^{m-1}}(t_0) = 0 \end{cases}$$

where $a_i$ are continuous functions and $Y(t)$ is a scalar m.s. continuous second order stochastic process. The solution is given by the following equation

$$Z(t) = \int_{t_0}^{t} G(t,s)Y(s)ds$$

where $G(t,s)$ is a Green function representing the response of the system to a unit impulse at time $s$. Denoting by $m_y$ and $C_y$ the mean and covariance function of $Y$, then the mean of $Z$ is given by

$$m_Z(t) = \int_{t_0}^{t} G(s,t)m_y(s)ds$$

and its covariance by

$$C_z(t_1,t_2) = \int_{t_0}^{t_1}\int_{t_0}^{t_1} G(t_1,s_1)G(t_2,s_2)C_Y(s_1,s_2)ds_1 ds_2$$

In the case when the $a_i$ are constants, the Green function is given in the form

$$G(t,s) = p(t-s)$$

and then the solution takes the form

$$Z(t) = \int_{t_0}^{t} p(t-s)Y(s)ds$$

If in addition $Y$ is a stationary process, the trend function of $Y$ constant and $C_y(t_1,t_2) = C_y(t_2 - t_1)$ then

$$m_z(t) = m_y \int_{t_0}^{t} p(t-s)ds$$
$$C_Z(t_1,t_2) = \int_{t_0}^{t_1}\int_{t_0}^{t_2} p(t_1-s_1)p(t_2-s_2)C_Y(s_2-s_1)ds_1 ds_2$$

# 8

# Fully Bayesian Approach

The usual proceeding in geostatistical analysis is to assume that the model is absolutely known, this means that we have full knowledge about the trend function and the covariance function. For example the trend is modelled by a low order polynomial and the covariance function is plugged in by some estimation and it is always assumed that the data are normally distributed, so that a decomposition of the random field makes sense. When we assume the data to be normally distributed, then the results of spatial interpolation crucially depend on the estimation of the covariance function. Instead of the true covariance function, we only use an empirical estimation of it, and the usual practice is then 'plug in Kriging', this is only a estimation of the unknown Kriging predictor, and the BLUP optimality is therefore no longer valid. It is known from Christensen (1991) that the actual mean squared error of the prediction is larger than the theoretical mean squared error of prediction resulting from the assumption that the covariance function has been specified exactly and also larger than the reported one. In this chapter we are going to present an alternative Bayesian approach which models the uncertainty by means of suitable posterior probability distributions for the parameters of a flexible class of nested Matern covariance functions.

## 8.1 Uncertainty due to the Model Assumptions

Implicit in all geostatistical applications is the assumption of ergodicity of the random function considered, but this assumption can not be verified in practice. So the covariance function can not be determined completely from the knowledge of some observations of the random field. Moreover a big problem is the decision on the scale, this means that we have to specify which part of the observed variable is in global scale modelled by some low order polynomial and the local scale modelled by the random field. An other source of uncertainty results from the various approaches to modelling anisotropy. After the choice of the model assumptions there is a large amount of uncertainty in the estimation of the covariance function and in the numerical fitting methods of the covariance functions. As we have seen in earlier chapters there is a large number of covariance models, different approaches to fitting and so on. One very critical evaluation is the choice of the so called nugget effect, since the nugget requires extrapolation of the variogram into the origin.

## 8.2 The full Bayesian approach

The Bayesian approach provides a general methodology for taking account of uncertainty with respect to the model and its components. This is specially important for the specification of the model parameters. We use the four parameter Matern variogram model $\gamma_M(h,\theta)$ with $\theta = (c_0, c, a, \nu)$ where $c_0$ is the nugget, $c$ the sill, $a$ the range and $\nu$ the smoothness parameter. The full Bayesian approach requires a completely specified distributional model, this means beside the data displayed in the likelihood function, we need to specify probability distributions for the trend parameter $'\beta$ and the covariance parameter $\theta$. The property of BLUB is only given when the covariance function $C$ is absolutely known, but in practice $C$ is estimated and plugged in as the true covariance function. Berger (1988) describes a way out of this dilemma by means of reference prior, but in practice they are hardly to compute and most of Bayesians use only non informative reference priors for the covariance parameters. The model parameters should be modelled by appropriate prior probability distributions and then the so called predictive density offered by the Bayesian paradigm gives us a complete probability distribution for the predicted values, instead of only having Kriging maps of the interpolation and the variance we then give the predictive density in each point and so the quantiles and other statistical measures. This is important for finding so called hot spots.
The Bayes optimal prediction of $Z(x_0)$ at an unobserved location $x_0$ in the domain based on the predictive density is given by

$$p(Z(x_0)|Z) = \int_{\Theta} \int_{\mathcal{B}} p(Z(x_0)|\beta, \theta, Z) p(\beta, \theta|Z) d\beta d\theta$$

where the first term under the integral sign is the conditional probability density and the second one is the posterior probability density. Above formula is the conditional probability density of $Z(x_0)$ given the data $Z = (Z(x_1), ..., Z(x_n))^T$ averaged over all trend parameters $\beta = (\beta_1, ..., \beta_n) \in \mathcal{B}$ and variogram parameters $\theta \in \Theta \subseteq \mathbb{R}^4$, where the averaging is done with respect to the posterior probability density of $\beta, \theta$ for given data $Z$. The posterior probability density of the parameters is obtained according to Bayes's theorem as

$$p(\beta, \theta|Z) = \frac{p(Z|\beta, \theta) p(\beta, \theta)}{\int_{\Theta} \int_{\mathcal{B}} p(Z|\beta, \theta) p(\beta, \theta) d\beta d\theta}$$

where $p(Z|\beta, \theta)$ is the likelihood function of the data and $p(\beta, \theta)$ is the prior probability density of the parameters $\beta$ and $\theta$.The density $p(Z(x_0)|\beta, \theta, Z)$ is a normal distribution with mean $\mu = f(x_0)^T\beta$ and variance $c_0+c-k_0^T K(\theta)^{-1}k_0$ where $k_0$ is given by $(c_0 + c - \gamma(x_0 - x_1; \theta), ..., c_0 + c - \gamma(x_0 - x_n; \theta))$. Our approach for specifying the prior probability distribution of the covariance parameters $\theta$ is done by simulation of an possible empirical covariance parameters. Instead of $p(\theta)$ we generate the posterior density given the data by conditional simulations. We factor the joint posterior density of the trend and covariance parameter according to

$$\begin{aligned} p(\beta, \theta|Z) &= p(\beta|\theta, Z)p(\theta|Z) \\ &= k\, p(z|\beta, \theta)p(\beta|\theta)p(\theta|Z) \end{aligned}$$

where $k$ denotes a normalisation constant.

### 8.2.1 Algorithm for the posterior Density

The posterior density $p(\theta|Z)$ of the variogram parameters is computed by

- Generate $m$ simulated data sets from the random function $Z(.)$, conditional on the actual observations $Z$.
- Fit a theoretical variogram $\gamma$ to the estimated variogram of the simulated random functions, we obtain sufficiently many realisations of the posterior density
- Compute the predictive density given in above formula
- Make a back transformation to the original scale by the Jacobian.

Assuming a normally distributed prior for $\beta$ the posterior density $p(Z(x_0)|Z, \theta)$ is normally distributed with mean ordinary Bayes Kriging predictor and variance ordinary Bayes Kriging variance.

# 9

# Practical Example

Based on the theoretical results in our work we are now going to present a practical example. The data set used in this chapter is from the region of Gomel (State of Belarus) ten years after the terrible Chernobyl catastrophe. In April 1986 the World's worst nuclear power accident occurred at Chernobyl in the former UdSSR (now Ukraine).

## 9.1 What happened in April 1986

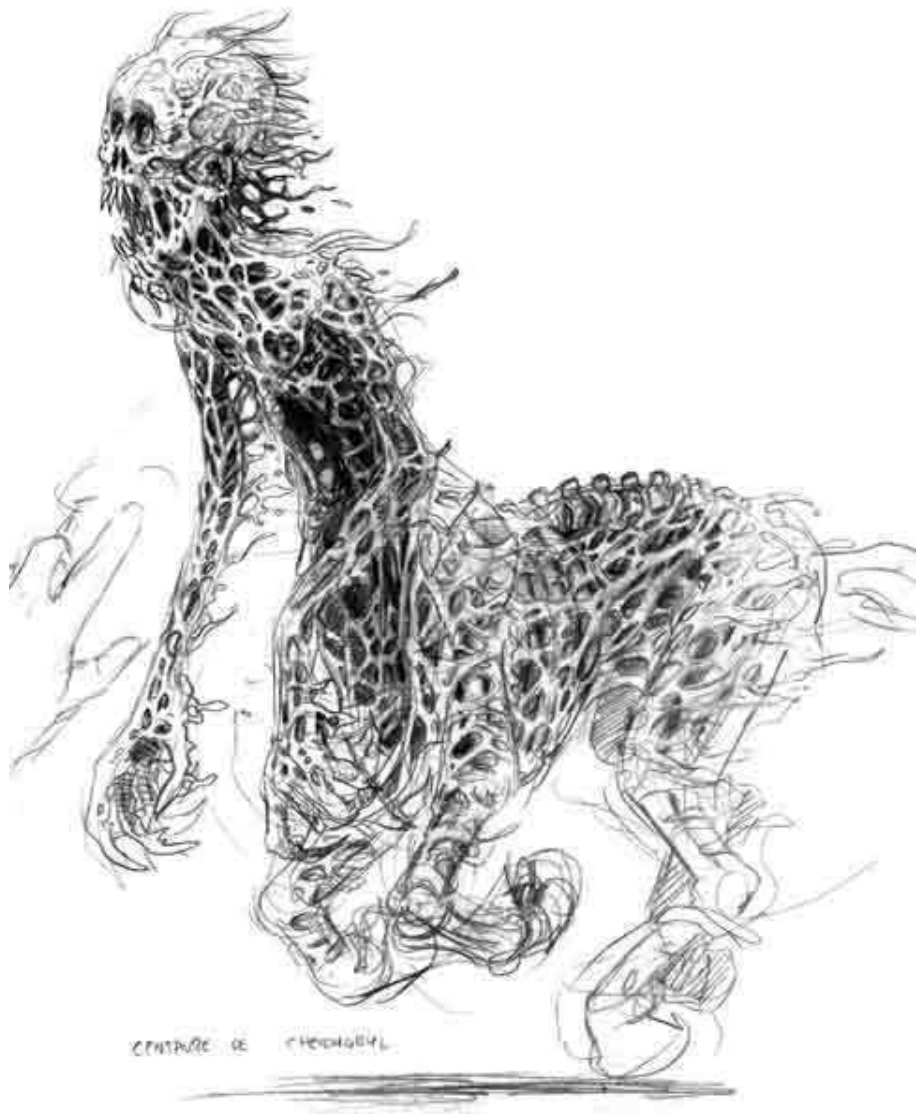

**Fig. 9.1.** The death of our science, by David Cochard

We believe that one picture tells more about the story, more than hundreds of lines in this context, when the nuclear death was coming to all the inhabitants of the former beautiful region of Gomel, and even now, he is not willing to leave that region.

## 9.2 The Chernobyl Accident

The Chernobyl nuclear power plant, located 100 kilometers north of Kiev, had 4 reactors and the testing reactor number 4, where numerous safety procedures were disregarded. At one April day at 1:23am the chain reaction in the reactor became out of control, creating explosions and a fireball which blew off the reactor's heavy steel and high concrete lid. The disaster destroyed the

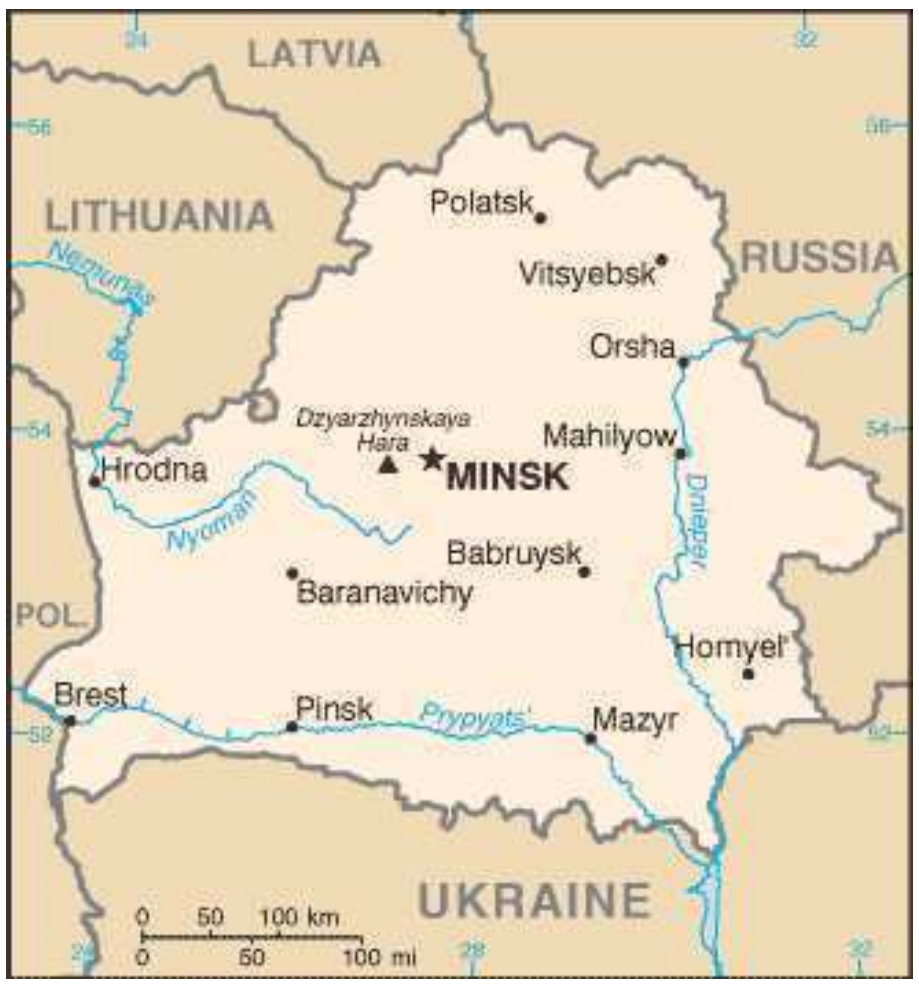

**Fig. 9.2.** The region of Gomel

Chernobyl-4 reactor and killed 30 people, including 28 from radiation exposure. A further 209 on site were treated for acute radiation poisoning and among these, 134 cases were confirmed. Large areas of Belarus, Ukraine, Russia and beyond were contaminated in varying degrees. The Chernobyl disaster was a unique event and the only accident in the history of commercial nuclear power where radiation-related fatalities occurred. But we should always mention, that there are lot of other nuclear power plants in our surrounding, hopefully working well. Now 18 years after the accident still over 3 millions of children are suffering on it, the area around the reactor is still highly contaminated, nature is dead and there is no wildlife.

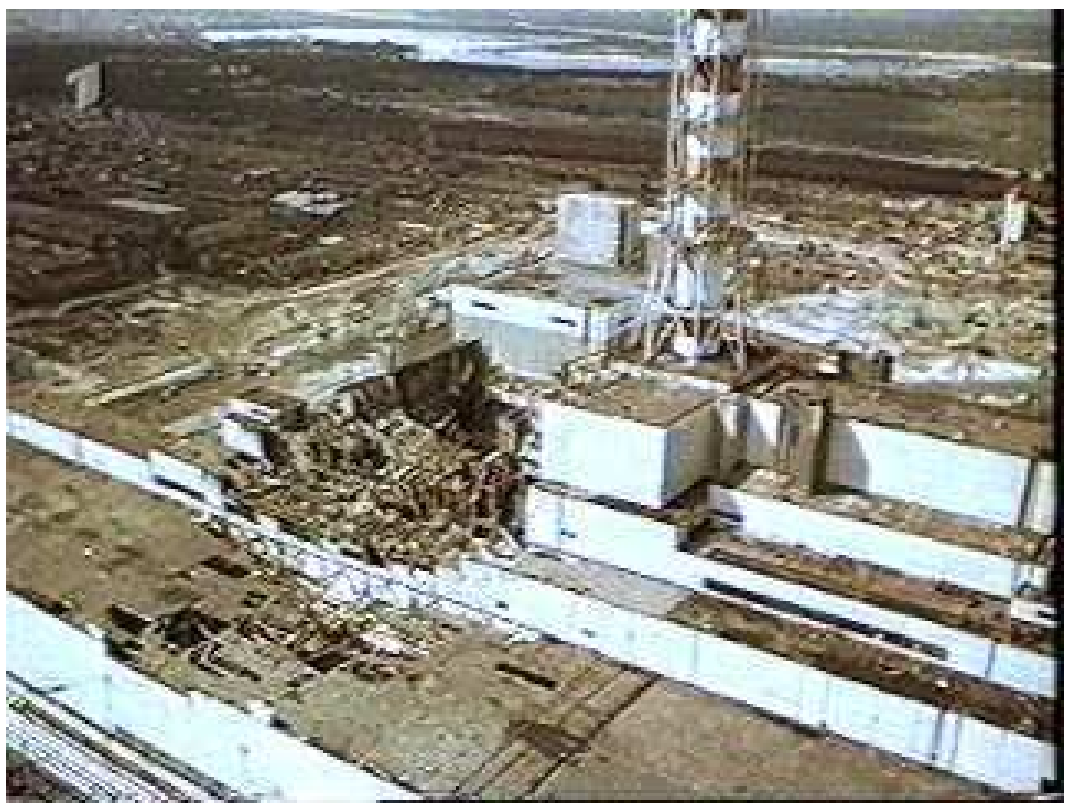

**Fig. 9.3.** The reactor number 4 after the accident

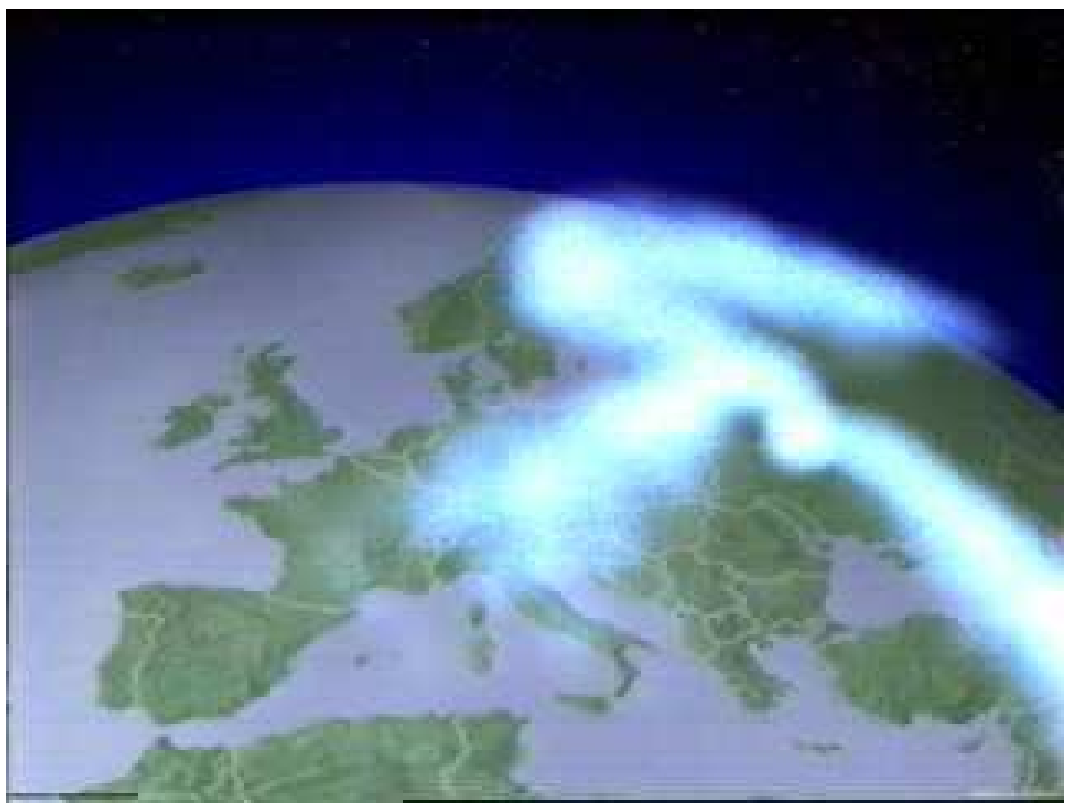

**Fig. 9.4.** Radioactive cloud in the northern hemisphere

## 9.3 Work on the Chernobyl Data Set

In the region of Gomel the Russian government built a network where the concentration of the radioactivity is still measured. The dataset we use are measurements of Cs137, a radioactive isotope.

### 9.3.1 Explorative Statistics

By using simple statistics we are able to conclude that the data set is not normally distributed.

```
> summary(cs137aver)
 CS137.AVER
 Min.   : 0.040
 1st Qu.: 0.710
```

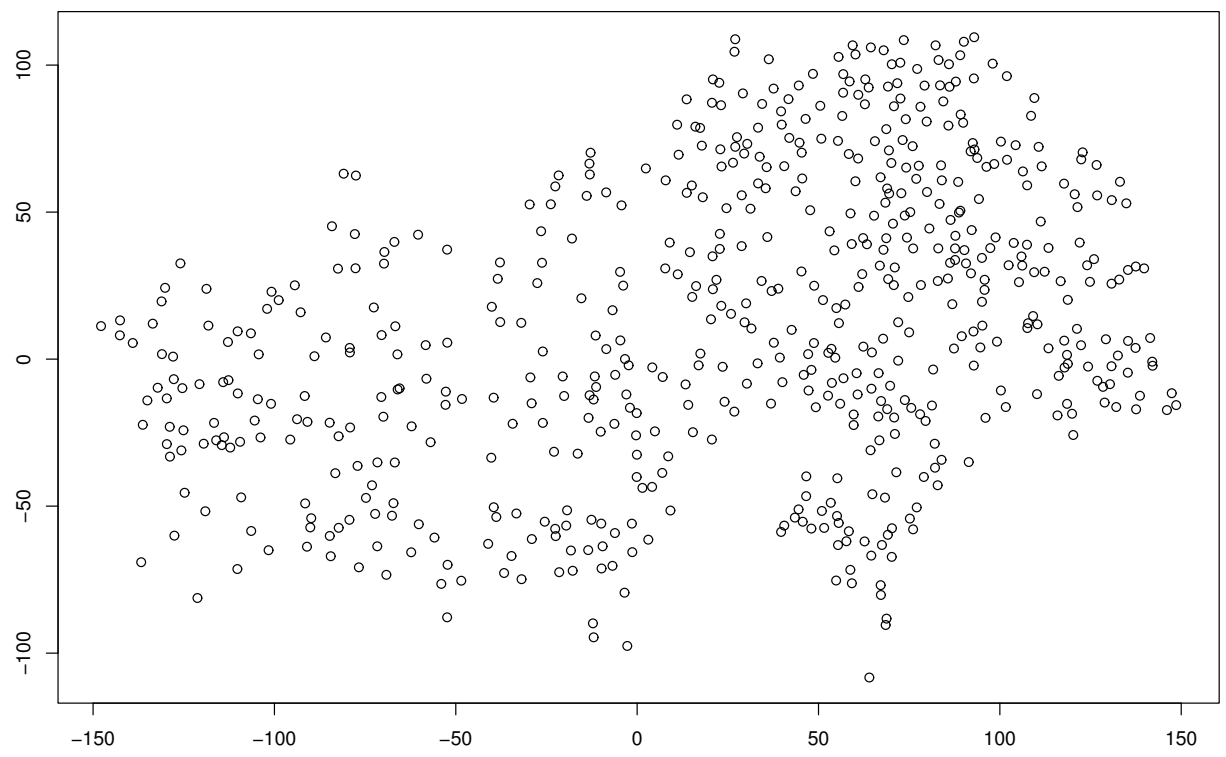

**Fig. 9.5.** Sampling network

```
 Median : 1.920
 Mean   : 4.995
 3rd Qu.: 5.855
 Max.   :61.390

> stem(cs137aver)

  The decimal point is 1 digit(s) to the right of the |

  0 | 000000000000000000000000000000000...0000000000000000000+331
  0 | 55555555555555...5666...66667777777777777777777777777+15
  1 | 0000000111111111222222233333444
  1 | 5556677778889999
  2 | 0111122222223444
  2 | 667788
  3 | 011223344
  3 | 67
  4 | 014
  4 |
  5 |
  5 | 5
  6 | 1
```

So we have to apply the Box Cox transformation and can then assume that the transformed data is normally distributed, where the concepts of spatial interpolation with Kriging makes sense.

```
> summary(log(cs137aver))
   CS137.AVER
 Min.   :-3.2189
```

```
 1st Qu.:-0.3425
 Median : 0.6523
 Mean   : 0.6641
 3rd Qu.: 1.7671
 Max.   : 4.1172

> ks.test(log(cs137aver),"pnorm", mean(log(cs137aver)),sd(log
                                                  (cs137aver)))

        One-sample Kolmogorov-Smirnov test

data:  log(cs137aver)
D = 0.0271, p-value = 0.7791
alternative hypothesis: two.sided
```

### 9.3.2 Geostatistics

So by the transformation to the log scale we can apply geostatistics to our data set. We use the 'geoR' package from Paulo J. Ribeiro Jr. and Peter J. Diggle which can be found on the CRAN archive.

```
> library(geoR)

---------------------------------------------------------
geoR - functions for geostatistical data analysis
geoR version 1.4-7 (built on 2004-05-04) is now loaded
---------------------------------------------------------

> data <- as.geodata(cbind(coord,cs137aver))
> str(data)
List of 2
 $ coords: num [1:591, 1:2]  64.01  -2.73 -11.99  68.51 -12.19 ...
  ..- attr(*, "dimnames")=List of 2
  .. ..$ : chr [1:591] "1" "2" "3" "4" ...
  .. ..$ : chr [1:2] "EASTING" "NORTHING"
 $ data  : num [1:591]  3.05  3.06  4.33  5.74 12.84 ...
 - attr(*, "class")= chr "geodata"
```

Now we are going to produce a plot of the data locations and values in each location. This is done by the function 'plot.geodata', and the plot shows a $2 \times 2$ window with data locations (top plots) and data versus coordinates (bottom plots). For an object of the class 'geodata' the plot is produced by the command:

```
> plot(data)
```

For a better graphical output we used following

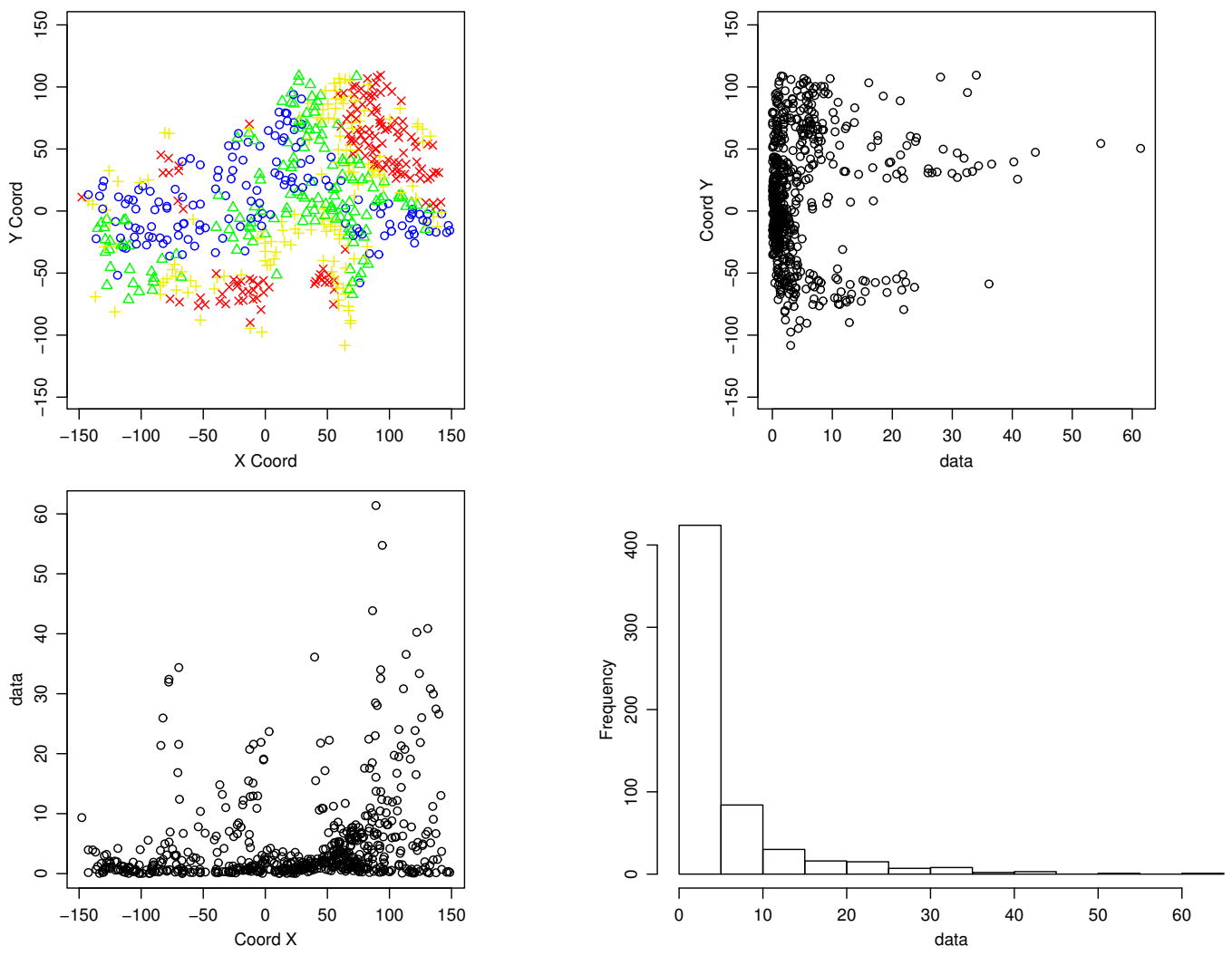

**Fig. 9.6.** Data locations and values

```
> par.ori <- par(no.readonly = TRUE)
> par(mfrow = c(2,2))
> points(s100, xlab = "Coord X", ylab = "Coord Y")
> points(s100, xlab = "Coord X", ylab = "Coord Y", pt.divide =
          "rank.prop")
> points(s100, xlab = "Coord X", ylab = "Coord Y", cex.max = 1.7,
+         col = gray(seq(1, 0.1, l=100)), pt.divide = "equal")
> points(s100, pt.divide = "quintile", xlab = "Coord X", ylab =
          "Coord Y")
> par(par.ori)
```

which produce figure 'Data location and values 2'

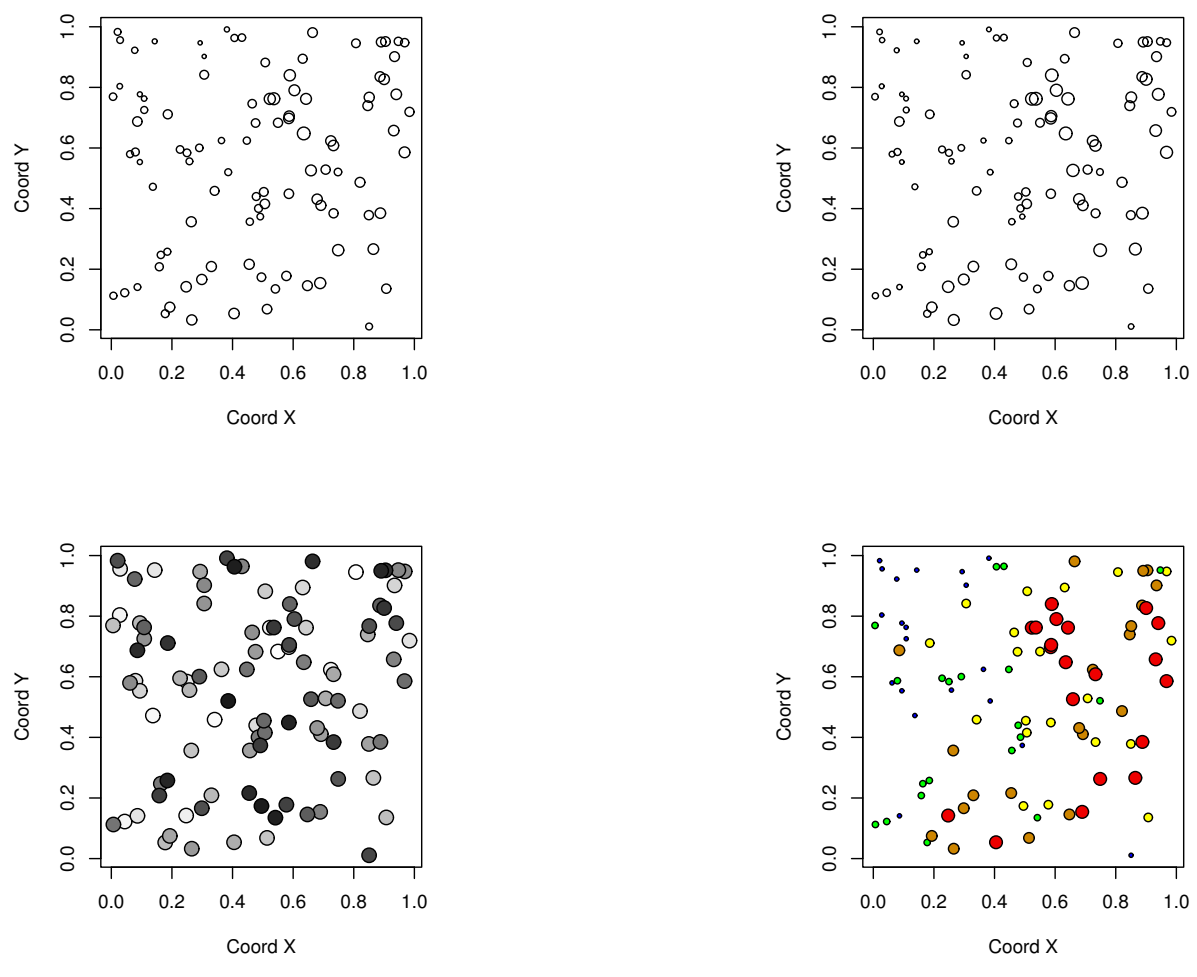

**Fig. 9.7.** Data location and values 2

### 9.3.3 Empirical Variogram Estimation, Theoretical Variograms

Now we are going to estimate the empirical variograms, that are calculated using the function 'variog' from the package 'geoR'. There are options for the classical or modulus estimator as introduced in earlier chapters. The results are returned as variogram clouds, binned and smoothed variograms.

```
> cloud1 <- variog(data, option = "cloud", max.dist=150)
> cloud2 <- variog(data, option = "cloud", estimator.type =
+                    "modulus", max.dist=150)
> bin1 <- variog(data, uvec=seq(0,150,l=11))
> bin2  <- variog(data, uvec=seq(0,150,l=11), estimator.type=
+                    "modulus")
>
>
> par(mfrow=c(2,2))
> plot(cloud1, main = "classical estimator")
> plot(cloud2, main = "modulus estimator")
> plot(bin1, main = "classical estimator")
> plot(bin2, main = "modulus estimator")
> par(par.ori)
>
```

Furthermore, the points of the variogram clouds can be grouped into classes of distances, in the geostatistical language called 'bins', those we are now going to display with a box-plot for each bin.

```
> bin1 <- variog(data,uvec = seq(0,150,l=11), bin.cloud = T)
> bin2 <- variog(data,uvec = seq(0,150,l=11), estimator.type =
```

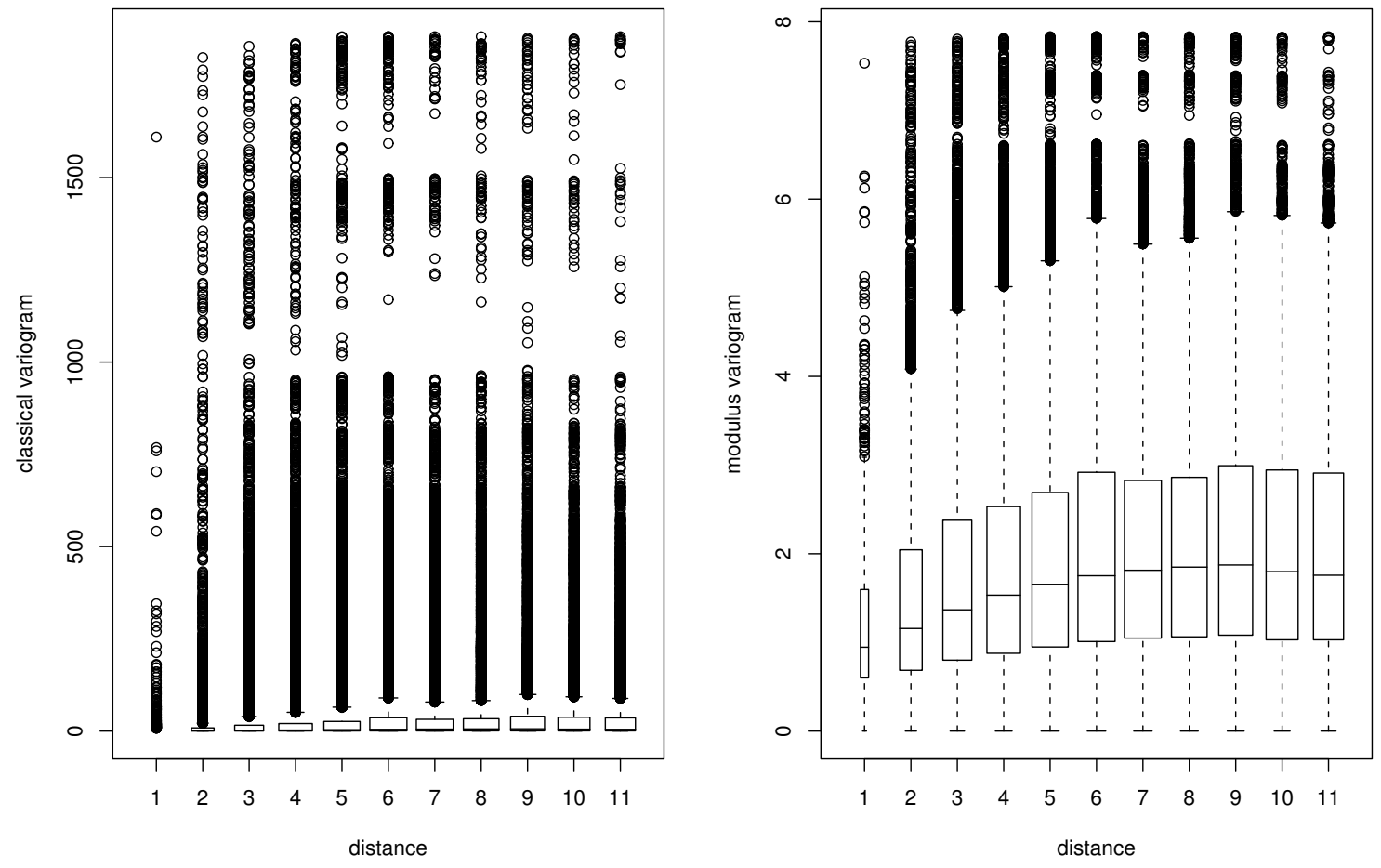

**Fig. 9.9.** Box plots for variogram clouds

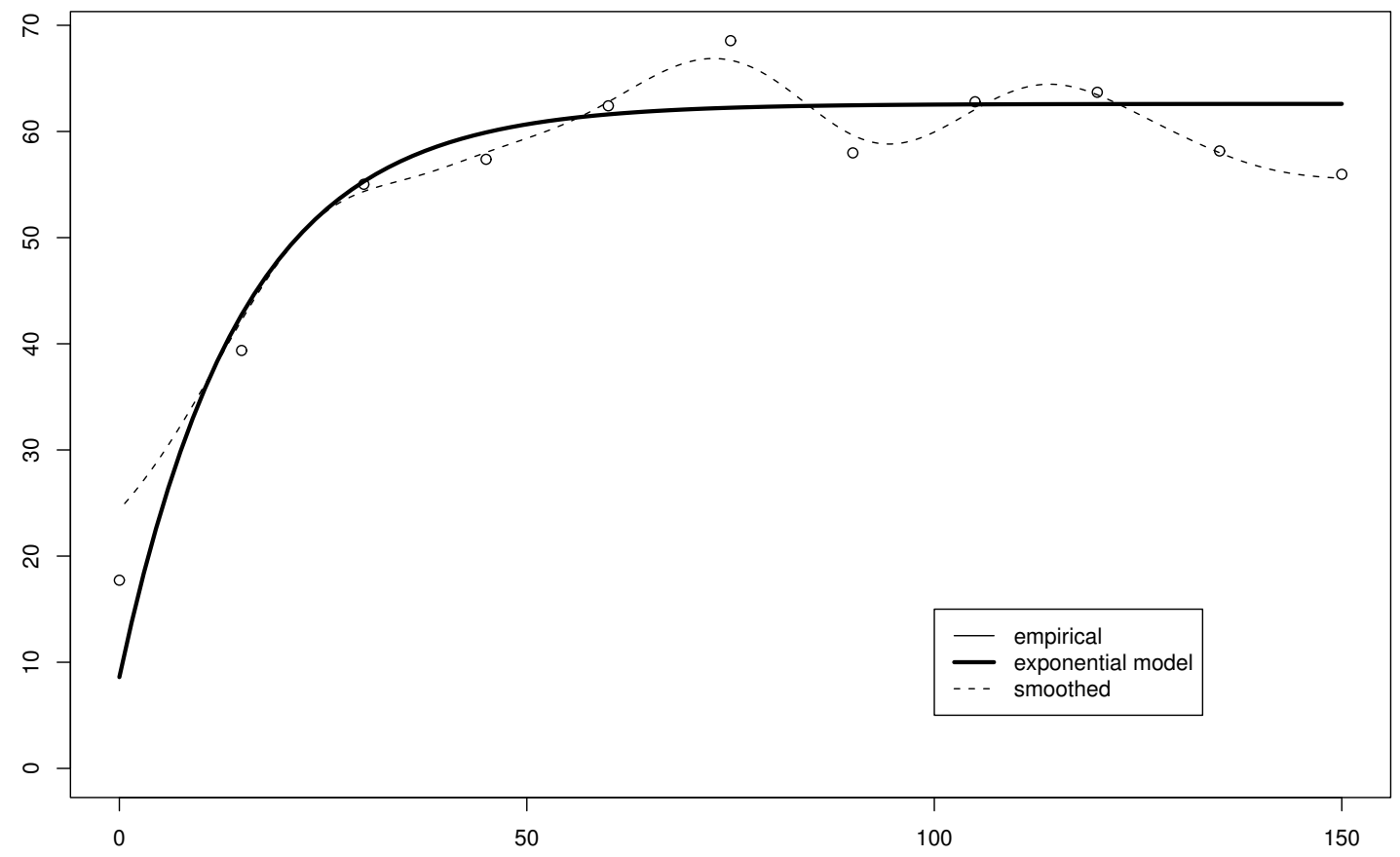

**Fig. 9.10.** Visual comparison of theoretical and empirical variogram

```
        (0.785 radians)
        tolerance angle = 22.5 degrees (0.393 radians)
variog: computing variogram for direction = 90 degrees
        (1.571 radians)
        tolerance angle = 22.5 degrees (0.393 radians)
variog: computing variogram for direction = 135 degrees
```

```
        (2.356 radians)
        tolerance angle = 22.5 degrees (0.393 radians)
variog: computing omnidirectional variogram
>
> plot(vario.4, lwd=2,xlab="",ylab="")
```

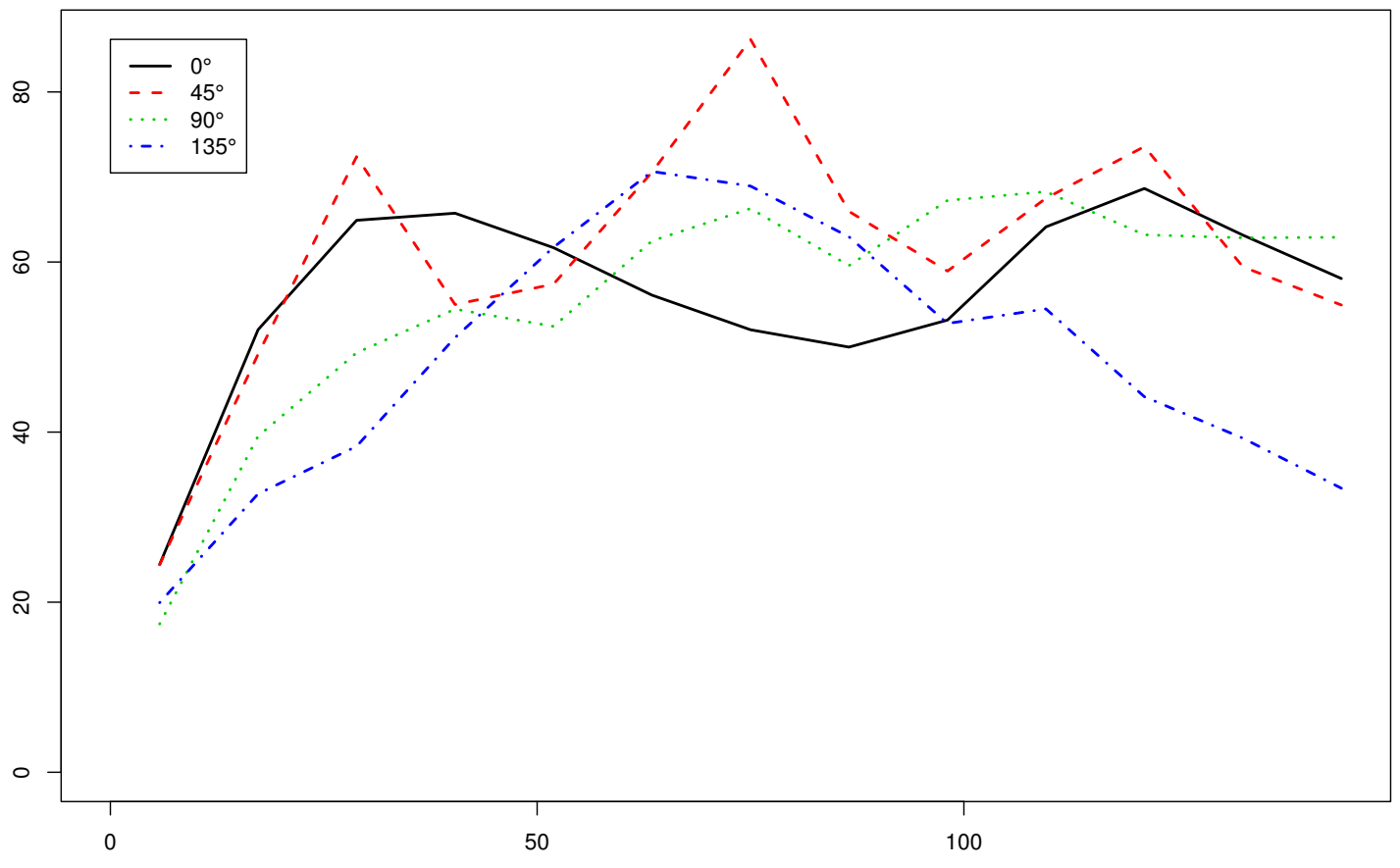

**Fig. 9.11.** Variograms in different directions

In the following we are interested in the possible 'bounds' of our variogram estimation, for that we used following code

```
> data.ml <- likfit(data, ini = c(0.5, 0.5), fix.nugget = TRUE)
---------------------------------------------------------------
likfit: likelihood maximisation using the function optimize.
likfit: Use control() to pass additional
         arguments for the maximisation function.
        For further details see documentation for optimize.
likfit: It is highly advisable to run this function several
        times with different initial values for the parameters.
likfit: WARNING: This step can be time demanding!
---------------------------------------------------------------
likfit: end of numerical maximisation.
>
> data.vario <- variog(data, max.dist = 150)
variog: computing omnidirectional variogram
>
> data.env <- variog.model.env(data, obj.v = data.vario,
+                                   model.pars = data.ml)
variog.env: generating 99 simulations (with  591 points each)
```

```
using the function grf
variog.env: adding the mean or trend
variog.env: computing the empirical variogram for the 99 simulations
variog.env: computing the envelops
>
> plot(data.vario, env = data.env)
```

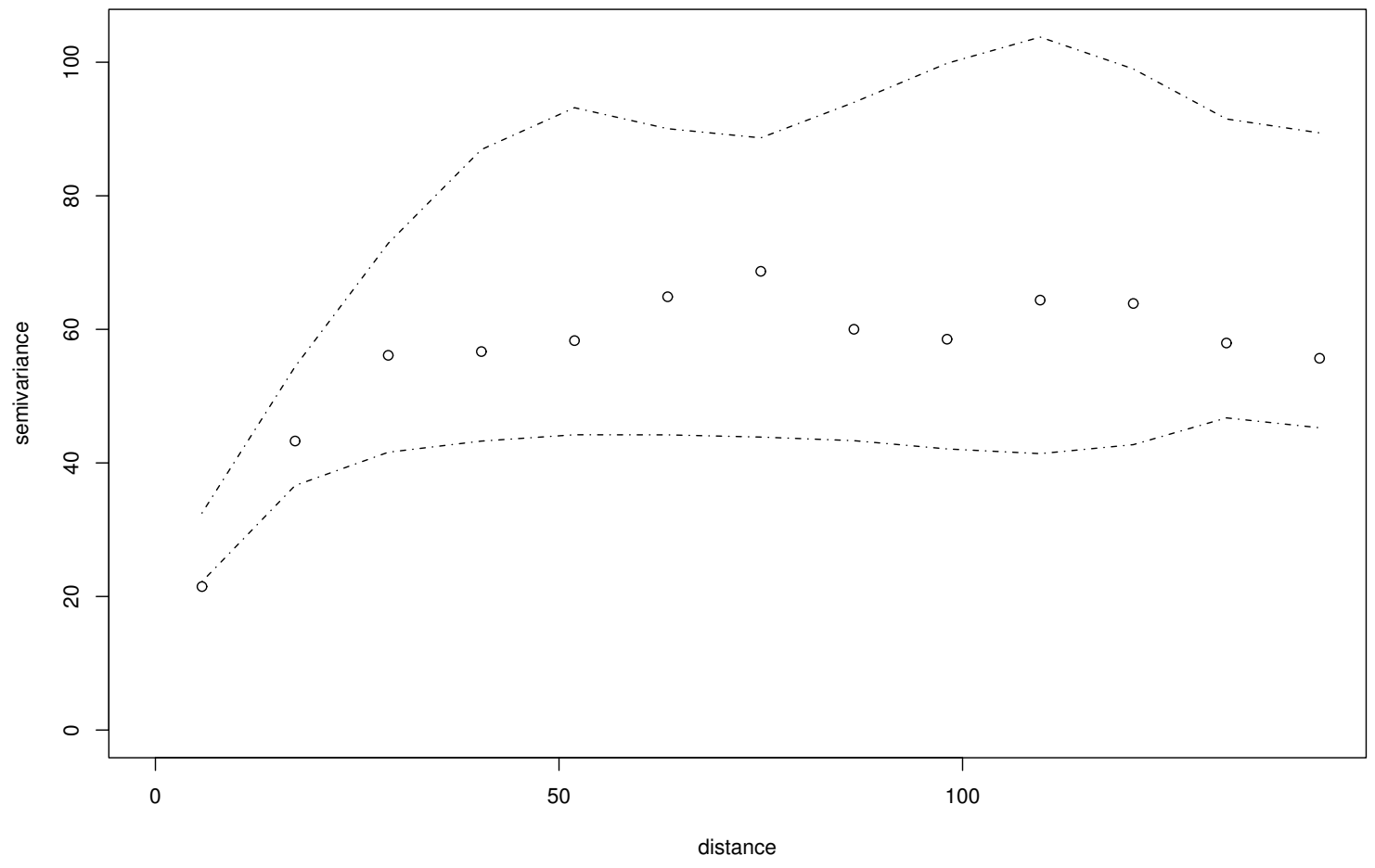

**Fig. 9.12.** Possible bounds for the variograms with geoR

In comparison to that we can use our own coding for the simulation of possible variograms, which is approximately ten times faster than the 'geoR' routines, we also need the library 'sgeostat':

```
>
> coord<-daten[,1:2]
> cs137aver<-daten[,3]
> gomel.pts<-point(logCS137aver,x="EASTING",y="NORTHING")
> gomel.prs<-pair(gomel.pts,num.lags=21,type='anisotropic',theta=0,
+                dtheta=22.5, maxdist=200)
> Gdist<-as.matrix(dist(coord))
> GCov<-matern.cov(Gdist,range=122.79,sill=2.4523,nugget=0.0661,
+                 nue=0.5)
> LW1<-chol(GCov)
> varsim<-simulate(daten[,1:2],L=LW1,gomel.pts,gomel.prs,n=100)
>
> simplot(varsim)
```

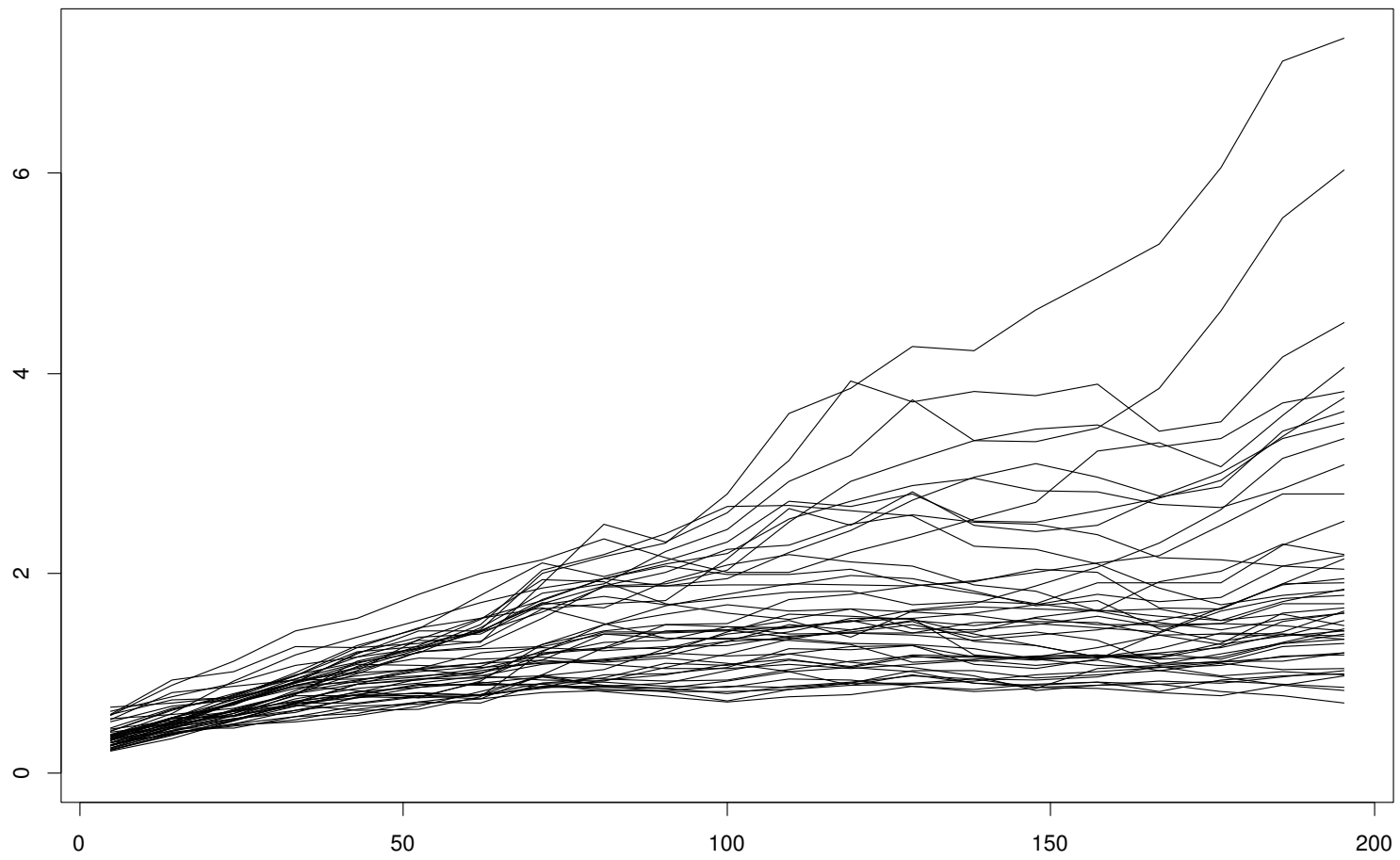

**Fig. 9.13.** Possible bounds for the variogram with our routines

### 9.3.4 Kriging Maps

In the following part of our practical example we are going to present the kriging maps under the use of the package 'sgeostat'. We first present ordinary kriging maps with the exponential variogram function and second we will compare these results with ordinary kriging based on the Matern variogram function.

```
> range(northing)
[1] -108.2907  109.4745
> range(easting)
[1] -147.8045  148.6105
> grid <- list(x=seq(min(easting),max(easting), length=50),
+              y=seq(min(northing),max(northing),length=37))
> grid$xy <- data.frame(cbind(c(matrix(grid$x,length(grid$x),
+   length(grid$y))),c(matrix(grid$y,length(grid$x),
+   length(grid$y),byrow=T))))
> colnames(grid$xy) <- c("x","y")
> library(sgeostat)
Loading required package: mva
Loading required package: tripack
>
> grid$pts <- point(grid$xy)
> grid$krige.expo <- krige(grid$pts, gomel.pts, 'cs137.aver',
+                    gomel.vmod.expo, maxdist= 45, extrap=T)

Using points within 45 units of prediction points.
  Predicting.................................

> image(grid$x, grid$y, matrix(grid$krige.expo$zhat,
+       length(grid$x), length(grid$y)), xlab= "easting (km)",
```

```
+       ylab= "northing (km)")
> contour(grid$x,grid$y,matrix(grid$krige.expo$zhat,
+        length(grid$x), length(grid$y)),add=T)

> range(grid$krige.expo$zhat, na.rm=T)
[1] -0.3504042 50.7699237
```

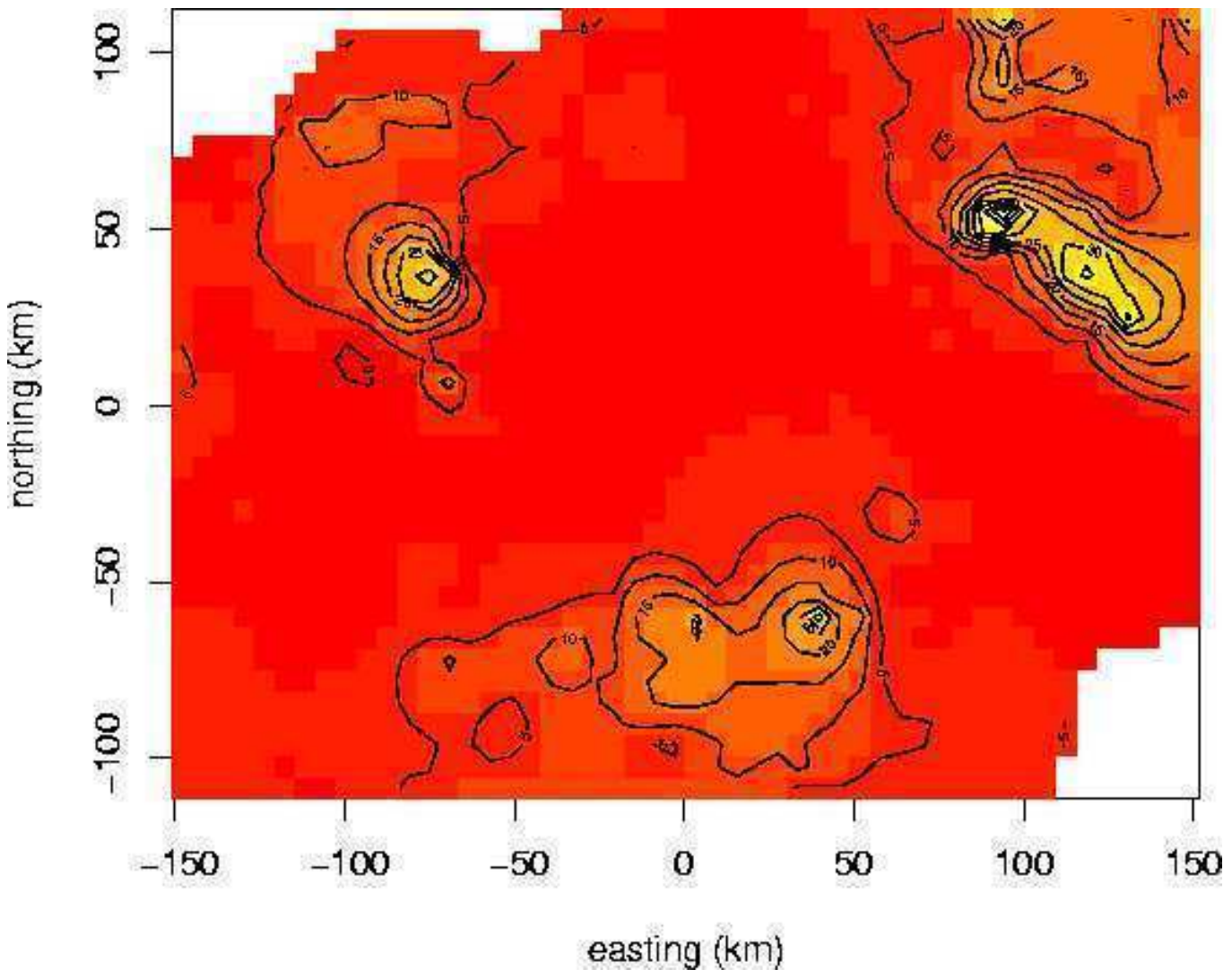

**Fig. 9.14.** Ordinary Kriging with exponential variogram function, see Breitenecker(2001)

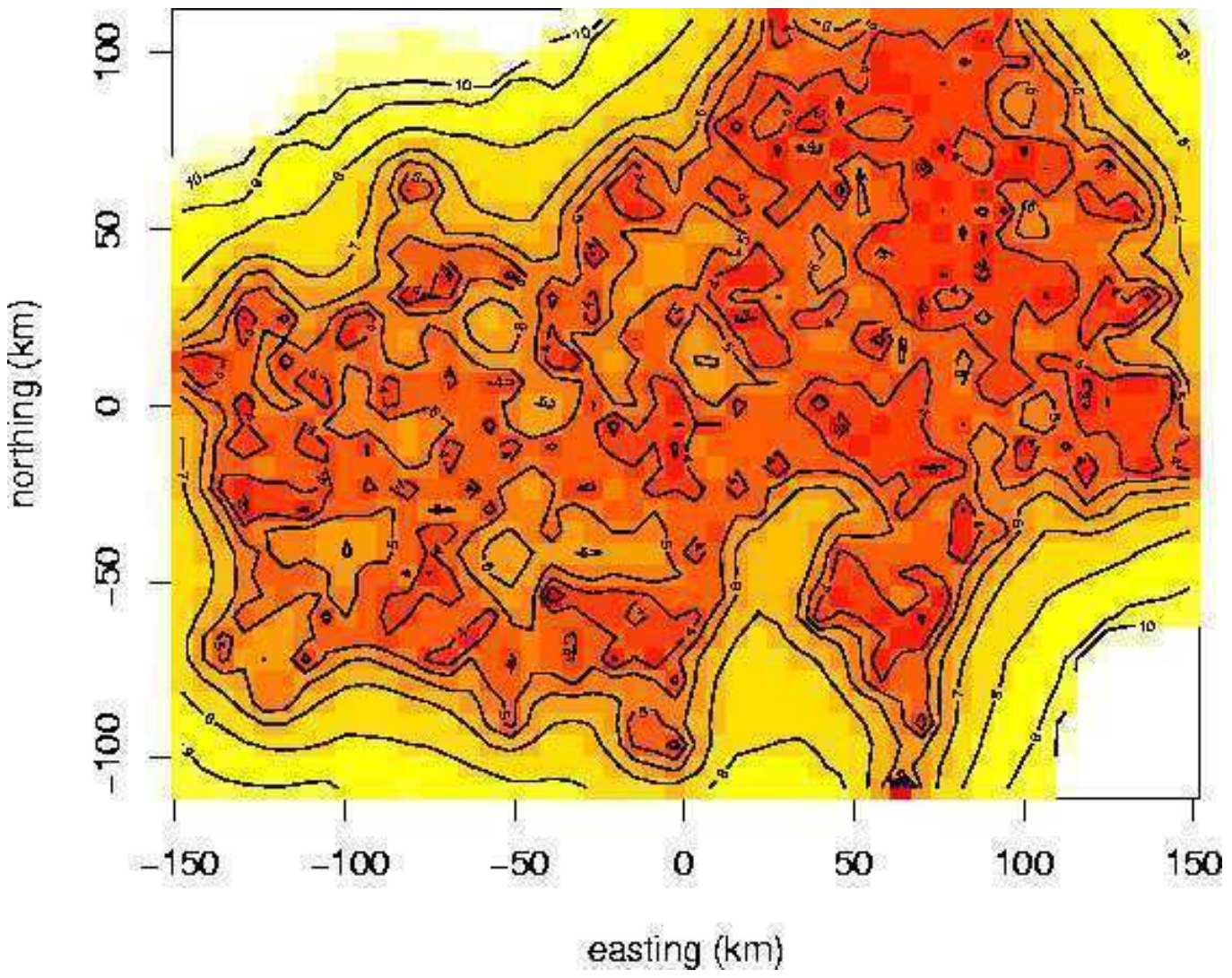

**Fig. 9.15.** Ordinary Kriging variance map with exponential variogram function, see Breitenecker(2001)

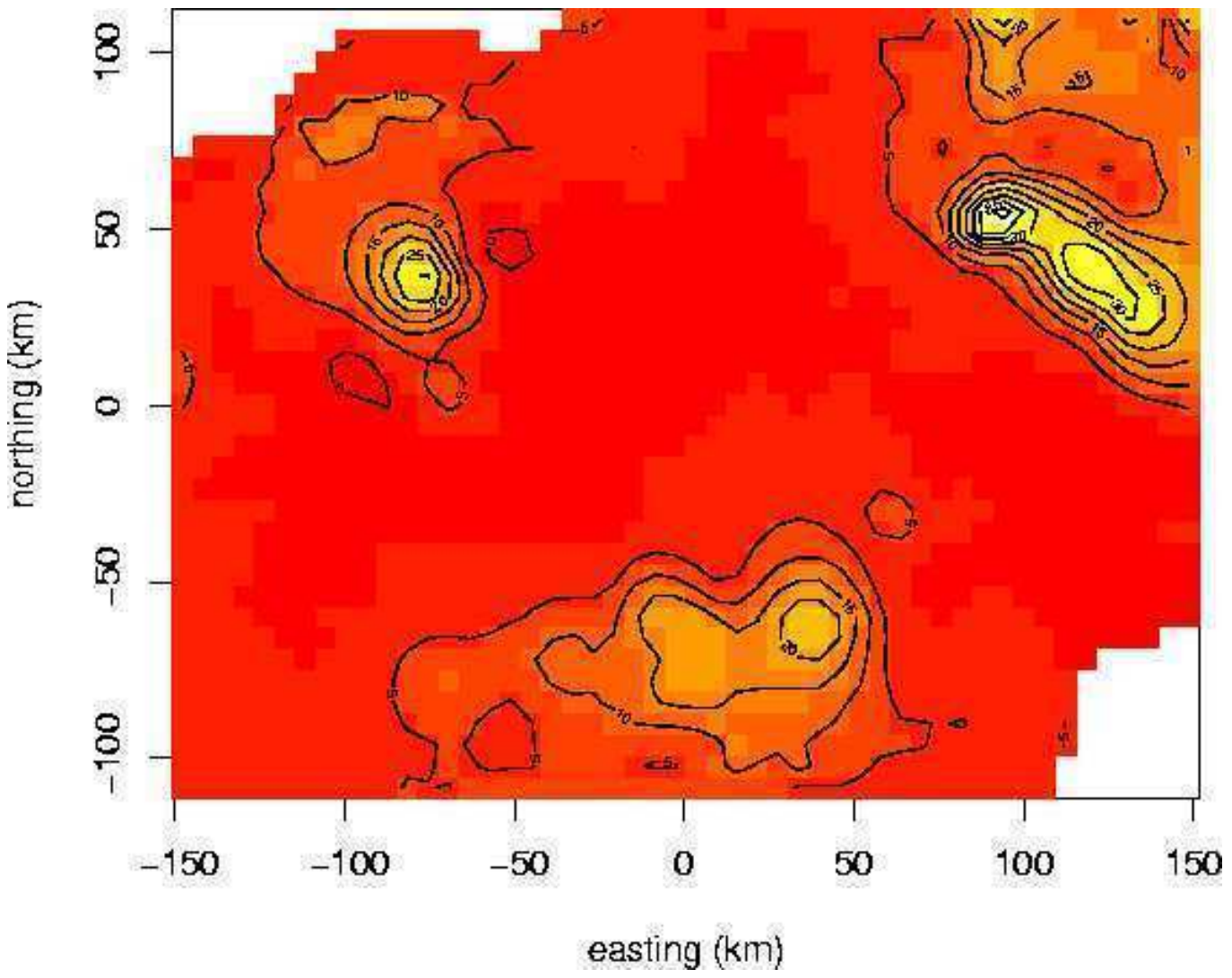

**Fig. 9.16.** Ordinary Kriging with the Matern variogram function, see Breitenecker(2001)

### 9.3.5 Bayesian Approach

Now we are going to produce Bayes maps by using the function our implementation. We can see that this method is a very smooth interpolation of the given data set. In comparison to Ordinary Kriging we can conclude that the variances in Bayes Kriging map is higher than in the Ordinary Kriging map. Based on the discussion in further chapters we are now going to use the predictive density for interpolating the data, as an example we are now going to evaluate the predictive density at a given hot spot, namely at $HS = (x_0, y_0) = (-80, 40)$.

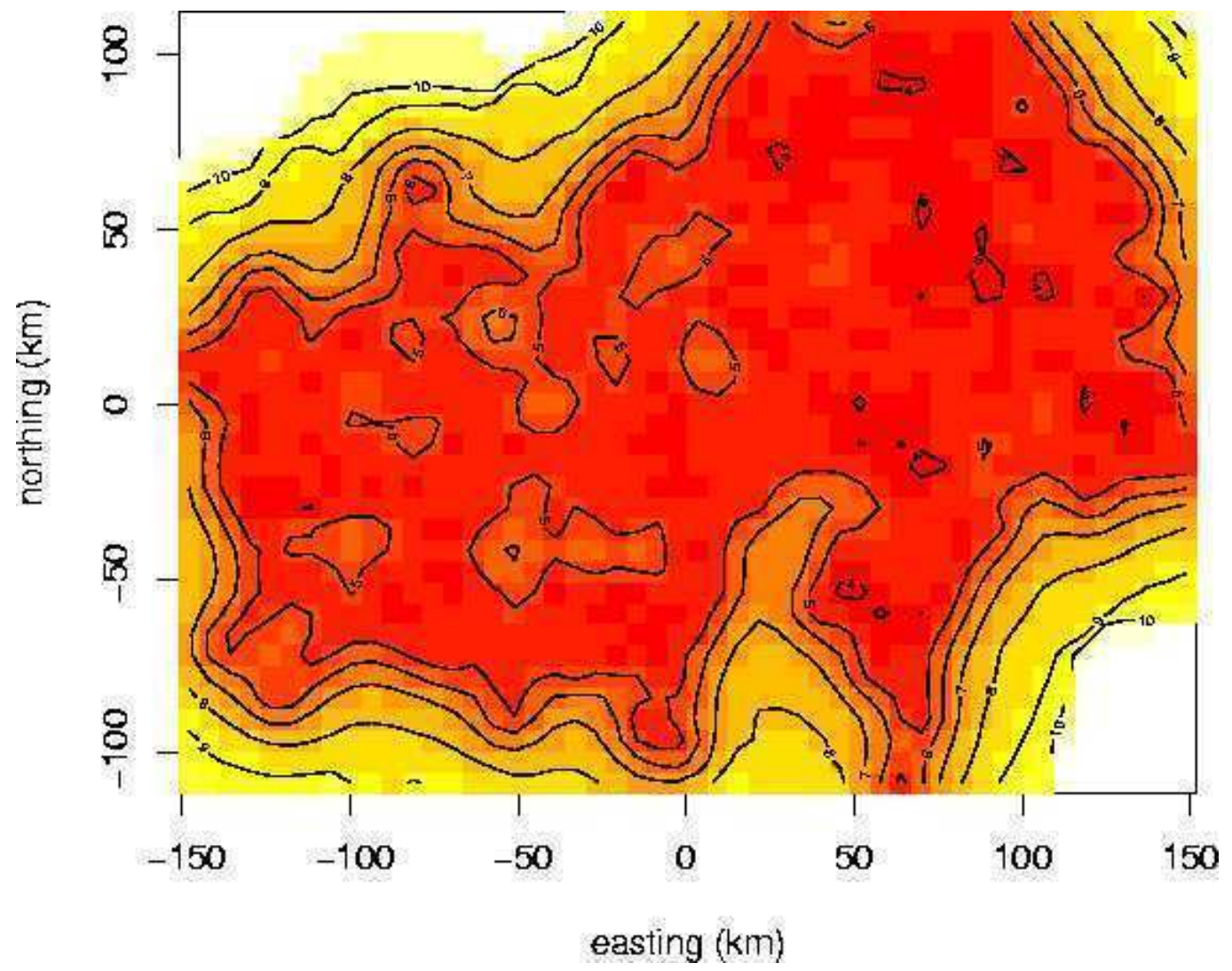

**Fig. 9.17.** Ordinary Kriging variance map with the Matern variogram function, see Breitenecker(2001)

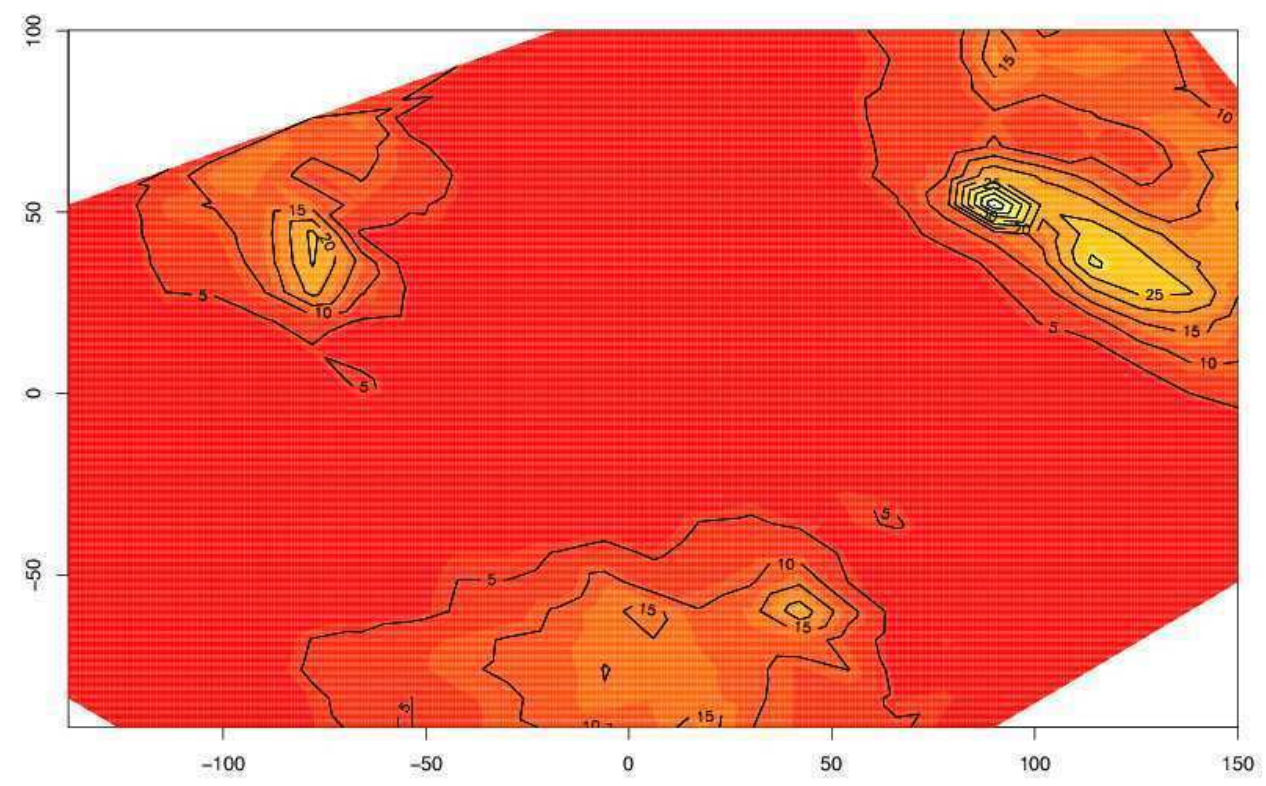

**Fig. 9.18.** Ordinary Bayes Kriging with Matern variogram function

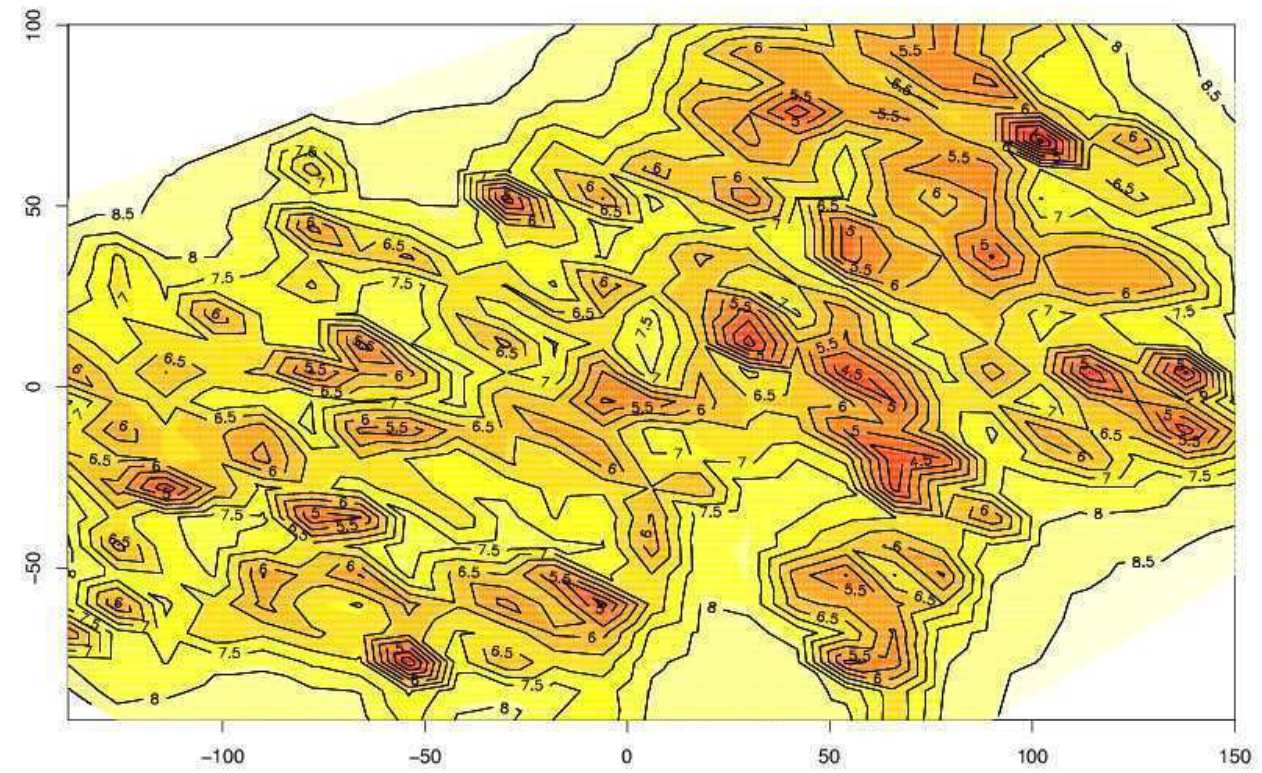

**Fig. 9.19.** TMSEP of Ordinary Bayes Kriging with Matern variogram function

```
val <- bayesKriging(-80,40,0,10000,delta0,40,Coord,obsVal)
```

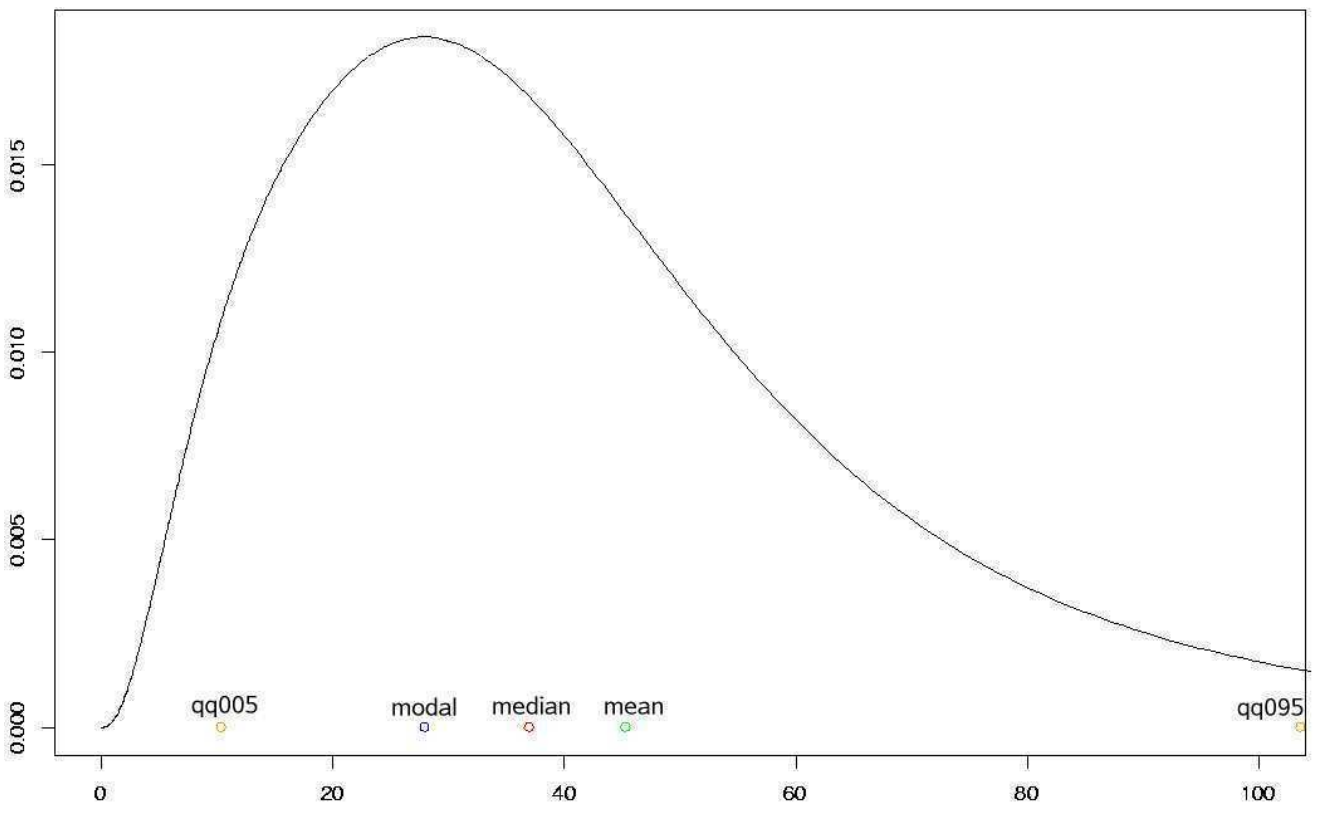

**Fig. 9.20.** Predictive Density at a hot Spot

Evaluating the predictive density at each point on a given grid over the domain, we are able to give maps of the different quantiles, like here illustrated by the 95 percent and 5 percent quantiles. So if we are interested in threshold values, e.g. the 75 percent threshold value, then we simply produce a map of the 75 percent quantiles from the predictive distributions computed at a corresponding grid of points. For a uncertainty map we use the approximation of the standard deviation by the interquartile range divided by $a = 1.45$. The interquartile range is defined as the difference between the 75 percent quantile and the 24 percent quantile.

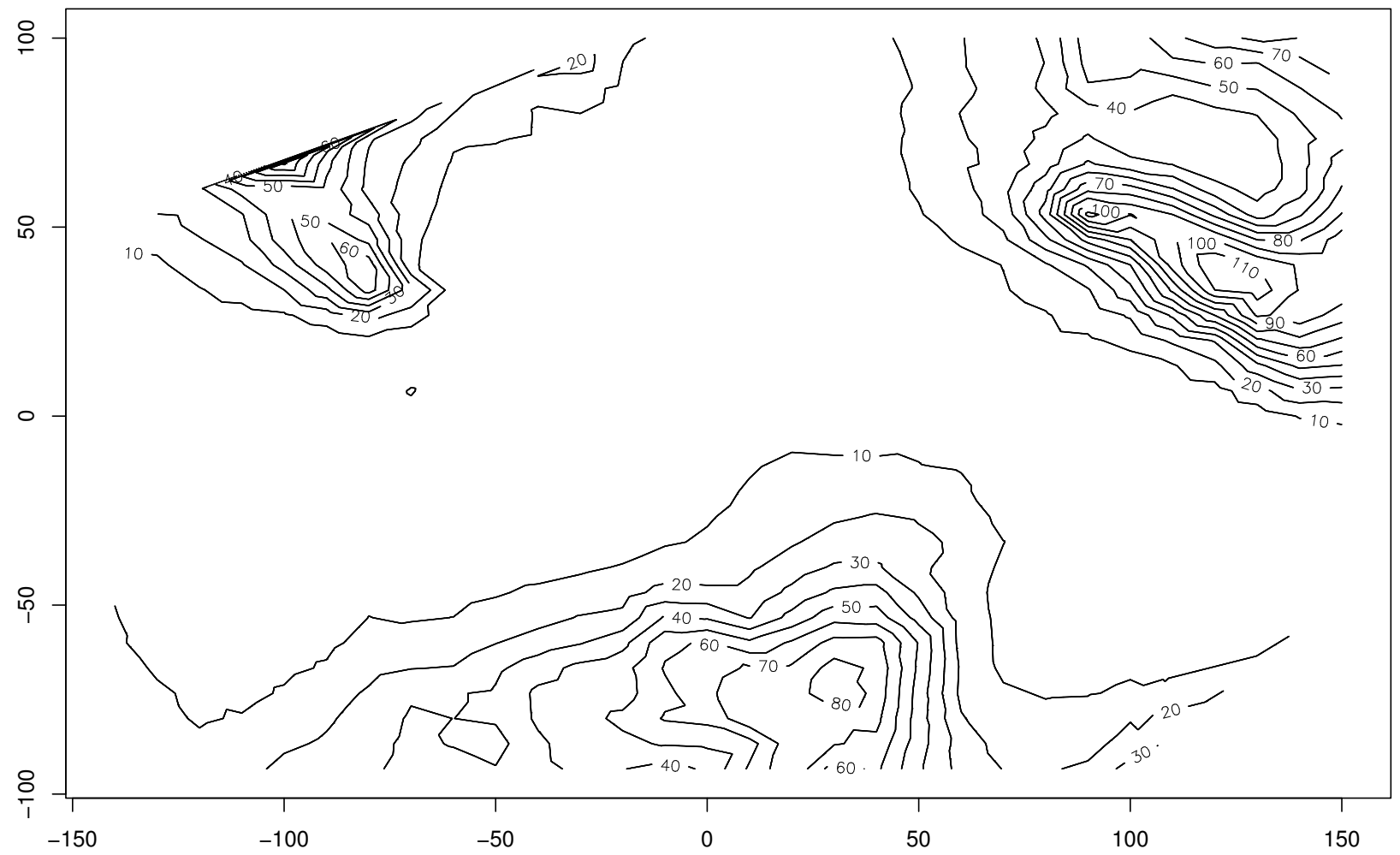

**Fig. 9.21.** Plot of the 95 percent Quantile

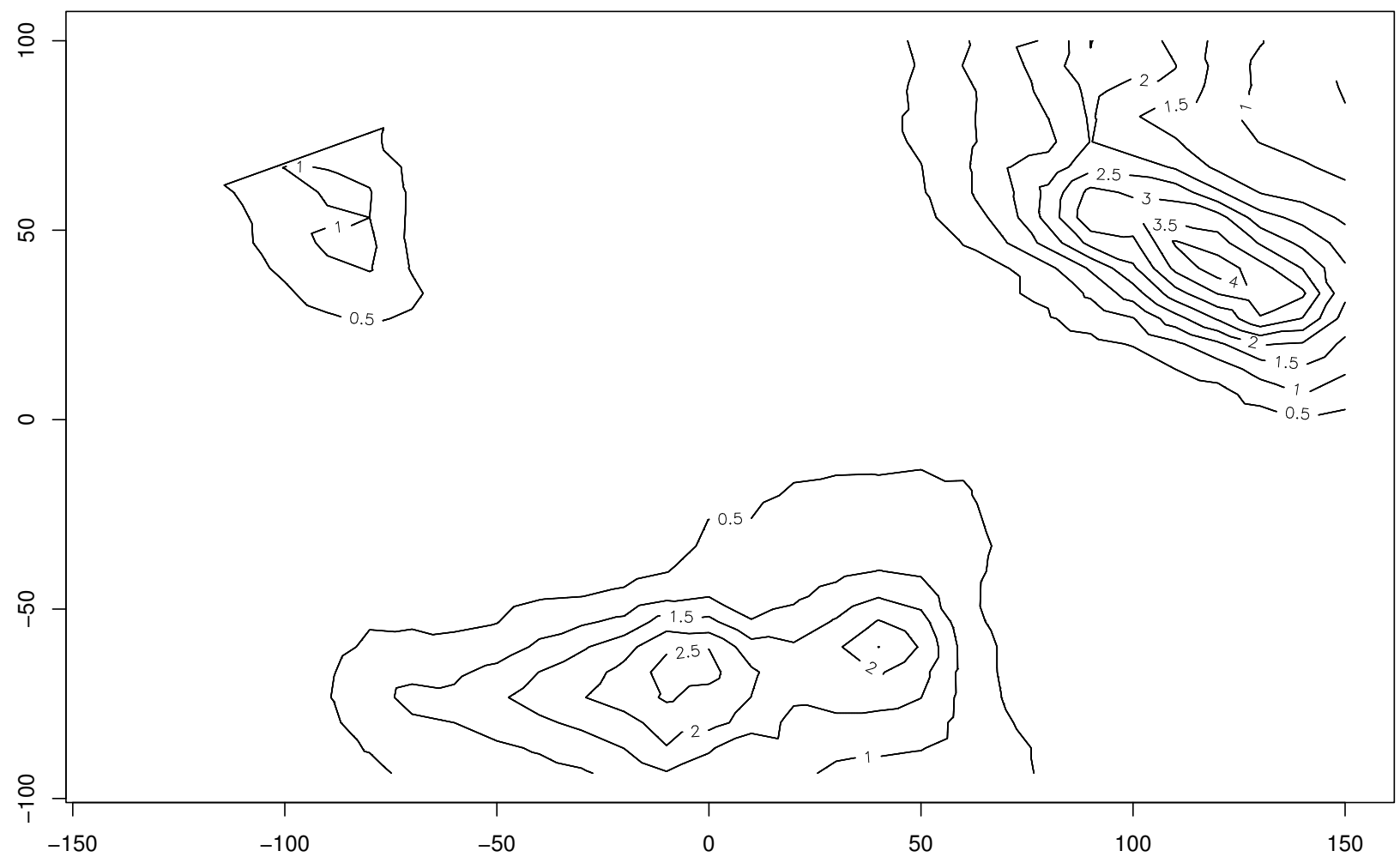

**Fig. 9.22.** Plot of the 5 percent Quantile

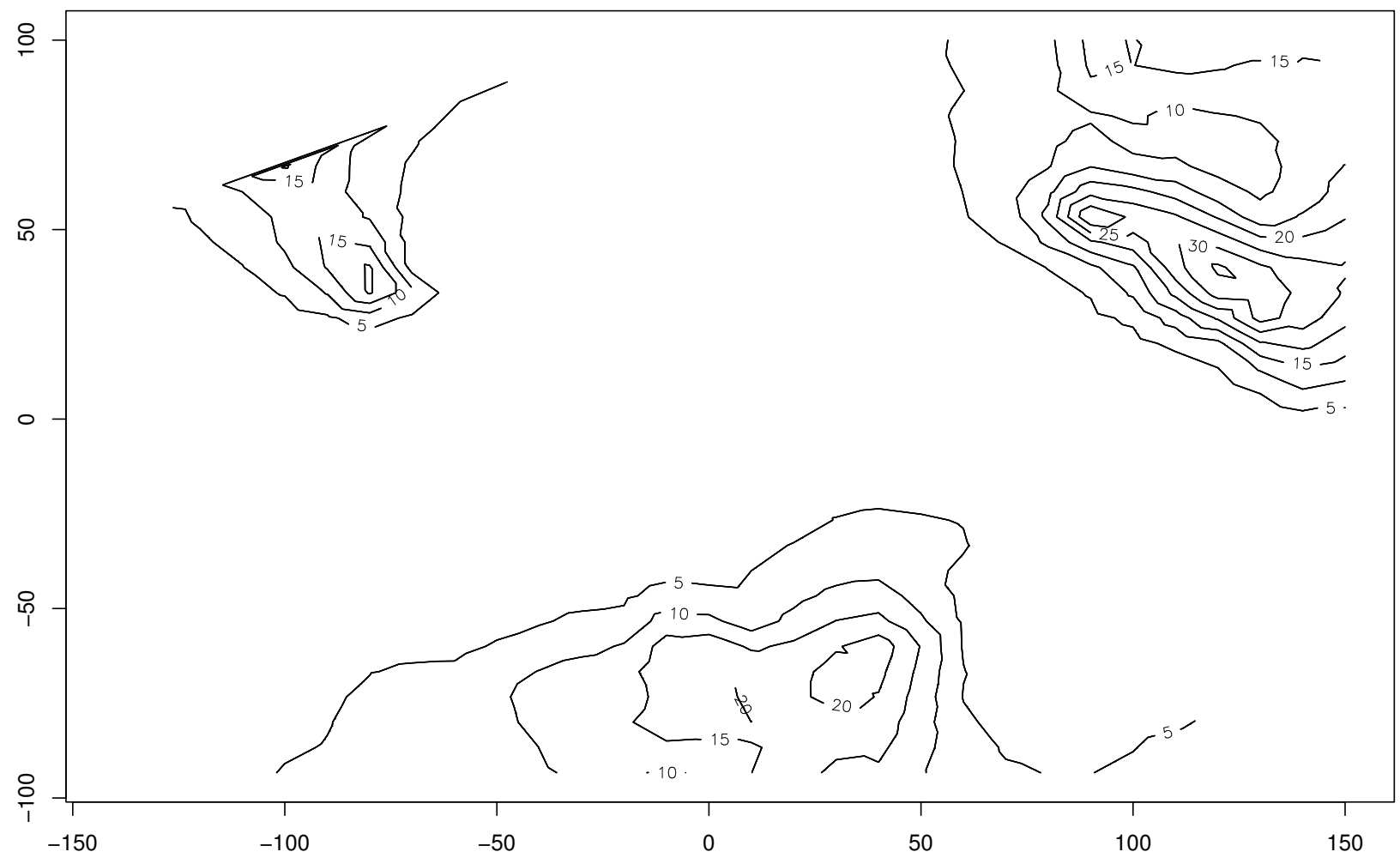

**Fig. 9.23.** Plot of the Inter Quartiles Distance

# 10

# R-Functions

In the following chapter we will give a short introduction to the programming language R and its development through time. The main part here will be the R codes for the implementation of theoretical results obtained in the chapters before.

## 10.1 S and R, the open source software

R is a language and environment for statistical computing, data analysis and graphics, it is also extendable, every user can write new functions and even can rewrite functions that are available, load complied codes from different languages like C, C++, FORTRAN, and the project provides tools for checking the code, writing own extensions and so on. This reduces or removes the barrier between the so called users and the development team and this way it is a perfect GNU project, which is similar to the S language, but in our opinion much better. Historically the environment was developed at Bell Laboratories (formerly AT&T, now Lucent Technologies) by John Chambers and colleagues, but now every R user is a part of R.
R is a language based on modern programming concepts and released under the GNU General Public License, this means that the source code is freely available to everyone and developed for everyone. The so called comprehensive R archive (CRAN) can be found on the R homepage and that is a network of ftp and web servers where documentation and the actual R version can be found. R permits the integration of program scripts with compiled dynamically loaded libraries of functions when computing speed is important. We are able to find lots of packages for spatial statistics, for example at the time, when this thesis is written, there are 335 hits for 'spatial statistics' on the R homepage search engine and also a web page for spatial projects maintained by Roger Bivand. At the moment there are the following spatial statistics packages available for download

- 'DCluster' by Virgilio Gomez-Rubio, Juan Ferrandiz and Antonio Lopez
  A set of functions for the detection of spatial clusters of disease using count data. Bootstrap is used to estimate sampling distributions of statistics.
- 'fields' by Doug Nychka
  Fields is a collection of programs for curve and function fitting with an emphasis on spatial data and spatial statistics. The major methods implemented include cubic, robust, and thin plate splines, universal Kriging and Kriging for large data sets. One main feature is that any covariance function implemented in R can be used for spatial prediction. There are also some useful functions for plotting and working with spatial data as images. This package also contains implementation of a smooth wavelet basis suitable for spatial models.
- 'geoR' by Paulo J. Ribeiro and Peter J. Diggle
  Geostatistical analysis including likelihood-based and Bayesian methods.
- 'geoRglm' by Ole F. Christensen and Paulo J. Ribeiro Jr
  Functions for inference in generalised linear spatial models. The posterior and predictive inference is based on Markov chain Monte Carlo methods.
- 'grasper' by F. Fivaz, A. Lehmann, J.R. Leathwick and J.McC. Overton
  Uses generalised regressions analyses to automate the production of spatial predictions.
- 'rgdal' by Timothy H. Keitt and Roger Bivand
  Provides bindings to Frank Warmerdam's Geospatial Data Abstraction Library (GDAL).
- 'spatialCovariance' by David Clifford
  Functions that compute the spatial covariance matrix for the Matern and power classes of spatial models, for data that arise on rectangular units. This code can also be used for the change of support problem and for spatial data that arise on irregularly shaped regions like countries or zipcodes by laying a fine grid of rectangles and aggregating the integrals in a form of Riemann integration.
- 'spatstat' by Adrian Baddeley and Rolf Turner
  Spatial Point Pattern data analysis, modelling and simulation including multitype/marked points and spatial covariates.
- 'spdep' by Roger Bivand, Luc Anselin, Marilia Carvalho, Stephane Dray, Nicholas Lewin-Koh, Hisaji Ono, Michael Tiefelsdorf and Danlin Yu.
  A collection of functions to create spatial weights matrix objects from polygon contiguities (including the creation of polygon lists from imported shapefiles), from point patterns by distance and tesselations, for summarising these objects, and for permitting their use in spatial data analysis; a collection of tests for spatial autocorrelation, including global Moran's I, Geary's C, Hubert/Mantel general cross product statistic, Empirical Bayes estimates and Assuncao / Reis Index, Getis/Ord G and multicoloured join count statistics, local Moran's I and Getis/Ord G, saddlepoint approximations for global and local Moran's I; and functions for estimating spatial simultaneous autoregressive (SAR) models.

- ’splancs’ by Barry Rowlingson and Peter Diggle
  Spatial and Space-Time Point Pattern Analysis Functions
- ’sgeostat’ by our colleague Albrecht Gebhardt
  An Object-oriented Framework for Geostatistical Modeling in S+
- ’gstat’ by Edzer J. Pebesma
  Variogram modelling; simple, ordinary and universal point or block (co)kriging, sequential Gaussian or indicator (co)simulation; variogram and map plotting utility functions.
- ’GRASS’ by Roger Bivand
  Interface between GRASS 5.0 geographical information system and R.
- ’akima’ by Albrecht Gebhardt
  An R interface to spatial spline interpolation Fortran code by H. Akima. Closely modelled on the interp function in S-PLUS.
- ’tripack’ by Albrecht Gebhardt
  Delaunay triangulation of spatial data. An R interface by Albrecht Gebhardt to Fortran code by R. J. Renka.
- ’spatial’ originally written by B. D. Ripley
  Now part of the VR bundle. Contains trendsurface analysis, kriging and point-process code originally written.
- ’splancs’ by Barry Rowlingson
  Originally commercial code for spatial and space-time point patterns by Barry Rowlingson.
- ’RandomFields’ by Martin Schlather
  Specialized code to simulate from continuous random fields
- ...

Something has to be said about the mother of R, the so called commercial S language. The history of S is relatively long, and as with so many other innovations in software, stems from researchers at Bell Laboratories where the codings for that software have been implemented. The two major sources on the language are Becker, Chambers, and Wilks (1988) and Chambers and Hastie (1992); Venables and Ripley provide a very useful book on an introduction to applied statistics using S (1997). S presents the data analyst with a rich toolbox of components, permitting both the routine processing of statistical tasks, and the programming of new functions not initially included in the language. It employs vectors as basic building blocks (like in R), both permitting the convenient use of linear algebra operations, and the application of standard or user-defined functions to data.
S is now only available commercially as S-PLUS, a fact which has concerned US users, who face a Federal requirement that they ought not to develop software in a language not available from multiple independent sources.
The S-PLUS spatial statistics module includes a fairly wide range of techniques for spatial data analysis, covering many of the key methods presented in Haining (1991) and Cressie (1993). Ripley acted as the main consultant on the design and development of the module. The work of Venables and Ripley on geostatistical and point pattern analyses in their book (1997) is reflected in that toolbox.
In the following the main procedures of the implementation of the so called predictive density is shown and also some coding for the Bayesian ordinary predictor. The codes are written in the way that they should work on R 1.9 and higher, and also in the way that they can be combined with all other spatial packages on CRAN. The user of that code should be able to go through all practical examples in that work and also be able by the documentation of the code to develop it in the sense 'free to all minds'.

## 10.2 Covariance and Variogram Models

```
#################################################################
## Function for evaluating the  exponential  covariance function
## and variogram function
#################################################################

exponentialCov <- function(distance, ae, ce, nugget){
   ifelse(distance>0, 2*nugget + ce*(2- exp(-distance/ae)),
           nugget + ce)
}

exponentialVar <- function(distance, ae, ce, nugget=0){
   ifelse(distance == 0, 0, nugget + ce *
                                (1 - exp(-distance/ae)))
}


#################################################################
## Function for evaluating the Gaussian covariance function
## and vriogram function
#################################################################

gaussianCov <-function(distance, ag, cg, nugget){
   nugget + cg*exp(-distance^2/ag^2)
}

gaussianVar <- function(distance, ag, cg, nugget=0){
   ifelse(distance==0,0, nugget + cg*(1-exp(-distance^2/ag^2)))
}


#################################################################
## Function for evaluating the spherical variogram function
#################################################################
```

```
sphericalVar <- function(distance, as, cs, nugget=0){
   ifelse(distance == 0, 0,
     ifelse(distance <= as, nugget + cs*(3/2* (distance/as) -1/2
           * (distance/as)^3), nugget + cs))
}


################################################################
## Function for evaluation of simple  Matern - covariance function
## and variogram function
################################################################

maternCov <- function(distance, range, sill, nugget, nue)
{

    as <- range/(2 * sqrt(nue))
    u <- distance/as

    ifelse(u>0, nugget+sill/(2^(nue-1)*grangema(nue)) * (u^nue)
                    * besselK(u,nue),nugget + sill)
}

maternVar <- function(distance, range, sill, nugget, nue){
   as <- range/(2 * sqrt(nue))
   u <- distance/as
   ifelse(u > 0, nugget + sill*(1-1/(2^(nue - 1) *
      gamma(nue)) * (u^nue) * besselK(u, nue)),0)
}


##############################################################
### Function for evaluating two Matern variogram functions
##############################################################


MixVariogram <- function(nugget, sill1, range1, nue1, sill2,
                                    range2, nue2, distance,...)
  {
    #######
    # nugget ... value for nugget
    # sill1 ... value for sill1
    # range1 ... value for range1
    # nue1 ... value for nue1
    # sill2 ... value for sill2
    # range2 ... value for range2
    # nue2 ... value for nue2
    #      ... all are parameters form the matern variogram
    #          function
    # distance ... distance in the empirical estimator
```

```
    #######

    k <- sill1+sill2+nugget-(sill1*1/(2^(nue1-1)*gamma(nue1))*
           (distance/(range1/2*(nue1)^(1/2)))^nue1*besselK(
            distance/(range1/2*(nue1)^(1/2)),nue1)+sill2*
            1/(2^(nue2-1)*gamma(nue2))*(distance/(range2/2*
            (nue2)^(1/2)))^nue2*besselK(distance/(range2/2*
            (nue2)^(1/2)),nue2))

      return(k)

  }



################################################################
### Function for calculating the covariance
################################################################

MixCovariance <- function(nugget, sill1, range1, nue1, sill2,
                          range2, nue2, distance)
  {

    ####
    # nugget ... nugget of the sum of two matern covariance
    #            functions
    # sill1/2 ... sills of the sum of two matern covariance
    #             functions
    # range1/2 ... ranges of the sum of two matern covariance
    #              functions
    # nue1/2 ... smoothingparamteters of the sum of two Matern
    #            covariance functions
    # distance ... distance
    ####

    # evaluation of the covariance
    # looking for the Kronecker symbol
    if(abs(distance)>0){

      m <- sill1*1/(2^(nue1-1)*gamma(nue1))*(abs(distance)/
           (range1/2*(nue1)^(1/2)))^(nue1)*besselK(abs(distance)/
           (range1/2*(nue1)^(1/2)),nue1)+sill2*1/(2^(nue2-1)*
           gamma(nue2))*(distance/(range2/2*(nue2)^(1/2)))^(nue2)
           *besselK(distance/(range2/2*(nue2)^(1/2)),nue2)

      # only finite values are given back
      if(is.finite(m))
        {
          return(m)
```

```
        }

      # if not finite we give back "1"
      else
        {
          return(1)
        }
    }

    # looking for the Kronecker symbol
    else
      {
        m<-sill1+sill2+nugget
        return(m)
      }
  }
```

## 10.3 Bayes Kriging

```
################################################################
### Function for calculating the total covariance
################################################################

totalCovarianceMat <- function(x0,y0,x,y,deltaa,aprioriVar)
  {
    ####
    # x0, y0, x, y ... coordinates
    # deltaa ... variogram parameters
    # aprioriVar ... aprioir variance
    ####

    x<-c(x0,x)
    y<-c(y0,y)

    # help variable
    ZeroMatrix <- c(rep(0,length(x)))%*%t(c(rep(0,length(x))))

    Cov <- ZeroMatrix

    # filling the covarinace matrix
    for(i in 1:length(x))
      {
        for(j in 1:length(y))
          {
```

```
          h <- ((x[i]-x[j])^2+(y[i]-y[j])^2)^(1/2)
          # evaluating at i,j
          Cov[i,j] <- MixCovariance(deltaa[1],deltaa[2],
                         deltaa[3],deltaa[4],deltaa[5],
                         deltaa[6],deltaa[7],h)
        }
      }

    # adding the apriori variance
    Cov <- Cov+aprioriVar

    # giving back
    return(Cov)
  }



################################################################
#### Function for Bayes Kriging at a point x0|y0
################################################################

bayesKriging <- function(x0,y0,aprioriMean, aprioriVar,delta0,
                         searchradius,Coord,obsVal)
  {
    #####
    # x0,y0 ... point where the prediction is evaluated
    # aprioriMean ... a priori mean
    # aprioriVar ... a priori variance
    # delta0 ... variogram parameters
    # searchradius ... defining the surrounding
    # Coord ... coordinates in the form x|y
    # obsVal ... observed data
    ####

    # reading out the coordinates
    x <- Coord[,1]
    y <- Coord[,2]

    # help variable
    help1 <- setNames(obsVal,"obsVal")

    # reading the data out
    z <- help1$obsVal

    # parameter we give back
    param <- NULL

    # reading out the variables in the surrounding
    xyz <- searchRad(x0,y0,x,y,z,searchradius)
```

```
    # evaluation of the covariance matrix
    cov <-  totalCovarianceMat(x0,y0,xyz[,1],xyz[,2],delta0,
                               aprioriVar)
    cov1 <- solve(cov)

    # Bayes Predictor
    z1 <- aprioriMean+cov[1,c(2:dim(cov)[2])]%*
              %solve(cov[-1,-1])%*%(xyz[,3]-aprioriMean)
    sqrtmsep <- (1/cov1[1,1])^(1/2)

    param <- matrix(c(z1,sqrtmsep),1,2)
    colnames(param) <- c("prediction","sqrtTMSEP")

    # give it back
    return(param)

  }



#################################################################
#### Function for Bayes Kriging on a grid
#################################################################

BayesKrigingMap <- function(grid,aprioriMean, aprioriVar,delta,
                     searchradius,Coord,obsVal, plotOpt=TRUE,
                     pred=TRUE, tmsep=FALSE, ... )
{

  #####
  # x0,y0 ... point where the prediction is evaluated
  # aprioriMean ... a priori mean
  # aprioriVar ... a priori variance
  # delta ... variogram parameters - Matern
  # searchradius ... defining the surrounding
  # Coord ... coordinates in the form x|y
  # obsVal ... observed data
  # plotOpt ... ask for plotting
  # pred ... plot the Kriging surface
  # tmsep ... plot the tmsep
  # library(akima) is needed
  ####


  # handling the input
  if(!(is.array(grid)|is.data.frame(grid)))
    {
      stop("\n grid must be a array or data.frame\n")
```

```
  }

if(dim(grid)[2]!=2)
  {
    stop("\n grid must be in the form X | Y \n")
  }

if(!(is.array(Coord)|is.data.frame(Coord)))
  {
    stop("\n Coord must be a array or data.frame\n")
  }

if(dim(Coord)[2]!=2)
  {
    stop("\n Coord must be in the form X | Y \n")
  }

if(dim(Coord)[1]!=dim(obsVal)[1])
  {
    stop("\n Number of obsVal and locations in Coord
          differs\n")
  }

if(!(is.finite(searchradius)))
  {
    stop("\n Finite value for searchradius is
          needed\n")
  }

# ask for library(akima)
if(library(akima, logical.return=T)[1]!=TRUE)
  {
    stop("\n \t library akima is needed - please instal
          from CRAN\n\n")
  }

#if it is there - we load it
else
  {
    library(akima)
  }


# help variable
help1 <- setNames(obsVal,"obsVal")

# reading the data out
z <- help1$obsVal
```

```
# value we give back
wert <- NULL

# going throw the grid
for(i in 1:dim(grid)[1])
  {
    # screen output where we are
    cat(i," ")

    #defining the surrounding
    xyz <- searchRad1(grid[i,1],grid[i,2],Coord[,1],
                      Coord[,2],z,searchradius)

    # only take surroundings with mor then 4 points
    if(length(xyz[,1])>4)
      {
        # calling bayesKriging
        wert1 <- matrix(bayesKriging(grid[i,1],grid[i,2],
                        aprioriMean, aprioriVar,delta,
                        searchradius,Coord,obsVal),1,2)

      }

    # less then 4 points - NaN back
    else
      {
        wert1 <- matrix(c(NaN,NaN),1,2)
      }

    # binding the values
    wert <- rbind(wert,wert1)

  }


wert2 <- cbind(wert[,1][(is.finite(wert[,1]))],wert[,2]
               [(is.finite(wert[,2]))])

# splitting the grid
a <- grid[,1]
b <- grid[,2]
grid1 <- a[is.finite(wert[,1])]
grid2 <- b[is.finite(wert[,1])]
grid <- cbind(grid1,grid2)
wert <- wert2

# starting the spline interpolation
x0 <- seq(min(grid[,1]),max(grid[,1]),length=length(grid[,1])
         *3)
```

```
  y0 <- seq(min(grid[,2]),max(grid[,2]),length=length(grid[,2])
           *3)

  interPrediction <- interp(grid[,1],grid[,2],wert[,1],x0,y0)
  intersqrtTMSEP <- interp(grid[,1],grid[,2],wert[,2],x0,y0)


  # ask graphical output
  if(plotOpt)
    {
      if(pred)
        {
          image(interPrediction$x,interPrediction$y,
                interPrediction$z)
          contour(interPrediction$x,interPrediction$y,
                  interPrediction$z,add=T)
        }
      if(tmsep)
        {
          image(intersqrtTMSEP$x,intersqrtTMSEP$y,
                intersqrtTMSEP$z)
          contour(intersqrtTMSEP$x,intersqrtTMSEP$y,
                  intersqrtTMSEP$z,add=T)
        }
    }

  # name the value
  colnames(wert) <- c("prediction","sqrtTMSEP")

  # give it back
  return(wert)
}
```

## 10.4 Simulation

```
################################################################
### Function for simulating Variograms at given points
################################################################


simulate <- function(x,L,point.obj,pairs.obj,n)
{
  #######
  # x ...  a 2d array with the location
```

```
  #       (northing, easting)
  # L ...  Cholesky decomposition of the covariance matrix
  # point.obj ... Given point object generated by
  #               library(sgeostat)
  # pair.obj  ... Given pair object generated by
  #               library(sgeostat)
  # n ... Number of variogarm simulations
  #######

        nd <- nrow(x)
        dimnames(point.obj)[[2]][3] <- "simz"

        cat("make",n," simulations ")
        ####
        # Begin of simulation
        ####
        for (i in 1:n){
          # Number on screen
          cat(i," ")
          z <- rnorm(nd,0,1)
          simz <- as.vector( t(L) %*% as.matrix(z))
          point.obj[,3] <- simz
          simd.var <- my.variogram(point.obj, pairs.obj,'simz')
          # Fist row in simulation output
          if (i==1){
            n.lag <- nrow(simd.var)
            vario <- matrix(0, nrow=n.lag, ncol=(n+3))
            vario[,1] <- simd.var[,1]
            vario[,2] <- simd.var[,2]
            vario[,3] <- simd.var[,4]
          }
          vario[,(i+3)] <- simd.var[,3]
        }
        cat("\n")
        return(invisible(vario))
}



################################################################
### Function for plotting the simulated variograms
################################################################

simplot <- function(vdata,col="black", main="variogram
                    simulation")
{
  #####
  # vdata ... empirical variograms form the simulation
  #           typically the output from simulate.R
```

```
  # col   ... colour of the plot
  # main  ... title of the plot
  #####

  col1 <- col
  main1 <- main

  n <- ncol(vdata)
  vmax <-  max(vdata[,-c(1:3)]/2)

  # stating the printing on screen
  # performed by lines

  plot(vdata[,2],vdata[,4]/2,ylim=c(0,vmax),type="l",main=main1,
       col=col1)
    for(i in 5:n){
        lines(vdata[,2],vdata[,i]/2,type="l")
    }
}
```

## 10.5 Predictive Density

```
#################################################################
### Function for defining the neighbourhood of the prediction
#################################################################

searchRad <- function(x0,y0,x,y,z,searchradius,...)
{
   ######
  # x0 ... x coordinate of the centre
  # y0 ... y coordinate of the centre
  # x,y,z ... coordinates and values in the surrounding
  # searchradius ... radius of the surrounding
  ######

  # help variable
  Rad <- searchradius

  # evaluating the points
  dist <- as.vector(((x0-x)^2+(y0-y)^2)^(1/2))

  # reading them out
  zz <- z[dist<=Rad]
  xx <- x[dist<=Rad]
  yy <- y[dist<=Rad]
```

```
  # value to give back
  xyz <- cbind(xx,yy,zz)

  # giving back
  return(xyz)
}



################################################################
### Function for fitting two Matern variogram functions to
### empirical simulated variograms
################################################################

fitTwoMatern <- function(vdata, point.obj, pair.obj, a1,
                         ordinary=TRUE, startvalues=NULL,
                         plotoption=TRUE,...)
{

  ######
  # vdata ... simulated empirical variograms, output varsim
  #           col1 <- Lag Number
  #           col2 <- distance
  #           col3 <- n(h)
  #           then the simulated
  # point.obj ... Given point object generated by
  #               library(sgeostat)
  # pair.obj  ... Given pair object generated by
  #               library(sgeostat)
  # a1 ... observed data in the random field
  # ordinary ... ask for OLS (true) or WLS (false)
  # startvalues ... optional the startingvalues for the
  #                 fitting
  # plotoption ... TRUE - plot of the fitting is shown
  #                FALSE - no output on screen
  ######

  # runindex
  nv <- ncol(vdata)
  nl <- nrow(vdata)

  # what we will give back
  parameter <- NULL

  # working variables
  point <- point.obj
  pair <- pair.obj
  a1 <- a1
```

```
#screen output OLS or WLS
if(ordinary)
  {
    cat("\n \t Performing OLS \n")
  }
else
  {
    cat("\n \t Performing WLS \n")
  }
#check input
if((class(point)[1]!="point"))
  {
    stop("Wrong input, \t point object is needed \n")
  }
if(class(pair)!="pair")
  {
    stop("Wrong input, \t pair object is needed \n")
  }

# try-catch - check distance
distance <- try(as.vector(vdata[,2]))
if(class(distance[1])=="try-error")
  {
    cat("\n Wrong Input in \n \t fitTwoMatern(***,...) \n")
    stop("New input on that position - error, stop evaluating
          \n")
  }

# runindex
m <- length(distance)

# check input
if(m==0)
  {
    stop("The distance has length 0 - error, stop evaluating
          \n")
  }
if(is.null(m))
  {
    stop("NULL object in distance! error, stop evaluating,
          \n \t Wrong input \n")
  }

# help variable for lag - use in WLS
runlag <- vdata[,3]

#new Col names - able now
```

```
a1 <- setNames(a1, "observ")
point <- setNames(point,c("x","y","observ"))

# output on screen
cat("\n")

# running up all simulated empirical variograms
for(i in 4:nv)
  {
    param2 <- NULL
    vdatai <- try(as.vector(vdata[,i]))

    #check vdatai
    if(class(vdatai[1])=="try-error")
      {
        stop("\n Error in Input, simulated variograms are
              false \n")
      }

    # output on screen
    cat((i-3),",")

    # starting values
    # if is.null - calculate them

    if(length(startvalues)!=7){
      # starting values
      # nugget - linear interpolation to the origin
      nuggetStart <- vdatai[1]-(vdatai[1]-vdatai[2])/
                     (vdata[,2][1]-vdata[,2][2])*vdata[,2][1]
      if(nuggetStart<0)
        {
          nuggetStart <- 0.01
        }
      #sill1 - mean of vdatai
      sill1Start <- var(vdatai)

      #range1 - where sill1 is at x
      vdata2 <- vdata[,2]
      range1Start <- vdata2[length(vdatai[vdatai<sill1Start])]

      #sill2 - q075
      sill2Start <- quantile(vdatai)[4][[1]]

      #range2 - where sill2 is at x
      range2Start <- vdata2[length(vdatai[vdatai<sill2Start])]

      # nue1 and nue2 are assumed to be 0.5
      startvalue <- c(nuggetStart, sill1Start, range1Start,
```

```
                 0.5, sill2Start, range2Start, 0.5)
}
# startingvalues are an input vector

else
  {
            startvalue <- startvalues
  }
# performing the fit(param,vdatai,distance,ordinary,runLag)
h1<-try(fitVal(param=startvalue, vdatai,distance,ordinary,
        runlag))

if(class(h1[1])=="try-error")
  {
    stop("\n Error in function fitVal \n contact:
          philipp.pluch@uni-klu.ac.at")
  }
# catch if error
# performed well than
if(class(h1[1])!="try-error")
  {
    #where we calculate the variogram
    x<-seq(0.00001,distance[length(distance)],length=
           length(distance))

    # evaluating the covariogram
    CO2<-h1[1]+h1[2]+h1[5]
    CH2<-h1[2]*1/(2^(h1[4]-1)*gamma(h1[4]))*
         (x/(h1[3]/2*(h1[4])^(1/2)))^h1[4]*
         besselK(x/(h1[3]/2*(h1[4])^(1/2)),
         h1[4])+h1[5]*1/(2^(h1[7]-1)*gamma(h1[7]))*
         (x/(h1[6]/2*(h1[7])^(1/2)))^h1[7]*
         besselK(x/(h1[6]/2*(h1[7])^(1/2)),h1[7])
    # evaluating the variogram
    var2<-CO2-CH2

    # if the variogram is finite
    if(is.finite(var2[1]))
      {
        # ask for plotting
        if(plotoption)
          {
            plot(x,var2,col="red",type="l",xlim=c(0,200),
                 ylim=c(0,20))
            lines(vdata[,2],vdata[,i])
          }
      }

    # give this back
```

```
          param2<-matrix(h1,1,7)

          # name the rows
          rownames(param2)<-paste("v",(i-3),sep="")

          #combine them
          parameter<-rbind(parameter,param2)
        }
      # error in fitting
      else
        {
          # bad thing is given back
          h1<-c(0,0,0,0,0,0,0)

          param2<-matrix(h1,1,7)

          # name them and give them back
          rownames(param2)<-paste("v",(i-3),sep="")
          parameter<-rbind(parameter,param2)
        }

    } # end of for

  # name the colums
  colnames(parameter)<-c("nugget","sill1","range1","nue1",
                       "sill2","range2","nue2"))

  cat("\n")

  return(parameter)

}




################################################################
### Function for fitting
################################################################

fitVal<-function(param,vdatai,distance,ordinary,runLag,...)
  {
    #########
    # param ... parameters for starting from Matern variogram
    # vdatai ...  simulated variograms
    # distance ... distance from lag class
    #########

    # fitting is started
```

```
# L-BFGS-B method is used for bounds
# relative tolerance is as high as possible

# warnings() are given out
g<- try(optim(par=param,fn=leastsquares,vdatai=vdatai,
        method="L-BFGS-B",lower=c(0.001,0.001,0.001,
        0.001,0.001,0.001,0.001),upper=c(10,18,400,
        40,18,400,40),gr=NULL, control = list(maxit=
        40000, reltol=0.0000000000000001), distance=
        distance, ordinary=ordinary, runLag=runLag))

# catch the error
if(class(g)!="try-error")
  {
    back <- as.vector(g$par)
    return(back)
   # return(cbind(g$par,g$value,g$convergence))
  }
# if error - give that back
else
  {
    g<-c(0.001,0.001,0.001,0.001,0.001,0.001,0.001)
    cat("fitting error")
    return(g)
  }

}




################################################################
### Least squares
################################################################
leastsquares <-function(param, vdatai, distance, ordinary,
                        runLag,...)
    {
      ########
      # param, vdatai, distance ... values for computation
      # ordinary ... ask for WLS (false) or OLS (true)
      # runLag ... number of points in each lag
      ########

      nugget<-param[1]
      sill1<-param[2]
      range1<-param[3]
      nue1<-param[4]
      sill2<-param[5]
      range2<-param[6]
```

```
      nue2<-param[7]

      # return value
      z <- NULL

      gamma <- MixVariogram(nugget, sill1, range1, nue1, sill2,
                            range2, nue2, distance)

      #performing ordinary least squares
      if(ordinary)
        {
          z <-sum((gamma-vdatai)^2)
        }
      #performing weighted least squares
      else
        {
          Var <- 2*((gamma)^2)/runLag
          z <-sum((gamma-vdatai)^2/Var)

        }

      # only a finite value is given back
      if(is.finite(z))
        {
          return(z)
        }
      else
        {
          z<-0.001
          return(z)
        }
    }


################################################################
### Function for evaluating the predictive density at
### a given point
################################################################


ourBayes <- function(x0,y0,aprioriMean, aprioriVar,delta,
                     searchradius, Coord, obsVal,
                     plotOptPred=FALSE, ...)
  {

    ####
    # x0,y0 ... point where the predictive density is evaluated
    # aprioriMean ... a priori mean
```

```
# aprioriVar ... a priori variance
# delta ... simulated fitted variogarm parameters
# searchradius ... parameter for defining the surrounding
# Coord ... given coordinates where the observations are
#           done in the form: x|y - array
# obsVal ... the observed values
# plotOptPred ... ask for plotting the predictive density
####

# reading out the coordinates
x <- Coord[,1]
y <- Coord[,2]

# look at the str(obsVal)
if(!(is.vector(obsVal)))
   {
     # help variable
     z1 <- setNames(obsVal, "observ")

     # reading out the column
     z <- z1$observ
   }
else
  {
    z <- obsVal
  }

# selecting the points and there values in the surrounding
xyz <- searchRad(x0,y0,x,y,z,searchradius)

# initial values for the distribution
predictiveDist <- 0

# initial values for the mean
predictiveMean <- 0

# initial values for the variance
predictiveVar <- 0

# runindex

k <- 0

# evaluation of the predictive density
for(i in 1:dim(delta)[1])
  {

        # run index for averaging
        k <- k+1
```

```
        #  print(k)
        # evaluating the covariance
        cov <- totalCovarianceMat(x0,y0,xyz[,1],xyz[,2],
                                  delta[i,],aprioriVar)

        # evaluating the inverse of cov
        cov1 <- solve(cov)

        # evaluating the a posteriori variance and mean
        aposterioriVar <- 1/cov1[1,1]
        # print(aposterioriVar)
        aposterioriMean <- aprioriMean+cov[1,
                           c(2:dim(cov)[2])]%*%solve
                           (cov[-1,-1])%*%(xyz[,3]-
                           aprioriMean)

        # distribution is given by
        distri <- dnorm(seq(-5,6.5,0.01),aposterioriMean,
                        (aposterioriVar)^(1/2))

        # predictive dist is given by
        predictiveDist <- predictiveDist+distri*(1/exp
                                        (seq(-5,6.5,0.01)))

  }


# simple averaging
predictiveDist <- predictiveDist/k

# cumulative predictive distribution
cumPredictiveDist <-cumsum(predictiveDist*(exp(seq(-4.99,
                            6.51, 0.01))-exp(seq(-5,6.5
                            0.01))))/sum(predictiveDist*
                            (exp(seq(-4.99,6.51,0.01))-
                            exp(seq(-5,6.5,0.01))))

# reading out the index for
# q005
ind1 <- length(cumPredictiveDist[cumPredictiveDist<=0.05])
#q095
ind2 <- length(cumPredictiveDist[cumPredictiveDist<=0.95])
# q001
ind3 <- length(cumPredictiveDist[cumPredictiveDist<=0.01])
# q025
ind4 <- length(cumPredictiveDist[cumPredictiveDist<=0.25])
# q075
ind5 <- length(cumPredictiveDist[cumPredictiveDist<=0.75])
# q099
```

```
ind6 <- length(cumPredictiveDist[cumPredictiveDist<=0.99])
# median
indMedian <- length(cumPredictiveDist[cumPredictiveDist <=
                    0.5])

# modal
modalWert <- max(predictiveDist)
indModal <- NULL
for(i in 1:length(predictiveDist))
  {
    if(predictiveDist[i]==modalWert)
      {
        indModal <- i
      }
  }

# for transformation
x <- exp(seq(-5,6.5,0.01))
Median <- x[indMedian]

Modal <- x[indModal]
qq005 <- x[ind1]
qq095 <- x[ind2]
qq001 <- x[ind3]
qq025 <- x[ind4]
qq075 <- x[ind5]
qq099 <- x[ind6]

Mean <- sum(predictiveDist*(exp(seq(-4.99,6.51,0.01))-
            exp(seq(-5,6.5,0.01)))*(exp(seq(-4.99,6.51,
            0.01))+exp(seq(-5,6.5,0.01)))/2)/
            sum(predictiveDist*(exp(seq(-4.99,6.51,0.01))
            -exp(seq(-5,6.5,0.01))))

# ask for plotting
if(plotOptPred)
  {
    plot(exp(seq(-5,6.5,0.01)),predictiveDist,type="l",
         xlim=c(0,qq095+10),xlab="",ylab="")
    # plot(predictiveDist)
    points(Median,0,col="red")
    points(Modal,0,col="blue")
    points(Mean,0,col="green")
    points(qq005,0,col="orange")
    points(qq095,0,col="orange")

    # adding a legend to the plot
    legend(qq095+5,0.75*max(predictiveDist),
     c("median","modal","mean","q005","q095"), pch=21,
```

```
        pt.bg="white", lty=0, col = c("red","blue","green",
        "orange","orange"))
      }

    predicted <- matrix(c(Modal,Median,Mean,qq001,qq005,qq025,
                        qq075,qq095,qq099),1,9)
    colnames(predicted) <- c("Modal","Median","Mean","qq001",
                             "qq005",,"qq025","qq075","qq095",
                             "qq099")

    # value to give back
    return(predicted)
  }




################################################################
### Function for making a map of the predictive density
################################################################

BayesMap <- function(grid,aprioriMean, aprioriVar,delta,
                     searchradius,Coord,obsVal, plotOpt=TRUE,
                     meanP=TRUE, medianP =FALSE, modalP=FALSE,
                     q005P=FALSE, q025P=FALSE, q075P=FALSE,
                     q095P=FALSE, approxVarP=FALSE,...)
{

  #######
  # grid ... grid where the predictive density is evaluated
  # aprioriMean ... apriori mean
  # aprioriVar ... apriori variance
  # delta ... fitted simulated variogram parameters
  # searchradius ... surounding for the prediction
  # Coord ... coordinates for the observed values in the form
  #           x|y
  # obsVal ... observed values
  # plotOpt ... ask for graphical output
  # meanP ... ask for graphical output of mean map
  # medianP ... ask for graphical output of median map
  # modalP ... ask for graphical output of modal map
  # q005P ... ask for graphical output of q005 map
  # q025P ... ask for graphical output of q025 map
  # q075P ... ask for graphical output of q075 map
  # q095P ... ask for graphical output of q095 map
  # approxVarP ... ask for the graphical output of apprx.
  #                variance (q075-q025)/1.45
  #
  # library(akima) is needed!
```

```
#######

# handling the input
if(!(is.array(grid)|is.data.frame(grid)))
  {
    stop("\n grid must be a array or data.frame\n")
  }

if(dim(grid)[2]!=2)
  {
    stop("\n grid must be in the form X | Y \n")
  }

if(!(is.array(Coord)|is.data.frame(Coord)))
  {
    stop("\n Coord must be a array or data.frame\n")
  }

if(dim(Coord)[2]!=2)
  {
    stop("\n Coord must be in the form X | Y \n")
  }

if(dim(Coord)[1]!=dim(obsVal)[1])
  {
    stop("\n Number of obsVal and locations in Coord
          differs\n")
  }

if(!(is.finite(searchradius)))
  {
    stop("\n Finite value for searchradius is needed\n")
  }

# ask for library(akima)
if(library(akima, logical.return=T)[1]!=TRUE)
  {
    stop("\n \t library akima is needed - please
          instal from CRAN\n\n")
  }

#if it is there - we load it
else
  {
    library(akima)
  }

# help variable
help1 <- setNames(obsVal,"obsVal")
```

```
# reading out and making a vector
obsVal <- as.vector(help1$obsVal)
# return value
wert <- NULL

# going throw the grid
for(i in 1:dim(grid)[1])
  {
    # output on screen
    cat(i," ")

    # selecting the points in the surrounding
    xyz <- searchRad(grid[i,1],grid[i,2],Coord[,1],Coord[,2],
                     obsVal,searchradius)

    # selecting more then 4 points - 'couse of inverse
    if(length(xyz[,1])>4)
      {

        # calling the function ourBayes

        wert1 <- matrix(ourBayes(grid[i,1], grid[i,2],
                 aprioriMean, aprioriVar, delta, searchradius,
                 Coord, obsVal),1,9)

      }

    # if there are less then 4 points we give NaN back
    else
      {
        wert1 <- matrix(c(NaN,NaN,NaN,NaN,NaN,NaN,NaN,NaN,
                        NaN),1,9)
      }

    # binding the values
    wert <- rbind(wert,wert1)

  }

# taking only the finite values
wert2 <- cbind(wert[,1][(is.finite(wert[,1]))],wert[,2]
               [(is.finite(wert[,2]))],wert[,3]
               [(is.finite(wert[,3]))],wert[,4]
               [(is.finite(wert[,4]))],wert[,5]
               [(is.finite(wert[,5]))],wert[,6]
               [(is.finite(wert[,6]))],wert[,7]
               [(is.finite(wert[,7]))],wert[,8]
               [(is.finite(wert[,8]))],wert[,9]
```

```
                [(is.finite(wert[,9]))])

# splitting up the grid
a <- grid[,1]
b <- grid[,2]

# just take those points with finite values on the grid
grid1 <- a[is.finite(wert[,1])]
grid2 <- b[is.finite(wert[,1])]

# building up the new grid
grid <- cbind(grid1,grid2)

# helping variable
wert <- wert2

# sequence for interpolation
x0 <- seq(min(grid[,1]),max(grid[,1]),length=length(grid[,1])
         *3)
y0 <- seq(min(grid[,2]),max(grid[,2]),length=length(grid[,2])
         *3)

# the interpolated values
interModal <- interp(grid[,1],grid[,2],wert[,1],x0,y0)
interMedian <- interp(grid[,1],grid[,2],wert[,2],x0,y0)
interMean <- interp(grid[,1],grid[,2],wert[,3],x0,y0)
inter005 <- interp(grid[,1],grid[,2],wert[,5],x0,y0)
inter025 <- interp(grid[,1],grid[,2],wert[,6],x0,y0)
inter075 <- interp(grid[,1],grid[,2],wert[,7],x0,y0)
inter095 <- interp(grid[,1],grid[,2],wert[,8],x0,y0)
interapproxVar <- interp(grid[,1],grid[,2],
                        (wert[,7]-wert[,6])/1.45)

# save approximative variance
approxVar <- (wert[,7]-wert[,6])/1.45

# ask for graphical output
if(plotOpt)
  {
    plot(Coord)

    if(modalP)
      {
        contour(interModal$x,interModal$y,interModal$z,
                col="green",add=T)
      }

    if(medianP)
      {
```

```
        contour(interMedian$x,interMedian$y,interMedian$z,
                add=T)
      }

    if(meanP)
      {
        contour(interMean$x,interMean$y,interMean$z,
                col="violet",add=T)
      }

    if(q005P)
      {
        contour(inter005$x,inter005$y,inter005$z,
                col="yellow",add=T)
      }

    if(q025P)
      {
        contour(inter025$x,inter025$y,inter025$z,
                col="orange",add=T)
      }

    if(q075P)
      {
        contour(inter075$x,inter075$y,inter075$z,
                col="red",add=T)
      }

    if(q095P)
      {
        contour(inter095$x,inter095$y,inter095$z,
                col="blue",add=T)
      }

    if(approxVarP)
      {
        contour(interapproxVar$x, interapproxVar$y,
                interapproxVar$z, col="ligthblue",add=T)

      }

    points(grid,col="red")
  }

wert <- cbind(wert,approxVar)
# name the columns

colnames(wert) <- c("Modal","Median","Mean","qq001","qq005",
                    "qq025","qq075","qq095","qq099",
```

```
                          "approxVar")

  # give the value back
  return(wert)
}
```

# List of Figures